\documentclass[12pt]{amsart}
\usepackage{amsmath,amsthm,amsfonts,amssymb,amsxtra,mathrsfs,txfonts,amscd,color,esint,mathtools}
\usepackage[colorlinks,citecolor=blue,pagebackref,hypertexnames=false]{hyperref}
\usepackage{cite}
\usepackage{needspace}
\usepackage[T1]{fontenc}
\allowdisplaybreaks
\newtheorem{theorem}{Theorem}[section]
\newtheorem{proposition}[theorem]{Proposition}
\newtheorem{lemma}[theorem]{Lemma}
\newtheorem{corollary}[theorem]{Corollary}
\newtheorem{remark}[theorem]{Remark}
\theoremstyle{plain}
\newtheorem{definition}[theorem]{Definition}
\theoremstyle{plain}

\numberwithin{equation}{section}

\newcommand{\C}{\mathbb C}
\newcommand{\F}[1]{F^{#1}(\varphi)}
\newcommand{\Lw}[1]{L^{#1}(\varphi)}
\newcommand{\rank}{\operatorname{rank}}
\newcommand{\As}[1]{\|#1\|_{\mathfrak A_s}}

\title[Approximation numbers of Hankel Operators]{Approximation numbers of Hankel Operators on Weighted Fock Spaces
	and the Berger--Coburn Phenomenon}
	
	\author{Guijun Liu,  Xiaofeng Wang and Jin Xia }
\address{Guijun Liu, Xiaofeng Wang and Jin Xia,
	School of Mathematics and Information Science,
	Guangzhou University, Guangzhou 510006, China}

\email{lgj19981998@163.com (G. Liu);
	wxf@gzhu.edu.cn (X. Wang);
2695931921@qq.com (J. Xia)}
\date{\today}

\begin{document}
\begin{abstract}
Let $\varphi\in C^2(\C^n)$ be a real-valued weight whose
real Hessian is uniformly positive definite and bounded.
For $0<p<\infty$, $1\le q<\infty$, and possibly unbounded
symbols $f$, we give necessary and sufficient conditions
for prescribed decay of the approximation numbers of
Hankel operators $H_f$ from weighted Fock spaces $F^p(\varphi)$
to weighted Lebesgue spaces $L^q(\varphi)$ in terms of
the non-increasing rearrangement of the IDA function of $f$.
We also characterize membership in the approximation ideals
$\mathfrak A_s(F^p(\varphi),L^q(\varphi))$ and their weak
counterparts $\mathfrak A_{s,\infty}(F^p(\varphi),L^q(\varphi))$
for $0<s<\infty$. The symbol criteria involve Lorentz,
Lorentz--Zygmund, and Copson--Lorentz spaces.
Our approach also yields the converse to the sufficient condition
in \cite{FWZ26}, thereby completing the characterization of
prescribed singular value decay for Hankel operators on a broad
class of weighted Bergman spaces.

For the pair $H_f$ and $H_{\bar f}$ on weighted Fock spaces,
we characterize simultaneous decay of their approximation
numbers and simultaneous membership in these ideals in terms
of the mean oscillation function of $f$.
For bounded symbols, we establish the Berger--Coburn
phenomenon for both ideals in sharp exponent ranges.
A new feature in the case $p<2<q$ is that the Berger--Coburn
equivalence holds for $s_0(p,q)<s\le1$, where the sharp
threshold $s_0(p,q)$ is strictly less than $1$.
\end{abstract}

\baselineskip 16.6pt
\hfuzz=6pt
\widowpenalty=10000
\maketitle
\setcounter{tocdepth}{1}
\tableofcontents
\medskip

\section{Introduction}
\subsection{Background and problem}

Hankel operators form an important class of operators on spaces
of analytic functions and have close connections with function
theory, harmonic analysis and approximation theory.
Peller's monograph \cite{Peller03} provides a comprehensive account
of their classical theory and applications.
Among the central questions concerning these operators are
boundedness, compactness and Schatten class membership; see, for example,
\cite{AFJP91,AFP88,Axler86,Bauer05,BEMN22,Dostanic04,ER09,
	FWZ26,HV22,HV23,HW18,Isralowitz11,Isralowitz13,Janson88,
	Luecking92,Pau16,SY13,Wallsten90,Xia01,Xia02,Xia23,
	XZ04,ZHW25,LWH25,LWZ26,Zhu90,Zhu91,Zhu07}
and the references therein.

For possibly unbounded symbols, Luecking \cite{Luecking92}
obtained a complete characterization of these three properties on the standard Bergman space in 1992. The corresponding problems on Fock spaces remained
open for more than two decades, until Hu and Virtanen \cite{HV22,HV23} resolved them using IDA spaces. We recall their characterization in the weighted Fock
setting considered here. For a domain $\Omega\subset\C^n$, let $\mathcal H(\Omega)$
denote the space of holomorphic functions on $\Omega$. Denote by $dv$ the Lebesgue measure on $\C^n$. For $0<p<\infty$, the weighted Fock space is defined by
\[
F^p(\varphi):=
\left\{
f\in\mathcal H(\C^n):
\|f\|_{p,\varphi}:=
\left(
\int_{\C^n}|f(z)|^p e^{-p\varphi(z)}\,dv(z)
\right)^{1/p}<\infty
\right\}.
\]
where $\varphi$ is a weight function satisfying
the assumptions specified in Section~2.  In particular, $\varphi(z)=\alpha|z|^2/2$, $\alpha>0$,
gives the classical Fock spaces $F_\alpha^p$;
see \cite{JPR87,Zhu12}. For $0<p<\infty$ and $f\in\mathcal S$
(see \eqref{eq:symbol1} for the precise definition),
the Hankel operator $H_f$ is densely defined on
$F^p(\varphi)$ by
\[
H_fg:=(I-P)(fg),
\]
where $P$ is the orthogonal projection from
$L^2(\varphi)$ onto $F^2(\varphi)$. Let $L^q(\varphi):=L^q(\mathbb C^n, e^{-q\varphi}dv)$.
The precise result is stated as follows.
\begin{enumerate}
	\item[\textup{(A)}]
	If $0<p\le q<\infty$ and $q\ge1$, then
	$H_f$ is bounded \textup{(compact)} from
	$F^p(\varphi)$ to $L^q(\varphi)$ if and only if
	$f\in\mathrm{BDA}^q$ \textup{($f\in\mathrm{VDA}^q$)}.
	
	\item[\textup{(B)}]
	If $1\le q<p<\infty$, then
	$H_f$ is bounded if and only if it is compact
	from $F^p(\varphi)$ to $L^q(\varphi)$,
	and this holds if and only if
	$f\in\mathrm{IDA}^{r,q}$,
	where $1/r=1/q-1/p$.
	
	\item[\textup{(C)}]
	For $0<s<\infty$, $H_f\in\mathcal S_s(F^2(\varphi), L^2(\varphi))$
	if and only if $f\in\mathrm{IDA}^{s,2}$.
\end{enumerate}
where the spaces $\mathrm{BDA}^q$, $\mathrm{VDA}^q$
and $\mathrm{IDA}^{s,q}$ are defined below.

Along this line, Fan, Wang, and Zeng \cite[Theorem~1.4]{FWZ26} studied
the singular value decay of Hankel operators on weighted Bergman
spaces $A_\omega^2$ with $\omega\in\mathcal W^*(\Omega)$
(see Subsection~\ref{sec:Bergman-converse} for the precise
definitions). For $0<\delta<\delta_\omega$ and
$f\in\mathcal S_\omega$ as in \eqref{BergS}, suppose that
$H_f:A_\omega^2\to L_\omega^2$ is compact and that
$\Psi:[1,\infty)\to(0,\infty)$ is nondecreasing. They proved
the following assertions.
\begin{enumerate}
\item[\textup{(D)}]
If $\Psi(t)/t^\gamma$ is nonincreasing for some $0<\gamma<1$, then
\begin{equation*}
s_k(H_f)\lesssim\Psi(k)^{-1}\quad(k\in\mathbb N)
\quad\Longleftrightarrow\quad
(G_{\delta,\omega}(f))^*(k)\lesssim\Psi(k)^{-1}\quad(k\in\mathbb N).
\end{equation*}

\item[\textup{(E)}]
If $\Psi(t)/t^\gamma$ is nonincreasing for some $\gamma>0$, then
\begin{equation*}
(G_{\delta,\omega}(f))^*(k)\lesssim\Psi(k)^{-1}\quad(k\in\mathbb N)
\quad\Longrightarrow\quad
s_k(H_f)\lesssim\Psi(k)^{-1}\quad(k\in\mathbb N),
\end{equation*}
\end{enumerate}
where $s_k(H_f)$ denotes the $k$th singular value of $H_f$,
and $(G_{\delta,\omega}(f))^*$ denotes the nonincreasing rearrangement
with respect to $d\lambda_\omega$ of the IDA function in \eqref{BergIDA}.

Beyond the Hilbert space setting, approximation numbers extend
the notion of singular values. Introduced by Pietsch \cite{Pietsch63},
they measure the error of approximation by operators of prescribed
finite rank and form part of the theory of operator ideals
\cite{Pietsch74,Pietsch80}. For $1\le p<\infty$, their decay has
been studied for composition operators
$C_\psi:H^p(\mathbb D)\to H^p(\mathbb D)$ on the Hardy spaces
of the unit disk \cite{LQR15} and
$C_\psi:\mathcal H^p\to\mathcal H^p$ on the Hardy spaces
of Dirichlet series \cite{BQS16}.
We first recall the relevant definitions.

\begin{definition}
	Let $X$ and $Y$ be quasi-Banach spaces, and let
	$\mathcal L(X,Y)$ denote the bounded linear operators
	from $X$ to $Y$.
	For $T\in\mathcal L(X,Y)$, its $k$th approximation
	number is
	\[
	a_k(T):=
	\inf\bigl\{\|T-R\|_{X\to Y}:
	R\in\mathcal L(X,Y),\ \rank R<k\bigr\},
	\qquad k\ge1.
	\]
	For $0<s<\infty$, the approximation ideal
	$\mathfrak A_s(X,Y)$ consists of those $T$ satisfying
	\[
	\|T\|_{\mathfrak A_s}
	:=
	\left(\sum_{k\ge1}a_k(T)^s\right)^{1/s}<\infty.
	\]
\end{definition}

For a compact operator $T:X\to Y$ between Hilbert spaces,
$a_k(T)=s_k(T)$ for every $k\ge1$. Consequently, the approximation
ideal $\mathfrak A_s(X,Y)$ coincides with the Schatten class
$\mathcal S_s(X,Y)$ for $0<s<\infty$.
It is therefore natural to ask how the preceding
singular value and Schatten class characterizations
extend to Hankel operators between
(quasi-)Banach spaces.
However, despite the progress in the Hilbert space
setting, approximation number decay and approximation
ideal membership for Hankel operators from
$F^p(\varphi)$ to $L^q(\varphi)$ remain largely
unexplored.

Another related question is whether ideal membership is
preserved under conjugation of the symbol.
For bounded symbols, Berger and Coburn \cite{BC87}
proved that $H_f$ is compact on $F_\alpha^2$
if and only if $H_{\bar f}$ is compact.
This equivalence is known as the Berger--Coburn phenomenon.
Its Schatten class counterpart was established by
Bauer \cite{Bauer04} for the Hilbert--Schmidt class
and by Hu and Virtanen \cite{HV22,HVC23} for
$1<s<\infty$.
Fan, Wang and Zeng \cite[Theorem 1.11]{FWZ26}
established the corresponding equivalence for weak
Schatten classes.
In particular, for $f\in L^\infty(\C)$ and $1<s<\infty$,
\[
H_f\in\mathcal S_s(F_\alpha^2,L^2(\C,e^{-\alpha|z|^2}dv))
\quad\Longleftrightarrow\quad
H_{\bar f}\in\mathcal S_s(F_\alpha^2,L^2(\C,e^{-\alpha|z|^2}dv)),
\]
and
\[
H_f\in\mathcal S_{s,\infty}(F_\alpha^2,L^2(\C,e^{-\alpha|z|^2}dv))
\quad\Longleftrightarrow\quad
H_{\bar f}\in\mathcal S_{s,\infty}(F_\alpha^2,L^2(\C,e^{-\alpha|z|^2}dv)).
\]
The Schatten class equivalence fails in general
for $0<s\le1$ \cite{Xia23,HVBC23}.
	More generally, Xia \cite{Xia26} characterizes the
	Berger--Coburn phenomenon for symmetric norm ideals
	on the classical Hilbert Fock space in terms of Boyd indices.
	Beyond the Hilbert space setting, Hu and Lv \cite{HL26}
	characterize absolutely summing Hankel operators
	from $F_\alpha^p$ to $L_\alpha^p$ for $1\le p<\infty$
	and investigate the corresponding Berger--Coburn phenomenon.
These developments motivate the study of the same phenomenon
for approximation ideals and their weak counterparts
for Hankel operators from $F^p(\varphi)$ to $L^q(\varphi)$.

\subsection{Our results}
Building on \textup{(A)}--\textup{(C)}, we study Hankel operators
$H_f:F^p(\varphi)\to L^q(\varphi)$ for $0<p<\infty$,
$1\le q<\infty$, and possibly unbounded symbols $f\in\mathcal S$.
Here $\varphi\in C^2(\C^n)$ is a real-valued weight satisfying
$\operatorname{Hess}_{\mathbb R}\varphi\simeq E$, as specified in
\eqref{eq:weight} in Section~\ref{sec:preliminaries}.
Specifically, we characterize prescribed decay of the approximation
numbers of $H_f$ and membership in approximation ideals and their
weak counterparts in terms of the IDA function. We also characterize
simultaneous decay and simultaneous ideal membership for $H_f$ and
$H_{\bar f}$ in terms of the mean oscillation of $f$. For bounded
symbols, we establish the Berger--Coburn phenomenon for both
approximation ideals and their weak counterparts in sharp exponent
ranges.

Among our results, we highlight three novel contributions.
First, to the best of our knowledge, this work provides
the first unified characterization of these properties
in terms of local symbol conditions beyond the
Hilbert space setting. Second, we obtain necessary and sufficient conditions for prescribed approximation-number decay allowing $\gamma>1$ in the growth condition on $\Psi$.
This is new even when $p=q=2$.
The new method developed for this purpose also applies to
the broad class of weighted Bergman spaces considered
in \cite{FWZ26} and yields the converse to \textup{(E)};
see Subsection~\ref{sec:Bergman-converse}. Third, when $p<2<q$,
the Berger--Coburn equivalence holds for $s_0(p,q)<s\le1$,
where the sharp threshold $s_0(p,q)$ is strictly less than $1$.

We first recall the IDA spaces introduced by Hu and Virtanen
\cite{HV23}, which can be regarded as a natural evolution
of BMO spaces; see \cite{JN61,Zhu12}.

\begin{definition}\label{def:IDA}
	Let $1\le q<\infty$ and $r>0$.
	For $f\in L^q_{\rm loc}(\C^n)$, define the IDA function by
	\begin{equation}\label{eq:G}
		G_{q,r}(f)(z)
		:=\inf_{h\in\mathcal H(B(z,r))}
		\left(\frac1{v(B(z,r))}
		\int_{B(z,r)}|f(w)-h(w)|^q\,dv(w)\right)^{1/q},
		\qquad z\in\C^n,
	\end{equation}
	where $B(z,r):=\{w\in\C^n:|w-z|<r\}$. For a general measurable function $f$, we use the same definition,
	allowing the value $+\infty$.
	For $0<s\le\infty$, the space $\mathrm{IDA}^{s,q}$
	consists of all $f\in L^q_{\rm loc}(\C^n)$ such that
	\[
	\|f\|_{\mathrm{IDA}^{s,q}}
	:=\|G_{q,1}(f)\|_{L^s}<\infty.
	\]
	Write $\mathrm{BDA}^q:=\mathrm{IDA}^{\infty,q}$.
	The space $\mathrm{VDA}^q$ consists of all
	$f\in\mathrm{BDA}^q$ satisfying
	\[
	\lim_{|z|\to\infty}G_{q,1}(f)(z)=0.
	\]
\end{definition}

Next we recall the definitions of the distribution
function and non-increasing rearrangement
from \cite[Section 5.3]{BL76}.

\begin{definition}
	Let $f$ be a measurable function on $\C^n$.
	Its non-increasing rearrangement with respect to
	Lebesgue measure $dv$ is defined by
	\[
	f^*(t):=\inf\{s>0:d_f(s)\le t\},
	\qquad t\ge0,
	\]
	where
	$
	d_f(t):=v\bigl(\{z\in\C^n:|f(z)|>t\}\bigr)
	$
	is the distribution function of $f$.
\end{definition}

Our first main result characterizes approximation
number decay in three regimes determined by $p$ and $q$.
For the case $p<2<q$, we use the following indices.

\begin{definition}\label{def:lorentz-index1}
	Let $0<p<2<q<\infty$ and put
	\[
	p^\dagger:=
	\begin{cases}
		\infty,&0<p\le1,\\
		p/(p-1),&1<p<2.
	\end{cases}
	\]
	Define
	\begin{equation}\label{eq:rates-window-exponent}
		\tau:=\tau_{p,q}:=\min\{p^\dagger,q\}>2,
		\qquad
		\varsigma:=
		\begin{cases}
			2,&p^\dagger=q,\\
			\tau,&p^\dagger\ne q.
		\end{cases}
	\end{equation}
\end{definition}

\begin{theorem}\label{thm:rates-rearranged-decay}
Suppose $\operatorname{Hess}_{\mathbb R}\varphi\simeq E$.
	Let $0<p<\infty$ and $1\le q<\infty$.
	Let $\Psi:[1,\infty)\to(0,\infty)$ be nondecreasing
	such that $\Psi(t)/t^\gamma$ is nonincreasing
	for some $\gamma>0$.
	Assume $f\in\mathcal S\cap\mathrm{VDA}^q$ when $p\le q$ and
	$f\in\mathcal S\cap\mathrm{IDA}^{r,q}$ when $p>q$, where
	$1/r=1/q-1/p$.
	Then the following equivalences hold for all $k\in \mathbb N$.
	\begin{enumerate}
		\item[\textup{(i)}]
	If $p\le q\le2$ or $2\le p\le q$, then
		\begin{equation}\label{eq:rates-rearranged-same-side}
			a_k(H_f)\lesssim\frac1{\Psi(k)}
			\quad\Longleftrightarrow\quad
			(G_{q,1}(f))^*(k)\lesssim\frac1{\Psi(k)}.
		\end{equation}
		
		\item[\textup{(ii)}]
		If $p<2<q$, then
		\begin{equation}\label{eq:rates-rearranged-crossing}
			a_k(H_f)\lesssim\frac1{\Psi(k)}
			\quad\Longleftrightarrow\quad
			\left(
			\int_{k+1}^{(k+1)^{\tau/2}}
			u^{\varsigma/\tau-1}
			\bigl[(G_{q,1}(f))^*(u)\bigr]^\varsigma\,du
			\right)^{1/\varsigma}
			\lesssim\frac{\sqrt{k}}{\Psi(k)}.
		\end{equation}
		
		\item[\textup{(iii)}]
		If $p>q$, then
		\begin{equation}\label{eq:rates-rearranged-tail}
			a_k(H_f)\lesssim\frac1{\Psi(k)}
			\quad\Longleftrightarrow\quad
			\left(
			\int_k^\infty
			\bigl[(G_{q,1}(f))^*(u)\bigr]^r\,du
			\right)^{1/r}
			\lesssim\frac1{\Psi(k)}.
		\end{equation}
	\end{enumerate}
The implicit constants are independent of $k$.
\end{theorem}

\begin{remark}
	We make the following observations concerning
	Theorem~\ref{thm:rates-rearranged-decay}.
	\begin{enumerate}
		\item[\textup{(i)}]
		By \textup{(A)} and \textup{(B)}, the symbol assumptions
		in the theorem ensure that
		$H_f:F^p(\varphi)\to L^q(\varphi)$ is compact.
		
		\item[\textup{(ii)}]
		All three criteria remain equivalent when required
		only for sufficiently large $k$.
		By monotonicity and the growth assumption on $\Psi$,
		the rearrangement conditions may equivalently be
		formulated for real arguments $t\ge2$.
		
		\item[\textup{(iii)}]
		Under additional growth assumptions on $\Psi$,
		Corollary~\ref{cor:rates-direct-rearrangement}
		reduces the integral conditions in
		Theorem~\ref{thm:rates-rearranged-decay}\textup{(ii), (iii)}
		to pointwise bounds for $(G_{q,1}(f))^*$.
	\end{enumerate}
\end{remark}

We next turn to summability of approximation numbers.
The corresponding symbol conditions involve Lorentz
spaces \cite{Lorentz51}, Lorentz--Zygmund spaces
\cite{BR80,OP99}, and Copson--Lorentz spaces
\cite{Krepela18}, whose definitions we recall below.

\begin{definition}\label{def:function-Lorentz}
	\label{def:Lorentz-Zygmund}\label{def:critical-Copson}
	Let $0<r,s<\infty$ and let $f$ be a measurable
	function on $\C^n$.
	The Lorentz space $\mathcal L^{r,s}(\C^n,dv)$
	consists of those $f$ for which
	\begin{equation}\label{eq:general-rearrangement}
		\|f\|_{\mathcal L^{r,s}}
		:=
		\left(
		\int_0^\infty [t^{1/r}f^*(t)]^s
		\,\frac{dt}{t}
		\right)^{1/s}<\infty.
	\end{equation}
	For $\beta\in\mathbb R$, the Lorentz--Zygmund space
	$\mathcal L^{r,s,\beta}(\C^n,dv)$ is defined by
	\begin{equation}\label{eq:function-Lorentz-Zygmund}
		\|f\|_{\mathcal L^{r,s,\beta}}
		:=
		\left(
		\int_0^\infty
		[t^{1/r}\lambda(t)^\beta f^*(t)]^s
		\,\frac{dt}{t}
		\right)^{1/s}<\infty,
	\end{equation}
	where $\lambda(t):=1+\max\{0,\log t\}$.
	In particular, $\mathcal L^{r,s,0}=\mathcal L^{r,s}$.
	
	For positive measurable weights $u,w$ on $(0,\infty)$,
	the Copson--Lorentz space
	$\mathcal{CL}^{r,s}_{u,w}(\C^n,dv)$ is defined by
	\begin{equation}\label{eq:critical-Copson}
		\|f\|_{\mathcal{CL}^{r,s}_{u,w}}
		:=
		\left[
		\int_0^\infty
		\left(\int_t^\infty f^*(a)^r u(a)\,da\right)^{s/r}
		w(t)\,dt
		\right]^{1/s}<\infty.
	\end{equation}
	We write $\mathcal{CL}^{r,s}:=\mathcal{CL}^{r,s}_{1,w}$
	when $w(t)=(1+t)^{-1}$.
\end{definition}

These conditions extend the IDA condition, since
\[
\|G_{q,1}(f)\|_{\mathcal L^{s,s}}
=\|f\|_{\mathrm{IDA}^{s,q}}.
\]
The relevant Lorentz exponent depends on $p,q,s$
as follows.

\begin{definition}\label{def:lorentz-index}
	Let $0<p<\infty$, $1\le q<\infty$ and $0<s<\infty$.
	Put
	\[
	\mathcal E:=\{(p,q,s):0<p<2<q<\infty,\ s=2\}
	\]
	and
	\[
	\eta=\eta_{p,q}:=
	\max\left\{0,\min\left\{
	\frac1p-\frac12,\frac12-\frac1q
	\right\}\right\}.
	\]
	For $(p,q,s)\notin\mathcal E$, define
	\begin{equation}\label{eq:lorentz-index}
		\rho=\rho_{p,q,s}:=
		\begin{cases}
			s,
			&p\le q\le2\quad\text{or}\quad 2\le p\le q,\\[5pt]
			\displaystyle\frac{s}{1-s\eta},
			& p<2<q,\quad 0<s<2,\\[7pt]
			\displaystyle\frac{s}{1-2\eta},
			& p<2<q,\quad 2<s<\infty,\\[7pt]
			\displaystyle\frac{sr}{s+r},
			& p>q,\quad r=r_{p,q}:=\displaystyle\frac{pq}{p-q}.
		\end{cases}
	\end{equation}
\end{definition}

Our next result characterizes approximation ideal
membership in terms of these spaces.

\begin{theorem}\label{thm:increasing}
	Let $0<p<\infty$, $1\le q<\infty$ and $0<s<\infty$.
	Suppose that $\operatorname{Hess}_{\mathbb R}\varphi\simeq E$
	and $f\in\mathcal S$.
	\begin{enumerate}
		\item[\textup{(i)}]
		If $(p,q,s)\notin\mathcal E$, then
		$H_f\in\mathfrak A_s(\F p,\Lw q)$ if and only if
		$G_{q,1}(f)\in\mathcal L^{\rho,s}(\C^n,dv)$.
		Moreover,
		\begin{equation}\label{eq:ascending-rearrangement-main}
			\|H_f\|_{\mathfrak A_s}
			\simeq\|G_{q,1}(f)\|_{\mathcal L^{\rho,s}}.
		\end{equation}
		
		\item[\textup{(ii)}]
		Suppose $(p,q,s)\in\mathcal E$.
		If $p^\dagger=q$, then
		$H_f\in\mathfrak A_2(\F p,\Lw q)$ if and only if
		$G_{q,1}(f)\in\mathcal L^{q,2,1/2}(\C^n,dv)$.
		Moreover,
		\begin{equation}\label{eq:critical-selfdual-main}
			\|H_f\|_{\mathfrak A_2}
			\simeq\|G_{q,1}(f)\|_{\mathcal L^{q,2,1/2}}.
		\end{equation}
		If $p^\dagger\ne q$, then
		$H_f\in\mathfrak A_2(\F p,\Lw q)$ if and only if
		$
		G_{q,1}(f)\in
		\mathcal L^{\tau,\tau}(\C^n,dv)
		\cap\mathcal{CL}^{\tau,2}(\C^n,dv).
		$
		Moreover,
		\begin{equation}\label{eq:critical-main-upper}
			\|H_f\|_{\mathfrak A_2}
			\simeq
			\|G_{q,1}(f)\|_{\mathcal L^{\tau,\tau}}
			+\|G_{q,1}(f)\|_{\mathcal{CL}^{\tau,2}}.
		\end{equation}
	\end{enumerate}
	The comparison constants are independent of $f$.
\end{theorem}

\begin{remark}\label{rem:Hu-comparison}
	We make the following observations concerning
	Theorem~\ref{thm:increasing}.
	\begin{enumerate}
	\item[\textup{(i)}] Either condition implies that $f\in L^q_{\rm loc}(\C^n)$;
	see the beginning of the proof.

		\item[\textup{(ii)}]
		If $p\le q\le2$ or $2\le p\le q$,
		then $\rho=s$, so the criterion reduces to
		\begin{equation*}
			H_f\in\mathfrak A_s(\F p,\Lw q)
			\quad\Longleftrightarrow\quad
			f\in\mathrm{IDA}^{s,q}.
		\end{equation*}
		This includes every $1\le p=q<\infty$.
		In particular, at $p=q=2$, it recovers
		\textup{(C)}.
		
		\item[\textup{(iii)}]
		In the critical regime $p<2<q$, $s=2$, let
		$I_N:\ell_N^p\to\ell_N^q$ be the identity.
		Lemma~\ref{lem:finite-crossing} gives
		\[
		\sum_{k=1}^N a_k(I_N)^2
		\simeq N^{2/\tau}\log(eN),
		\]
		whereas the squared $\ell^{\tau,2}$-norm of $N$ unit
		entries is $\sum_{k=1}^N k^{2/\tau-1}\simeq N^{2/\tau}$.
		For general diagonal sequences,
		Lemma~\ref{lem:same-side} establishes a Lorentz--Zygmund
		criterion when $p^\dagger=q$ and a tail criterion when
		$p^\dagger\ne q$.
		The sampling comparisons \eqref{eq:logarithmic-sampling}
		and \eqref{eq:critical-Copson-sampling} translate these
		into the function spaces in
		Theorem~\ref{thm:increasing}\textup{(ii)}.
	\end{enumerate}
\end{remark}

We next consider weak approximation ideals, extending
the weak Schatten class characterization in
\cite[Theorem 1.10]{FWZ26} to the present setting.
For $0<s<\infty$, define
\[
\mathfrak A_{s,\infty}(X,Y):=
\left\{
T\in\mathcal L(X,Y):
\|T\|_{\mathfrak A_{s,\infty}}
:=\sup_{k\ge1}k^{1/s}a_k(T)<\infty
\right\}.
\]
The corresponding weak Lorentz space is defined by
\[
\mathcal L^{\sigma,\infty}(\C^n,dv):=
\left\{
f\text{ measurable}:
\|f\|_{\mathcal L^{\sigma,\infty}}
:=\sup_{t>0}t^{1/\sigma}f^*(t)<\infty
\right\},
\qquad 0<\sigma<\infty.
\]
For $\tau>2$ and $1\le\varsigma<\infty$, we also define
$\mathcal X^{\tau,\varsigma}(\C^n,dv)$ by
\[
\|f\|_{\mathcal X^{\tau,\varsigma}}
:=
\|f\|_{\mathcal L^{\tau,\infty}}
+\sup_{t\ge2}
\left(
\int_t^{t^{\tau/2}}
u^{\varsigma/\tau-1}f^*(u)^\varsigma\,du
\right)^{1/\varsigma}<\infty.
\]

\begin{theorem}\label{cor:weak-approximation-ideals}
	Let $0<p<\infty$, $1\le q<\infty$ and $0<s<\infty$.
	Suppose $\operatorname{Hess}_{\mathbb R}\varphi\simeq E$
	and $f\in\mathcal S$.
	\begin{enumerate}
		\item[\textup{(i)}]
		If $(p,q,s)\notin\mathcal E$, then
		$H_f\in\mathfrak A_{s,\infty}(\F p,\Lw q)$
		if and only if
		$G_{q,1}(f)\in\mathcal L^{\rho,\infty}(\C^n,dv)$.
		Moreover,
		\[
		\|H_f\|_{\mathfrak A_{s,\infty}}
		\simeq\|G_{q,1}(f)\|_{\mathcal L^{\rho,\infty}}.
		\]
		
		\item[\textup{(ii)}]
		If $(p,q,s)\in\mathcal E$, then
		$H_f\in\mathfrak A_{2,\infty}(\F p,\Lw q)$
		if and only if
		$G_{q,1}(f)\in\mathcal X^{\tau,\varsigma}(\C^n,dv)$.
		Moreover,
		\[
		\|H_f\|_{\mathfrak A_{2,\infty}}
		\simeq\|G_{q,1}(f)\|_{\mathcal X^{\tau,\varsigma}}.
		\]
	\end{enumerate}
	The comparison constants are independent of $f$.
\end{theorem}

Our results also have consequences for the pair
$H_f$ and $H_{\bar f}$, whose simultaneous boundedness,
compactness, Schatten class membership, and singular
value decay have been studied in various settings
\cite{FWZ26,HV22,HV23,Isralowitz13,Pau16,Xia01,Xia02,
	ZHW25,Zhu91,Zhu07}.
Specifically, applying Theorems~\ref{thm:rates-rearranged-decay},
\ref{thm:increasing}, and \ref{cor:weak-approximation-ideals}
to $H_f$ and $H_{\bar f}$ yields criteria for simultaneous
decay of their approximation numbers and simultaneous membership
in approximation ideals and their weak counterparts. These
criteria are expressed in terms of the mean oscillation function
$\operatorname{MO}_{q,r}(f)$ (see Subsection~\ref{subsec:paired}).

For bounded symbols, we establish the Berger--Coburn
equivalence for
$\mathfrak A_s(F^p(\varphi),L^q(\varphi))$
and $\mathfrak A_{s,\infty}(F^p(\varphi),L^q(\varphi))$
in the sharp range
\begin{equation}\label{eq:BC-range}
	\begin{aligned}
		\mathcal B:= {}&
		\left\{(p,q,s):0<p\le q<\infty,\ q\ge1,\
		\frac1{1+\eta}<s<\infty\right\}\\
		&{}\cup
		\left\{(p,q,s):1\le q<p<\infty,\
		\frac r{r-1}<s<\infty\right\}.
	\end{aligned}
\end{equation}
Notably, when $p<2<q$, this range includes
$$(1+\eta)^{-1}<s\le1,$$ extending
the Berger--Coburn phenomenon below $s=1$.
The precise result is stated as follows.

\begin{theorem}\label{thm:BC-approx}
	Suppose that $\operatorname{Hess}_{\mathbb R}\varphi\simeq E$
	and $f\in L^\infty(\C^n)$.
	Let $(p,q,s)\in\mathcal B$.
	\begin{enumerate}
		\item[\textup{(i)}]
		$H_f\in\mathfrak A_s(\F p,\Lw q)$
		if and only if
		$H_{\bar f}\in\mathfrak A_s(\F p,\Lw q)$.
		
		\item[\textup{(ii)}]
		$H_f\in\mathfrak A_{s,\infty}(\F p,\Lw q)$
		if and only if
		$H_{\bar f}\in\mathfrak A_{s,\infty}(\F p,\Lw q)$.
	\end{enumerate}
\end{theorem}

\begin{remark}\label{rem:BC-Hilbert}\label{rem:BC-coverage}
	We make the following observations concerning
	Theorem~\ref{thm:BC-approx}.
	\begin{enumerate}
		\item[\textup{(i)}]
		If $p\le q\le2$ or $2\le p\le q$, then $\eta=0$, so the range is $s>1$.
		In particular, at $p=q=2$, the theorem recovers
		the Berger--Coburn equivalences for
		$\mathcal S_s(F^2(\varphi),L^2(\varphi))$ and
		$\mathcal S_{s,\infty}(F^2(\varphi),L^2(\varphi))$.
		
		\item[\textup{(ii)}]
		The lower bounds on $s$ in \eqref{eq:BC-range} are sharp
	in every complex dimension for both $\mathfrak A_s$
	and $\mathfrak A_{s,\infty}$;
	see Theorem~\ref{thm:BC-counter}.
	\end{enumerate}
\end{remark}

\subsection{Difficulties and strategy}
Two difficulties arise in passing from singular values
to approximation numbers.
First, the convex trace method in \cite{FWZ26}
does not provide the reverse implication for general
fast decay.
Second, the Schatten class arguments in \cite{HV22}
use the positive operator $H_f^*H_f$ and orthonormal
systems, tools that are not directly available
for operators between non-Hilbert spaces.

Our approach compares $H_f$ with a diagonal operator
whose entries are sampled values of the IDA function.
For the upper estimates, Lemma~\ref{lem:transfer}
uses local Taylor approximations with geometrically
decaying errors and polynomially growing ranks
to control multiplication operators.
Combined with the symbol decomposition and
$\bar\partial$ estimates in \cite{HV23}, this transfers
decay bounds from the diagonal operator to $H_f$.
For the reverse estimates, local kernel estimates
recover the diagonal approximation numbers from
those of $H_f$, up to a fixed change of index.
Together, these arguments yield
Proposition~\ref{thm:rates-diagonal-transfer}
for arbitrary $\gamma>0$.
A summable majorant then gives the ideal comparison
in Proposition~\ref{prop:diagonal-equivalence}.

The remaining problem is to characterize the diagonal
approximation numbers.
We combine finite-dimensional estimates
\cite{Gluskin84,ZF14} with dyadic block arguments,
treating the critical regime $p<2<q$, $s=2$
separately.
Sampling and rearrangement comparisons then translate
the resulting sequence conditions into the stated
decay criteria and Lorentz, Lorentz--Zygmund,
and Copson--Lorentz conditions.

\subsection{Organization of the paper}
Section~\ref{sec:preliminaries} collects the preliminaries
on weighted Fock spaces, Hankel operators, lattice coverings,
reproducing kernels, and $\bar\partial$ estimates.
Section~\ref{sec:diagonal} establishes decay and approximation
ideal criteria for diagonal operators.
Section~\ref{sec:sampling} develops radius and sampling
comparisons for the IDA function and approximation estimates
for multiplication operators.
Section~\ref{sec:Hankel} establishes
Proposition~\ref{thm:rates-diagonal-transfer} and proves
Theorems~\ref{thm:rates-rearranged-decay},
\ref{thm:increasing}, and
\ref{cor:weak-approximation-ideals}, together with
the corresponding results for paired Hankel operators.
Section~\ref{sec:BC} proves
Theorem~\ref{thm:BC-approx} and establishes the sharpness
of its exponent ranges in Theorem~\ref{thm:BC-counter}.

\section{Preliminaries}\label{sec:preliminaries}
\subsection{Notation}
The letter $C$ denotes a positive constant whose value may vary
from line to line.
We write $A\lesssim B$ if $A\le CB$ for some constant $C>0$,
and $A\simeq B$ if $A\lesssim B$ and $B\lesssim A$.

For $F\subset\C^n$, let $\chi_F$ denote its indicator function.
For quasi-Banach spaces $X$ and $Y$, let $\mathcal L(X,Y)$ denote
the space of bounded linear operators from $X$ to $Y$. For
$T\in\mathcal L(X,Y)$, we write $\|T\|_{X\to Y}$ for its operator
norm and abbreviate it to $\|T\|$ when $X$ and $Y$ are clear from
the context. If $X$ and $Y$ are Banach spaces, the adjoint
$T^*:Y^*\to X^*$ is defined by $T^*y^*=y^*\circ T$ for $y^*\in Y^*$.
For a matrix $T$, the notation $T^*$ denotes its conjugate transpose.

For $1\le p\le\infty$, we denote by $p'$ the conjugate
exponent determined by $1/p+1/p'=1$, with the convention
that $1/\infty=0$. For $t\in\mathbb R$, let $\lfloor t\rfloor$ and $\lceil t\rceil$
denote the greatest integer not exceeding $t$ and the least
integer not smaller than $t$, respectively. 	
\subsection{Weighted Fock spaces}

For $z=(z_1,\ldots,z_n)$ and $w=(w_1,\ldots,w_n)$ in $\C^n$,
write
\[
z\cdot\bar w:=\sum_{j=1}^n z_j\bar w_j,
\qquad |z|:=(z\cdot\bar z)^{1/2}.
\]
The spaces $L^p(\varphi)$ and $F^p(\varphi)$ are defined
in the introduction. Let $\varphi\in C^2(\C^n)$ be real-valued. We write
$\operatorname{Hess}_{\mathbb R}\varphi$ for the real Hessian of
$\varphi$. In the real coordinates
$z_j=x_{2j-1}+\mathrm i x_{2j}$, $1\le j\le n$, it is given by
\[
\operatorname{Hess}_{\mathbb R}\varphi
=
\left(
\frac{\partial^2\varphi}{\partial x_j\partial x_k}
\right)_{j,k=1}^{2n}.
\]
We assume that there exist constants $0<m\le M<\infty$ such that
\begin{equation}\label{eq:weight}
mE\le \operatorname{Hess}_{\mathbb R}\varphi(z)\le ME,
\qquad z\in\C^n,
\end{equation}
where $E$ is the $2n\times2n$ identity matrix and the inequalities
are in the sense of quadratic forms. We write
$\operatorname{Hess}_{\mathbb R}\varphi\simeq E$ for \eqref{eq:weight}. This assumption includes the Gaussian
weights $\varphi(z)=\alpha|z|^2/2$, $\alpha>0$, which give the
classical Fock spaces $F_\alpha^p$ \cite{Zhu12}, and the weight
$\varphi(z)=|z|^2-\frac12\log(1+|z|^2)$, which gives a
Fock--Sobolev space \cite{CZ12}. It also allows nonradial weights
such as $\varphi(z)=|z|^2+\sin(\operatorname{Re}z_1)$.

The following weighted submean estimate is
\cite[Lemma 2.1]{HV23}.

\begin{lemma}\label{lem:weighted-local}
	Suppose $\operatorname{Hess}_{\mathbb R}\varphi\simeq E$,
	$0<p<\infty$, and $r>0$. Then
	\[
	|g(z)|^pe^{-p\varphi(z)}
	\le C\int_{B(z,r)}
	|g(\xi)|^pe^{-p\varphi(\xi)}\,dv(\xi),
	\qquad g\in F^p(\varphi),\quad z\in\C^n,
	\]
	where $C$ is independent of $g$ and $z$.
\end{lemma}
In particular, for $0<p\le q<\infty$,
\begin{equation}\label{eq:Fock-inclusion}
	\|g\|_{q,\varphi}\lesssim\|g\|_{p,\varphi},
	\qquad F^p(\varphi)\subseteq F^q(\varphi).
\end{equation}

Let $K(z,w)$ denote the reproducing kernel of $F^2(\varphi)$.
We recall the following estimates from
\cite[Lemma 2.2]{HV23}.

\begin{lemma}\label{lem:known-kernel}
	Suppose $\operatorname{Hess}_{\mathbb R}\varphi\simeq E$.
	There exist $C_1,C_2,\theta,r_0>0$ such that
	\begin{equation}\label{eq:kernel-upper}
		|K(z,w)|
		\le C_1e^{\varphi(z)+\varphi(w)-\theta|z-w|},
		\qquad z,w\in\C^n,
	\end{equation}
	and
	\begin{equation}\label{eq:kernel-lower}
		|K(z,w)|
		\ge C_2e^{\varphi(z)+\varphi(w)},
		\qquad z\in\C^n,\quad w\in B(z,r_0).
	\end{equation}
	The constants are independent of $z$ and $w$.
\end{lemma}

For $z\in\C^n$, define the normalized reproducing kernel by
\[
k_z(w):=\frac{K(w,z)}{\sqrt{K(z,z)}}.
\]
Taking $w=z$ and integrating
\eqref{eq:kernel-upper} and \eqref{eq:kernel-lower} gives,
for every $0<p<\infty$,
\begin{equation}\label{eq:normalized-kernel}
	K(z,z)\simeq e^{2\varphi(z)},\qquad
	\|K(\cdot,z)\|_{p,\varphi}\simeq e^{\varphi(z)},\qquad
	\|k_z\|_{p,\varphi}\simeq1.
\end{equation}
The comparison constants are independent of $z$.

\subsection{The Bergman projection}

Let $P:L^2(\varphi)\to F^2(\varphi)$ be the orthogonal
projection, given by
\begin{equation}\label{eq:Bergman-projection}
	Pu(z):=\int_{\C^n}
	K(z,w)u(w)e^{-2\varphi(w)}\,dv(w).
\end{equation}
By \eqref{eq:kernel-upper}, the same integral defines
a bounded operator on $L^p(\varphi)$ for $1\le p<\infty$,
and
\begin{equation}\label{eq:projection-bounded}
	\|Pu\|_{p,\varphi}\lesssim\|u\|_{p,\varphi},
	\qquad u\in L^p(\varphi).
\end{equation}
It also reproduces every $g\in F^p(\varphi)$:
\begin{equation}\label{eq:projection-identity}
	Pg=g,\qquad g\in F^p(\varphi),\quad 0<p<\infty;
\end{equation}
see \cite[Section 2B]{HV23}.

\subsection{Hankel operators}

To define Hankel operators with possibly unbounded symbols, put
\begin{equation}\label{eq:symbol}
	\Gamma:=
	\left\{
	\sum_{j=1}^N c_jK(\cdot,a_j):
	N\in\mathbb N,\ c_j\in\C,\ a_j\in\C^n
	\right\}
\end{equation}
and
\begin{equation}\label{eq:symbol1}
	\mathcal S:=
	\left\{
	f\text{ measurable on }\C^n:
	fg\in L^1(\varphi)\text{ for every }g\in\Gamma
	\right\}.
\end{equation}
For $f\in\mathcal S$, the Hankel operator
$H_f:=(I-P)M_f$ is well-defined on $\Gamma$.
By \cite[Proposition 2.5]{HV20}, $\Gamma$ is dense in
$F^p(\varphi)$ for every $0<p<\infty$.
Thus $H_f$ is densely defined on $F^p(\varphi)$.

\subsection{Lattices in \texorpdfstring{$\C^n$}{complex Euclidean space}}

We recall the lattice notation from \cite[Section 2D]{HV23}.
For $r>0$, a sequence $\{a_k\}_{k\ge1}\subset\C^n$
is called an $r$-lattice if
\begin{equation}\label{eq:lattice-definition}
	\C^n=\bigcup_{k\ge1}B(a_k,r),
	\qquad
	B\left(a_j,\frac{r}{2\sqrt n}\right)
	\cap B\left(a_k,\frac{r}{2\sqrt n}\right)
	=\varnothing
	\quad(j\ne k).
\end{equation}
A typical example is
\begin{equation}\label{eq:lattice-model}
	\left\{
	\frac r{\sqrt n}
	(m_1+\mathrm i\ell_1,\ldots,m_n+\mathrm i\ell_n):
	m_j,\ell_j\in\mathbb Z,\ 1\le j\le n
	\right\}.
\end{equation}
The following overlap property is given in \cite{Zhu12}
for $n=1$ and extends to $\C^n$ by the same argument.
There exists an integer $N$, independent of $r$ and
the lattice, such that
\begin{equation}\label{eq:lattice-overlap}
	1\le\sum_{k\ge1}\chi_{B(a_k,2r)}(z)\le N,
	\qquad z\in\C^n.
\end{equation}
\subsection{The \texorpdfstring{$\bar\partial$}{d-bar} equation}

Following \cite[(2-12)]{HV23}, with the normalization below, for a smooth compactly
supported $(0,1)$-form $\omega$, define
\begin{equation}\label{eq:A-formula}
	A_\varphi(\omega)(z)
	:=
	-\frac{1}{(2\pi i)^n}\int_{\C^n}
	e^{\langle2\partial\varphi(\xi),z-\xi\rangle}
	\sum_{j=0}^{n-1}
	\frac{
		\omega(\xi)\wedge\partial|\xi-z|^2
		\wedge(2\bar\partial\partial\varphi(\xi))^j
		\wedge(\bar\partial\partial|\xi-z|^2)^{n-1-j}
	}{j!\,|\xi-z|^{2n-2j}},
\end{equation}
where
\[
\langle\partial\varphi(\xi),z-\xi\rangle
:=\sum_{\ell=1}^n
\frac{\partial\varphi}{\partial\xi_\ell}(\xi)
(z_\ell-\xi_\ell).
\]
Here $L^p_{(0,1)}(\varphi)$ denotes the space of
$(0,1)$-forms whose pointwise Euclidean norm belongs
to $L^p(\varphi)$.
We recall the following consequences of
\cite[Lemma 2.4 and Corollary 2.5]{HV23}.

\begin{lemma}\label{lem:known-solution}
	Suppose $\operatorname{Hess}_{\mathbb R}\varphi\simeq E$.
	\begin{enumerate}
		\item[(i)]
		For every $1\le p\le\infty$, $A_\varphi$ extends
		to a bounded operator from $L^p_{(0,1)}(\varphi)$
		to $L^p(\varphi)$, and
		\begin{equation}\label{eq:A}
			\|A_\varphi\omega\|_{p,\varphi}
			\le C\|\omega\|_{p,\varphi},
			\qquad \omega\in L^p_{(0,1)}(\varphi).
		\end{equation}
		The constant is independent of $\omega$.
	
		\item[(ii)]
		If $u\in C^2(\C^n)\cap\mathcal S$ and
		$|\bar\partial u|\in L^r(\C^n,dv)$ for some
		$1\le r\le\infty$, then
		\begin{equation}\label{eq:Aidentity}
			H_ug=(I-P)A_\varphi(g\bar\partial u),
			\qquad g\in\Gamma.
		\end{equation}
	\end{enumerate}
\end{lemma}

\section{Approximation numbers of diagonal operators}\label{sec:diagonal}

In this section, we characterize summability of approximation
numbers of diagonal operators in terms of Lorentz sequence spaces.
These characterizations will be used to prove the main theorems. 

Throughout this section and the subsequent sections, we use the notation
$p^\dagger$, $\tau$, $\varsigma$, $\eta$, $\rho$, $r$,
and $\mathcal E$ as defined in
Definitions~\ref{def:lorentz-index1} and~\ref{def:lorentz-index}.

\subsection{Basic estimates for approximation numbers}

We recall the following standard estimates for approximation
numbers; see \cite[Section 1, pp.~84--86]{Linde85} and
\cite[Section 2, properties (PS1)--(PS3)]{ZFH14}
for the quasi-Banach setting. 

\begin{lemma}\label{lem:ideal-properties}
	Let $0<s<\infty$. Let $A:X_0\to X$, $B:Y\to Y_0$,
	and $S,T:X\to Y$ be bounded linear operators, where
	$X_0,X$ are quasi-Banach spaces and $Y,Y_0$ are Banach spaces.
	Then the following hold:
	\begin{enumerate}
		\item[(i)]
		$a_k(BTA)\le\|B\|\|A\|a_k(T)$
		for $k\in\mathbb N$.
		
		\item[(ii)]
		$a_{k+l-1}(S+T)\le a_k(S)+a_l(T)$
		for $k,l\in\mathbb N$.
		
		\item[(iii)]
		$\As{S+T}^{s}\le2(\As S^{s}+\As T^{s})$
		for $0<s\le1$ and $S,T\in\mathfrak A_s(X,Y)$.
		
		\item[(iv)]
		$|a_k(S)-a_k(T)|\le\|S-T\|$
		for $k\in\mathbb N$.
		
		\item[(v)]
		$\As{S+T}\le C_s(\As S+\As T)$
		for $S,T\in\mathfrak A_s(X,Y)$, where
		\[
		C_s=
		\begin{cases}
			2^{2/s-1},&0<s\le1,\\
			2,&s>1.
		\end{cases}
		\]
	\end{enumerate}
\end{lemma}

For a nonempty finite or countable index set $J$ and $0<s<\infty$,
let $\ell_J^s$ be the space of sequences $\mathbf{x}=(x_j)_{j\in J}$ such that
\[
\|\mathbf{x}\|_{\ell_J^s}
:=\left(\sum_{j\in J}|x_j|^s\right)^{1/s}<\infty.
\]
 We also use
$\|\mathbf{x}\|_{\ell_J^\infty}:=\sup_{j\in J}|x_j|$. We let $\mathbb N$ denote the set of positive integers. We write $\ell^s$ when $J=\mathbb N$ and $\ell_N^s$ when
$J=\{1,\ldots,N\}$. We denote by $c_0$ the space of sequences converging to zero.

For $N\in\mathbb N$, let $I_N$ denote the $N\times N$ identity
matrix. When written as $I_N:\ell_N^p\to\ell_N^q$, it denotes
the map $\mathbf x\mapsto\mathbf x$ between the indicated
sequence spaces. For a finite index set $I$, we write $|I|$ for the number of
indices in $I$. More generally, $\#E$ denotes the cardinality
of a set $E$, so $|I|=\#I$ for finite index sets.

Let $\mathbf b:=(b_j)_{j\in\mathbb N}$ be a nonnegative sequence.
The diagonal operator $D_{\mathbf b}$ acts on sequences
$\mathbf{c}=(c_j)_{j\in\mathbb N}$ by
\begin{equation}\label{diaope}
	D_{\mathbf b}\mathbf{c}:=(b_jc_j)_{j\in\mathbb N}.
\end{equation}
The following norm formulas for diagonal operators are given in \cite[p.~5, (4)]{Fischer20}.

\begin{lemma}\label{lem:diag-facts}
Let $0<p<\infty$ and $1\le q<\infty$. Suppose that $\mathbf b=(b_j)_{j\in\mathbb N}$ is nonnegative and $D_{\mathbf b}:\ell^p\to\ell^q$ is bounded.
\begin{enumerate}
	\item[(i)] If $p\le q$, then $\|D_{\mathbf b}\|=\|\mathbf b\|_{\ell^\infty}$.
	\item[(ii)] If $p>q$, then $\|D_{\mathbf b}\|=\|\mathbf b\|_{\ell^r}$, where $1/r=1/q-1/p$.
\end{enumerate}
\end{lemma}

The next lemma follows from Fischer's result.
\begin{lemma}\label{lem:diag-domination}
Let $0<p<\infty$, $1\le q<\infty$, and $0<s<\infty$. Suppose that $\mathbf b=(b_j)_{j\in\mathbb N}$ and $\mathbf g=(g_j)_{j\in\mathbb N}$ are nonnegative sequences and that $D_{\mathbf b},D_{\mathbf g}:\ell^p\to\ell^q$ are bounded. If $b_j\le Cg_j$ for every $j\in\mathbb N$ and some $C>0$, then $\As{D_{\mathbf b}}\le C\As{D_{\mathbf g}}$.
\end{lemma}
\begin{proof}
	Define $\boldsymbol{\eta}=(\eta_j)_{j\in\mathbb N}$ by
	\[
	\eta_j=
	\begin{cases}
		b_j/g_j,&g_j>0,\\
		0,&g_j=0.
	\end{cases}
	\]
	Then $D_{\mathbf b}=D_{\boldsymbol{\eta}}D_{\mathbf{\mathbf g}}$ and
	$\|D_{\boldsymbol{\eta}}\|=\|\boldsymbol{\eta}\|_{\ell^\infty}\le C$.
	By Lemma~\ref{lem:ideal-properties}\textup{(i)},
	\[
	a_k(D_{\mathbf b})\le C a_k(D_{\mathbf{\mathbf g}}),\qquad
	\As{D_{\mathbf b}}\le C\As{D_{\mathbf{\mathbf g}}}.
	\]
	This completes the proof.
\end{proof}

For $A=(a_{ij})\in\C^{N\times N}$, its Hilbert--Schmidt norm is
\[
\|A\|_{\mathrm{HS}}
:=\left(\sum_{i,j=1}^N|a_{ij}|^2\right)^{1/2}
=\left(\sum_{\nu=1}^N\|Ax_\nu\|_{\ell_N^2}^2\right)^{1/2},
\]
where $(x_\nu)_{\nu=1}^N$ is any Euclidean orthonormal basis
of $\C^N$. We denote by $e_j$ the $j$th canonical unit vector in
$\C^N$ when $1\le j\le N$, and in the sequence spaces
when $j\in\mathbb N$, according to the context.

\begin{lemma}\label{lem:finite-crossing}
	Let $0<p<2<q<\infty$.
	For $N\in\mathbb N$, let $I_N:\ell_N^p\to\ell_N^q$
	be the identity map. Then
	\begin{equation}\label{eq:finite-crossing-width}
		a_k(I_N)\simeq\min\{1,N^{1/\tau}k^{-1/2}\},
		\qquad 1\le k\le N/2.
	\end{equation}
	Moreover, for every $0<s<\infty$,
	\begin{equation}\label{eq:finite-crossing-sum}
		\sum_{k=1}^{N}a_k(I_N)^s\simeq
		\begin{cases}
			N^{1-s(1/2-1/\tau)},&0<s<2,\\
			N^{2/\tau}(1+\log N),&s=2,\\
			N^{2/\tau},&s>2.
		\end{cases}
	\end{equation}
	The comparison constants are independent of $N$ and $k$.
\end{lemma}

\begin{proof}
	Let $A\in\mathbb C^{N\times N}$, regarded as an operator
	from $\ell_N^p$ to $\ell_N^q$.
	For $0<p\le1$ and
	$\mathbf{x}=(x_j)_{j=1}^N\in\ell_N^p$,
	the triangle inequality gives
	\[
	\|A\mathbf{x}\|_{\ell_N^q}
	\le\max_{1\le j\le N}\|Ae_j\|_{\ell_N^q}
	\sum_{j=1}^N|x_j|
	\le\max_{1\le j\le N}\|Ae_j\|_{\ell_N^q}
	\|\mathbf{x}\|_{\ell_N^p}.
	\]
 Testing $\mathbf{x}=e_j$ yields
	\begin{equation}\label{eq:small-source-matrix}
		\|A\|
		=\max_{1\le j\le N}\|Ae_j\|_{\ell_N^q}
		=\|A\|_{\ell_N^1\to\ell_N^q}.
	\end{equation}
	Thus the approximation numbers of $I_N$ for $0<p\le1$
	agree with those for $p=1$.
	By \cite[Lemma 3.8]{ZF14},
	\begin{equation}\label{eq:finite-crossing-upper}
		a_k(I_N)\le C_{p,q}\min\{1,N^{1/\tau}k^{-1/2}\},
		\qquad 1\le k\le N.
	\end{equation}
	The reverse estimate for $1\le k\le N/4$ follows from
	\cite[Lemma 3.5(ii)]{ZF14}
	and \eqref{eq:small-source-matrix}.
	
	For $N/4<k\le N/2$, let
	$S:\ell_N^p\to\ell_N^q$ satisfy $\rank S<k$,
	and put $A=I_N-S$ and $d=\dim\ker S$.
	Choose a Euclidean orthonormal basis
	$(x_\nu)_{\nu=1}^N$ of $\mathbb C^N$
	with $x_1,\ldots,x_d\in\ker S$.
	Then
	\begin{equation}\label{eq:crossing-HS}
		\|A\|_{\mathrm{HS}}^2
		=\sum_{\nu=1}^N\|Ax_\nu\|_{\ell_N^2}^2
		\ge d=N-\rank S\ge N/2.
	\end{equation}
	Since $q>2$,
	\begin{equation}\label{eq:crossing-frobenius}
		\|A\|
		\ge\max_{1\le j\le N}\|Ae_j\|_{\ell_N^q}
		\ge N^{1/q-1}\|A\|_{\mathrm{HS}}
		\ge2^{-1/2}N^{1/q-1/2}.
	\end{equation}
	 If $1<p<2$, then $p^\dagger=p'>2$. Applying the argument
of \eqref{eq:crossing-frobenius} to
$A^*:\ell_N^{q'}\to\ell_N^{p'}$ gives
	\begin{equation}\label{eq:crossing-frobenius-dual}
		\|A\|
		=\|A^*\|
		\ge N^{1/p^\dagger-1}\|A^*\|_{\mathrm{HS}}
		\ge2^{-1/2}N^{1/p^\dagger-1/2}.
	\end{equation}
	Put $\eta:=1/2-1/\tau$.
	Equations \eqref{eq:crossing-frobenius}
	and \eqref{eq:crossing-frobenius-dual} imply
	\[
	\|I_N-S\|
	\ge2^{-1/2}N^{-\eta},
	\qquad
	\min\{1,N^{1/\tau}k^{-1/2}\}\le2N^{-\eta}.
	\]
	Taking the infimum over $\rank S<k$ and using
	\eqref{eq:finite-crossing-upper} proves
	\eqref{eq:finite-crossing-width}.
	
	For \eqref{eq:finite-crossing-sum}, put $K_N=N^{2/\tau}$.
	By \eqref{eq:finite-crossing-upper},
	\begin{equation}\label{eq:finite-crossing-geometric}
		\sum_{k=1}^N a_k(I_N)^s
		\le C_{p,q,s}\left(
		K_N+N^{s/\tau}\sum_{K_N<k\le N}k^{-s/2}
		\right).
	\end{equation}
	Since $\tau>2$, comparison with integrals gives
	\[
	K_N+N^{s/\tau}\sum_{K_N<k\le N}k^{-s/2}
	\lesssim 
	\begin{cases}
		N^{1-s\eta},&0<s<2,\\
		N^{2/\tau}(1+\log N),&s=2,\\
		N^{2/\tau},&s>2,
	\end{cases}
	\]
	where the implicit constant is independent of $N$. 

For the reverse estimate, let $N\ge2$,
	$M=\lfloor N/2\rfloor$, and $L=\lfloor M^{2/\tau}\rfloor$.
	By \eqref{eq:finite-crossing-width},
	\[
	\sum_{k=1}^N a_k(I_N)^s
	\ge c_{p,q,s}\left(
	L+M^{s/\tau}\sum_{k=L+1}^M k^{-s/2}
	\right).
	\]
	The same integral comparison, together with $M\simeq N$,
	gives the lower estimate.
	For $N=1$, use $a_1(I_1)=1$.
	This completes the proof.
\end{proof}

\begin{lemma}\label{lem:critical-weak-approximation}
	Let $0<p<2<q<\infty$ with $p^\dagger\ne q$.
	There exist $C>0$ and $K\in\mathbb N$ such that,
	for every $\mathbf b\in\ell^\tau$ and $m\in\mathbb N$,
	the diagonal operator $D_{\mathbf b}:\ell^p\to\ell^q$ satisfies
	\begin{equation}\label{eq:critical-weak-approximation}
		a_{Km+1}(D_{\mathbf b})
		\le Cm^{-1/2}\|\mathbf b\|_{\ell^\tau}.
	\end{equation}
	The constants $C$ and $K$ are independent of $\mathbf b$ and $m$.
\end{lemma}

\begin{proof}
Suppose first that $1<p<2<q$.
For $\mathbf b\in\mathbb C^N$, we construct
$R:\ell_N^p\to\ell_N^q$ with $\rank R\le m$ such that
\begin{equation}\label{eq:critical-Gaussian-approximation}
	\|D_{\mathbf b}-R\|
	\le C\left[
	\frac{\|\mathbf b\|_{\ell_N^q}
		+N^{1/p'}\|\mathbf b\|_{\ell_N^\infty}}{\sqrt m}
	+\frac{N^{1/p'}\|\mathbf b\|_{\ell_N^q}}{m}
	\right].
\end{equation}
Write $D_{\mathbf b}=ZD_{|\mathbf b|}$, where $Z$ is diagonal
with entries $b_j/|b_j|$ when $b_j\ne0$ and $1$ otherwise.
Since $Z$ is an isometry on $\ell_N^q$,
\[
\|D_{\mathbf b}-ZR\|=\|D_{|\mathbf b|}-R\|,
\qquad \rank(ZR)=\rank R.
\]
Thus we may assume $\mathbf b\in[0,\infty)^N$.
It suffices to construct $R\in\mathbb R^{N\times N}$, since
its rank is the same over $\mathbb R$ and $\mathbb C$. Put
$A=D_{\mathbf b}-R$ and write $\mathbf x=\mathbf u+\mathrm i\mathbf v$
with $\mathbf u,\mathbf v\in\mathbb R^N$. Then
\[
\|A\mathbf x\|_{\ell_N^q}
\le \|A\mathbf u\|_{\ell_N^q}+\|A\mathbf v\|_{\ell_N^q}
\le 2\|A\|_{\ell_N^p(\mathbb R)\to\ell_N^q(\mathbb R)}
\|\mathbf x\|_{\ell_N^p}.
\]
The factor $2$ is absorbed into the constant in
\eqref{eq:critical-Gaussian-approximation}.

Let $G,H\in\mathbb R^{m\times N}$ be independent standard Gaussian
matrices. Write $\mathbb E_X$ for expectation with respect to $X$,
and $\mathbb E_{X,Y}$ for joint expectation in $X$ and $Y$.
Put $Q(G)=G^*G-mI_N$. Since $\mathbb E_HQ(H)=0$,
Jensen's inequality gives
\[
\mathbb E_G\|D_{\mathbf b}Q(G)\|
\le\mathbb E_{G,H}\|D_{\mathbf b}(G^*G-H^*H)\|.
\]
Set $U:=(G+H)/\sqrt2$ and $V:=(G-H)/\sqrt2$.
Then $U,V$ are independent standard Gaussian matrices and
$G^*G-H^*H=U^*V+V^*U$. By symmetry,
\begin{equation}\label{eq:gaussian-decoupling}
	\mathbb E_G\|D_{\mathbf b}Q(G)\|
	\le 2\mathbb E_{U,V}\|D_{\mathbf b}U^*V\|.
\end{equation}

Let $g_m\in\mathbb R^m$ and $h_N\in\mathbb R^N$
be independent standard Gaussian vectors, independent of $U,V$.
Fix $U$. Apply the real Gaussian Chevet inequality
\cite[inequality (1.1)]{ALLPT12} to the symmetric convex sets
$B_{\ell_N^p}$ and $UD_{\mathbf b}(B_{\ell_N^{q'}})$.
Its constant is absolute and independent of $m,N$; since $V$ is
independent of $U$, it gives
\[
\mathbb E_V\|D_{\mathbf b}U^*V\|
\le C\Bigl(
\|I_N\|_{\ell_N^p\to\ell_N^2}
\,\mathbb E_{g_m}\|D_{\mathbf b}U^*g_m\|_{\ell_N^q}
+\mathbb E_{h_N}\|h_N\|_{\ell_N^{p'}}
\,\|D_{\mathbf b}U^*\|_{\ell_m^2\to\ell_N^q}
\Bigr).
\]
Using $\|I_N\|_{\ell_N^p\to\ell_N^2}=1$ and
$\mathbb E_{h_N}\|h_N\|_{\ell_N^{p'}}\le CN^{1/p'}$, and then
taking expectation with respect to $U$, we obtain
\begin{equation}\label{eq:gaussian-product}
\mathbb E_{U,V}\|D_{\mathbf b}U^*V\|
\le C\left[
\mathbb E_{U,g_m}\|D_{\mathbf b}U^*g_m\|_{\ell_N^q}
+N^{1/p'}\mathbb E_U\|D_{\mathbf b}U^*\|_{\ell_m^2\to\ell_N^q}
\right].
\end{equation}
Gaussian moments give
$\mathbb E\|D_{\mathbf b}h_N\|_{\ell_N^q}\le C_q\|\mathbf b\|_{\ell_N^q}$.
Conditionally on $g_m$, the vector $U^*g_m$ has the distribution
of $\|g_m\|_{\ell_m^2}h_N$, so
\begin{equation}\label{eq:gaussian-vector-moment}
\mathbb E_{U,g_m}\|D_{\mathbf b}U^*g_m\|_{\ell_N^q}
\le C_q\sqrt m\,\|\mathbf b\|_{\ell_N^q}.
\end{equation}
Since $q'<2$, $\|D_{\mathbf b}\|_{\ell_N^{q'}\to\ell_N^2}
=\|\mathbf b\|_{\ell_N^\infty}$.
Apply the same Chevet inequality to $B_{\ell_m^2}$ and
$D_{\mathbf b}(B_{\ell_N^{q'}})$ to obtain
\[
\mathbb E_U\|D_{\mathbf b}U^*\|_{\ell_m^2\to\ell_N^q}
\le C\bigl(\sqrt m\,\|\mathbf b\|_{\ell_N^\infty}
+\mathbb E\|D_{\mathbf b}h_N\|_{\ell_N^q}\bigr).
\]
Hence
\begin{equation}\label{eq:gaussian-left-factor}
	\mathbb E_U\|D_{\mathbf b}U^*\|
	\le C\left(
	\sqrt m\,\|\mathbf b\|_{\ell_N^\infty}
	+\|\mathbf b\|_{\ell_N^q}
	\right).
\end{equation}
Combining \eqref{eq:gaussian-decoupling},
\eqref{eq:gaussian-product},
\eqref{eq:gaussian-vector-moment}, and
\eqref{eq:gaussian-left-factor}, we obtain
\[
\mathbb E_G\|D_{\mathbf b}-m^{-1}D_{\mathbf b}G^*G\|
\le C\left[
\frac{\|\mathbf b\|_{\ell_N^q}
	+N^{1/p'}\|\mathbf b\|_{\ell_N^\infty}}{\sqrt m}
+\frac{N^{1/p'}\|\mathbf b\|_{\ell_N^q}}{m}
\right].
\]
Choose $G$ with approximation error at most twice the bound above,
and set $R=m^{-1}D_{\mathbf b}G^*G$. Then $\rank R\le m$, proving
\eqref{eq:critical-Gaussian-approximation}.

Suppose that $p'>q$ and $\mathbf b$ is finitely supported. By taking
absolute values and permuting coordinates, we may assume
$b_1\ge b_2\ge\cdots\ge0$. Write $\theta=q/p'\in(0,1)$ and choose
$J\in\mathbb N$ with $(q/2)\theta^J\le1$. For $m\ge2$, set
$N_0=\lceil m^{q/2}\rceil$ and $N_{j+1}=\lceil N_j^\theta\rceil$
for $0\le j<J$. Since $\lceil t\rceil\le2t$ for $t\ge1$, iteration gives
\begin{equation}\label{eq:dimension-contraction}
N_j\le 2^{1+\theta+\cdots+\theta^j}m^{(q/2)\theta^j},
\qquad N_J\le Cm.
\end{equation}
Monotonicity gives $b_i\le i^{-1/q}\|\mathbf b\|_{\ell^q}$. Since
$N_0\ge m^{q/2}$,
\begin{equation}\label{eq:discarded-tail}
\sup_{i>N_0}b_i\le m^{-1/2}\|\mathbf b\|_{\ell^q}.
\end{equation}
For $0\le j<J$, write $I_j=\{N_{j+1}+1,\ldots,N_j\}$ and
$\mathbf b|_{I_j}=(b_i)_{i\in I_j}$. For each nonempty $I_j$,
we have $\|\mathbf b|_{I_j}\|_{\ell_{I_j}^q}\le\|\mathbf b\|_{\ell^q}$.
Using $N_{j+1}\ge N_j^{q/p'}$, we also have
\[
|I_j|^{1/p'}\|\mathbf b|_{I_j}\|_{\ell_{I_j}^\infty}
\le N_j^{1/p'}(N_{j+1}+1)^{-1/q}\|\mathbf b\|_{\ell^q}
\le \|\mathbf b\|_{\ell^q}.
\]
Moreover, $|I_j|^{1/p'}\le N_0^{1/p'}\le Cm^{q/(2p')}\le C\sqrt m$.
By \eqref{eq:critical-Gaussian-approximation}, there is an operator
$R_j:\ell_{I_j}^p\to\ell_{I_j}^q$ of rank at most $m$ such that
\begin{equation}\label{eq:critical-weak-block-error}
\|D_{\mathbf b|_{I_j}}-R_j\|\le Cm^{-1/2}\|\mathbf b\|_{\ell^q}.
\end{equation}

Define $R:\ell^p\to\ell^q$ to agree with $D_{\mathbf b}$ on the
first $N_J$ coordinates, act as $R_j$ on each nonempty block $I_j$,
and vanish on coordinates beyond $N_0$. By
\eqref{eq:dimension-contraction}, $\rank R\le N_J+Jm\le Km$
for a fixed integer $K$. Since $p\le q$, the norm of the
block-diagonal error is bounded by the largest block norm.
Thus \eqref{eq:discarded-tail} and \eqref{eq:critical-weak-block-error}
give
\[
a_{Km+1}(D_{\mathbf b})
\le\|D_{\mathbf b}-R\|\le Cm^{-1/2}\|\mathbf b\|_{\ell^q}.
\]
For $m=1$, take $R=0$ and use
$\|D_{\mathbf b}\|=\|\mathbf b\|_{\ell^\infty}\le\|\mathbf b\|_{\ell^q}$.

If $p'<q$ and $\mathbf b$ is finitely supported, choose $N$ such that
$\operatorname{supp}\mathbf b\subseteq\{1,\ldots,N\}$.
Apply the preceding estimate to
$D_{\mathbf b}^*:\ell_N^{q'}\to\ell_N^{p'}$.
Indeed, $1<q'<2<p'$ and $(q')'=q>p'$, so the preceding case applies.
In finite dimensions, taking adjoints preserves rank and operator norm;
applying this in both directions gives
$a_k(D_{\mathbf b})=a_k(D_{\mathbf b}^*)$.
Thus \eqref{eq:critical-weak-approximation} holds for finitely supported
$\mathbf b$ in both cases. For arbitrary $\mathbf b\in\ell^\tau$, put
$\mathbf b^{(N)}=(b_1,\ldots,b_N,0,\ldots)$. Since
\[
\|D_{\mathbf b}-D_{\mathbf b^{(N)}}\|
=\sup_{i>N}|b_i|\longrightarrow0,
\qquad
\|\mathbf b^{(N)}\|_{\ell^\tau}\le\|\mathbf b\|_{\ell^\tau},
\]
Lemma~\ref{lem:ideal-properties}\textup{(iv)} allows us to pass
to the limit in \eqref{eq:critical-weak-approximation}.

Finally, for $0<p\le1$, choose $1<p_0<q'$. The inclusion
$I:\ell^p\to\ell^{p_0}$ has norm $1$ and
$D_{\mathbf b}=\widetilde D_{\mathbf b}I$, where
$\widetilde D_{\mathbf b}:\ell^{p_0}\to\ell^q$ has the same
diagonal entries. Since $p_0'>q$ and $\tau=q$,
Lemma~\ref{lem:ideal-properties}\textup{(i)} gives
\[
a_{Km+1}(D_{\mathbf b})
\le a_{Km+1}(\widetilde D_{\mathbf b})
\le Cm^{-1/2}\|\mathbf b\|_{\ell^\tau}.
\]
All constants are independent of $\mathbf b$, $m$, and $N$.
This completes the proof.
\end{proof}

\subsection{Sequence criteria for diagonal operators}

\begin{definition}\label{def:sequence-Lorentz}
	Let $\mathbf b=(b_j)_{j\in\mathbb N}$ be a complex-valued sequence.
	Its distribution function and non-increasing rearrangement
	are defined by
	\begin{equation}\label{eq:discrete-rearrangement}
		d_{\mathbf b}(u):=\#\{j\in\mathbb N:|b_j|>u\},\qquad
		b_k^*:=\inf\{u\ge0:d_{\mathbf b}(u)<k\},\quad k\in\mathbb N.
	\end{equation}
	For $0<\rho,s<\infty$, the Lorentz sequence space
	$\ell^{\rho,s}$ consists of the sequences $\mathbf b$ such that
	\begin{equation}\label{eq:sequence-Lorentz}
		\|\mathbf b\|_{\ell^{\rho,s}}
		:=\left(\sum_{k\in\mathbb N}
		k^{s/\rho-1}(b_k^*)^s\right)^{1/s}<\infty.
	\end{equation}
	For $\beta\in\mathbb R$, the Lorentz--Zygmund sequence space
	$\ell^{\rho,s,\beta}$ consists of the sequences $\mathbf b$ such that
	\begin{equation}\label{eq:sequence-Lorentz-Zygmund}
		\|\mathbf b\|_{\ell^{\rho,s,\beta}}
		:=\left(\sum_{k\in\mathbb N}
		k^{s/\rho-1}(1+\log k)^{\beta s}(b_k^*)^s\right)^{1/s}
		<\infty.
	\end{equation}
	In particular, $\ell^{\rho,s,0}=\ell^{\rho,s}$ and $\ell^{s,s}=\ell^s$ with
	\[
	\|\mathbf b\|_{\ell^{s,s}}
	=\left(\sum_{k\in\mathbb N}(b_k^*)^s\right)^{1/s}
	=\left(\sum_{j\in\mathbb N}|b_j|^s\right)^{1/s}
	=\|\mathbf b\|_{\ell^s}.
	\]
\end{definition}

Let $J$ be a nonempty subset of the coordinate index set.
Define $P_J:\ell^u\to\ell_J^u$ and
$\iota_J:\ell_J^u\to\ell^u$ by
\begin{align}\label{Pi}
	P_J\mathbf x=(x_j)_{j\in J},\qquad
	(\iota_J\mathbf y)_j=
	\begin{cases}
		y_j,&j\in J,\\
		0,&j\notin J.
	\end{cases}
\end{align}
Then $D_{\mathbf b|_J}=P_JD_{\mathbf b}\iota_J$, where $\mathbf b|_J:=(b_j)_{j\in J}$.

To describe the critical case $s=2$, we use the discrete counterpart
of the Copson--Lorentz tail integral. For $2<\nu<\infty$, set
\begin{equation}\label{eq:critical-sequence-tail}
	\mathcal T_\nu(\mathbf b)
	:=\left[
	\sum_{k\ge1}\frac1k
	\left(\sum_{j\ge k}(b_j^*)^\nu\right)^{2/\nu}
	\right]^{1/2}.
\end{equation}

The following lemmas characterize approximation ideal membership
for diagonal operators in terms of their diagonal sequences.

\begin{lemma}\label{lem:same-side}
	Let $0<p\le q<\infty$, $q\ge1$, and $0<s<\infty$.
	Let $\mathbf b\in\ell^\infty$ be nonnegative and
	$D_{\mathbf b}:\ell^p\to\ell^q$ be the corresponding diagonal operator.
	\begin{enumerate}
		\item[(i)]
		If $(p,q,s)\notin\mathcal E$, then
		$D_{\mathbf b}\in\mathfrak A_s(\ell^p,\ell^q)$ if and only if
		$\mathbf b\in\ell^{\rho,s}$. Moreover,
		\begin{equation}\label{eq:noncritical-diagonal}
			\As{D_{\mathbf b}}\simeq\|\mathbf b\|_{\ell^{\rho,s}}.
		\end{equation}
		
		\item[(ii)]
		If $(p,q,s)\in\mathcal E$, then
		$D_{\mathbf b}\in\mathfrak A_2(\ell^p,\ell^q)$ if and only if
		the corresponding expression on the right-hand side of
		\eqref{eq:critical-diagonal-sufficiency} is finite, and
		\begin{equation}\label{eq:critical-diagonal-sufficiency}
			\|D_{\mathbf b}\|_{\mathfrak A_2}
			\simeq
			\begin{cases}
				\|\mathbf b\|_{\ell^{q,2,1/2}},&p^\dagger=q,\\[3pt]
				\mathcal T_\tau(\mathbf b),&p^\dagger\ne q.
			\end{cases}
		\end{equation}
	\end{enumerate}
	The comparison constants are independent of $\mathbf b$.
\end{lemma}

\begin{remark}
	If $q\le2$ or $p\ge2$, then $\rho=s$ and
	\begin{equation}\label{eq:same-side}
		2^{-1/2}b_{2k}^*\le a_k(D_{\mathbf b})\le b_k^*,
		\qquad k\in\mathbb N.
	\end{equation}
	Consequently, $\As{D_{\mathbf b}}\simeq\|\mathbf b\|_{\ell^s}$.
\end{remark}

\begin{proof}[Proof of Lemma~\ref{lem:same-side}]
	For $J\subset\mathbb N$, write
	$B_J=P_JD_{\mathbf b}\iota_J=D_{\mathbf b|_J}$.
	Since $P_J$ and $\iota_J$ are contractions,
	Lemma~\ref{lem:ideal-properties}\textup{(i)} gives
	$a_k(B_J)\le a_k(D_{\mathbf b})$.
	For finite $J$, the operator $\iota_JB_JP_J$ has rank at most $|J|$.
	By Lemma~\ref{lem:diag-facts}\textup{(i)},
	$\|D_{\mathbf b}-\iota_JB_JP_J\|=\sup_{j\notin J}b_j$.
	Hence
	\begin{equation}\label{eq:diagonal-truncation}
		a_k(D_{\mathbf b})\le\inf_{|J|<k}\sup_{j\notin J}b_j=b_k^*.
	\end{equation}
	
	\textup{(i)}
	Suppose first that $q\le2$ or $p\ge2$. If $b_{2k}^*>0$,
	fix $0<\varepsilon<b_{2k}^*$ and choose $J$ with $|J|=2k$
	such that $b_j\ge b_{2k}^*-\varepsilon$ for $j\in J$.
	For $S:\ell_J^p\to\ell_J^q$ with $\rank S<k$, put $T=B_J-S$.
	Choose a Euclidean orthonormal basis $(x_\nu)_{\nu=1}^{2k}$
	with $x_1,\ldots,x_{k+1}\in\ker S$. Then
	\begin{equation}\label{eq:diagonal-HS}
		\|T\|_{\mathrm{HS}}^2
		=\sum_{\nu=1}^{2k}\|Tx_\nu\|_{\ell_J^2}^2
		\ge\sum_{\nu=1}^{k+1}\|B_Jx_\nu\|_{\ell_J^2}^2
		\ge(k+1)(b_{2k}^*-\varepsilon)^2.
	\end{equation}
	If $q\le2$, then
	\begin{equation}\label{eq:frobenius}
		\|T\|\ge\max_{j\in J}\|Te_j\|_{\ell_J^q}
		\ge(2k)^{-1/2}\|T\|_{\mathrm{HS}}
		\ge2^{-1/2}(b_{2k}^*-\varepsilon).
	\end{equation}
	If $p\ge2$, then $q'\le p'\le2$, and the same column estimate
	for $T^*:\ell_J^{q'}\to\ell_J^{p'}$ gives
	\begin{equation}\label{eq:diagonal-duality}
		\|T\|=\|T^*\|\ge(2k)^{-1/2}\|T^*\|_{\mathrm{HS}}
		\ge2^{-1/2}(b_{2k}^*-\varepsilon).
	\end{equation}
	Taking the infimum over $S$ and letting $\varepsilon\downarrow0$
	proves the lower bound in \eqref{eq:same-side}. When $b_{2k}^*=0$,
	the bound is immediate. Since $a_1(D_{\mathbf b})=b_1^*$ and $\rho=s$,
	\begin{equation}\label{eq:rearrange-pair}
		\|\mathbf b\|_{\ell^{s,s}}^s
		\le(b_1^*)^s+2\sum_{k\ge1}(b_{2k}^*)^s
		\le(1+2^{1+s/2})\sum_{k\ge1}a_k(D_{\mathbf b})^s.
	\end{equation}
	Together with \eqref{eq:diagonal-truncation}, this proves
	\eqref{eq:noncritical-diagonal} in this case.
	
	Suppose now that $0<p<2<q<\infty$. For finite $J$ with $|J|=m$
	and $d_J=\min_{j\in J}b_j>0$, the identity on $\ell_J^p$
	viewed as a map into $\ell_J^q$ factors as
	$I_m=D_{(b_j^{-1})_{j\in J}}B_J$, where the first factor
	has norm $d_J^{-1}$ on $\ell_J^q$.
	By Lemma~\ref{lem:ideal-properties}\textup{(i)} and
	\eqref{eq:finite-crossing-width},
	\begin{equation}\label{eq:finite-prefix-lower}
		a_k(D_{\mathbf b})\ge d_Ja_k(I_m)
		\ge cd_J\min\{1,m^{1/\tau}k^{-1/2}\},\qquad m\ge2k.
	\end{equation}
	If $\#\{j:b_j\ge d\}=\infty$ for some $d>0$, we can choose
	$J\subset\{j:b_j\ge d\}$ with arbitrarily large cardinality.
	Letting $m\to\infty$ in \eqref{eq:finite-prefix-lower} gives
	\begin{equation}\label{eq:diagonal-compact-coordinates}
		a_k(D_{\mathbf b})\ge cd,\qquad k\ge1.
	\end{equation}
	Thus $D_{\mathbf b}\in\mathfrak A_s$ implies $\mathbf b\in c_0$.
	For $J=\{j:b_j>0\}$, the identities
	\[
	D_{\mathbf b}=\iota_JB_JP_J,\qquad B_J=P_JD_{\mathbf b}\iota_J
	\]
	and the ideal property give
	\begin{equation}\label{eq:diagonal-active-coordinates}
		a_k(D_{\mathbf b})=a_k(B_J).
	\end{equation}
	We may therefore list the positive entries in non-increasing order
	and append zeros if there are only finitely many. Taking the
	$m$ largest entries in \eqref{eq:finite-prefix-lower} gives
	\begin{equation}\label{eq:prefix-crossing-lower}
		a_k(D_{\mathbf b})
		\ge cb_m\min\{1,m^{1/\tau}k^{-1/2}\},\qquad m\ge2k.
	\end{equation}
	If $b_m=0$, the same bound is immediate.
	
	If $0<s<2$, then $s/\rho-1=-s\eta_{p,q}$.
	Taking $m=2k$ in \eqref{eq:prefix-crossing-lower} gives
	\begin{equation}\label{eq:crossing-small-s-lower}
		\|\mathbf b\|_{\ell^{\rho,s}}^s
		=\sum_{j\ge1}j^{-s\eta_{p,q}}b_j^s
		\le C\left(b_1^s+\sum_{k\ge1}k^{-s\eta_{p,q}}b_{2k}^s\right)
		\le C\sum_{k\ge1}a_k(D_{\mathbf b})^s.
	\end{equation}
	If $s>2$, put $\alpha:=2/\tau\in (0,1)$ and $m_k=\lceil4k^{1/\alpha}\rceil$.
	By \eqref{eq:prefix-crossing-lower},
	\begin{equation}\label{eq:crossing-large-s-prefix}
		m_k\ge2k,\qquad m_k^{1/\tau}k^{-1/2}\ge1,\qquad
		a_k(D_{\mathbf b})\ge cb_{m_k}.
	\end{equation}
	Since $m_{k+1}^{\alpha}-m_k^\alpha\le4^\alpha+1$,
	\begin{equation}\label{eq:crossing-large-s-interval}
		\sum_{j=m_k}^{m_{k+1}-1}j^{\alpha-1}
		\le m_k^{\alpha-1}+\int_{m_k}^{m_{k+1}}x^{\alpha-1}\,dx
		\le 1+\frac{4^\alpha+1}{\alpha}.
	\end{equation}
	As $s/\rho-1=\alpha-1$ and $m_1=4$,
	\eqref{eq:crossing-large-s-prefix} and
	\eqref{eq:crossing-large-s-interval} yield
\begin{align}\label{eq:crossing-large-s-lower}
			\|\mathbf b\|_{\ell^{\rho,s}}^s
			\le Cb_1^s+\sum_{k\ge1}b_{m_k}^s
			\sum_{j=m_k}^{m_{k+1}-1}j^{\alpha-1}
			\le C\left(b_1^s+\sum_{k\ge1}b_{m_k}^s\right)
			\le C\sum_{k\ge1}a_k(D_{\mathbf b})^s.
		\end{align}

	For the reverse estimate, suppose $\mathbf b\in\ell^{\rho,s}$.
	Since $s/\rho>0$,
	\begin{equation}\label{eq:lorentz-coordinate-decay}
		m^{s/\rho}(b_m^*)^s
		\le C\sum_{k=1}^m k^{s/\rho-1}(b_k^*)^s
		\le C\|\mathbf b\|_{\ell^{\rho,s}}^s.
	\end{equation}
	Hence $\mathbf b\in c_0$, and \eqref{eq:diagonal-active-coordinates}
	again permits non-increasing order. Set
	\begin{equation}\label{eq:dyadic-crossing-blocks}
		I_j=\{2^j,\ldots,2^{j+1}-1\},\qquad
		T_j=D_{\mathbf b|_{I_j}},\qquad T^{(L)}=\bigoplus_{j=0}^LT_j,
	\end{equation}
	where $T_j:\ell_{I_j}^p\to\ell_{I_j}^q$ and $T^{(L)}$ is extended
	by zero to $\ell^p\to\ell^q$. By diagonal domination and
	\eqref{eq:finite-crossing-sum},
	\begin{equation}\label{eq:dyadic-crossing-cost}
		\sum_{k\ge1}a_k(T_j)^s
		\le b_{2^j}^s\sum_{k=1}^{2^j}a_k(I_{2^j})^s
		\le C2^{js/\rho}b_{2^j}^s.
	\end{equation}
	For $\varepsilon>0$, write
	$N_T(\varepsilon)=\#\{k\ge1:a_k(T)>\varepsilon\}$.
	Since $a_{N_{T_j}(\varepsilon)+1}(T_j)\le\varepsilon$,
	choose $R_j:\ell_{I_j}^p\to\ell_{I_j}^q$ with
	\begin{equation}\label{eq:block-approximants}
		\rank R_j\le N_{T_j}(\varepsilon),\qquad
		\|T_j-R_j\|<2\varepsilon,\qquad 0\le j\le L.
	\end{equation}
	Put $R^{(L)}=\bigoplus_{j=0}^LR_j$ and extend it by zero.
	By disjointness of the blocks and $p\le q$,
	\begin{equation}\label{eq:block-error}
		\|(T^{(L)}-R^{(L)})\mathbf x\|_{\ell^q}^q
		\le(2\varepsilon)^q\sum_{j=0}^L
		\|\mathbf x|_{I_j}\|_{\ell_{I_j}^p}^q
		\le(2\varepsilon)^q\|\mathbf x\|_{\ell^p}^q.
	\end{equation}
	Thus $a_{\rank R^{(L)}+1}(T^{(L)})\le2\varepsilon$, so
	\begin{equation}\label{eq:block-counting}
		N_{T^{(L)}}(2\varepsilon)
		\le\rank R^{(L)}\le\sum_{j=0}^L N_{T_j}(\varepsilon).
	\end{equation}
	By Tonelli's theorem,
	\begin{equation}\label{eq:approximation-layer-cake}
		\sum_{k\ge1}a_k(T)^s
		=s\int_0^\infty\varepsilon^{s-1}N_T(\varepsilon)\,d\varepsilon.
	\end{equation}
	Integrating \eqref{eq:block-counting} with this formula gives
	\begin{equation}\label{eq:finite-block-summability}
		\sum_{k\ge1}a_k(T^{(L)})^s
		\le2^s\sum_{j=0}^L\sum_{k\ge1}a_k(T_j)^s.
	\end{equation}
	Summing over $2^j\le k<2^{j+1}$ and using
	$b_{2^{j+1}}\le b_k\le b_{2^j}$ gives
	\begin{equation}\label{eq:dyadic-weight-comparison}
		\sum_{j\ge0}2^{js/\rho}b_{2^j}^s
		\simeq\sum_{k\ge1}k^{s/\rho-1}b_k^s.
	\end{equation}
	Combining \eqref{eq:dyadic-crossing-cost},
	\eqref{eq:finite-block-summability}, and
	\eqref{eq:dyadic-weight-comparison}, we obtain
	\begin{equation}\label{eq:dyadic-crossing-upper}
		\sum_{k\ge1}a_k(T^{(L)})^s\le C\|\mathbf b\|_{\ell^{\rho,s}}^s.
	\end{equation}
	By Lemmas~\ref{lem:diag-facts}\textup{(i)} and
	\ref{lem:ideal-properties}\textup{(iv)},
	\begin{equation}\label{eq:diagonal-block-convergence}
		\|D_{\mathbf b}-T^{(L)}\|=b_{2^{L+1}}\longrightarrow0,
		\qquad a_k(T^{(L)})\longrightarrow a_k(D_{\mathbf b}).
	\end{equation}
	Fatou's lemma and \eqref{eq:dyadic-crossing-upper} give
	\begin{equation}\label{eq:noncritical-diagonal-limit}
		\sum_{k\ge1}a_k(D_{\mathbf b})^s
		\le\liminf_{L\to\infty}\sum_{k\ge1}a_k(T^{(L)})^s
		\le C\|\mathbf b\|_{\ell^{\rho,s}}^s.
	\end{equation}
	Together with \eqref{eq:crossing-small-s-lower} for $0<s<2$
	and \eqref{eq:crossing-large-s-lower} for $s>2$,
	this proves \textup{(i)}.
	
\textup{(ii)}
We first prove a lower bound for $D_{\mathbf b}:\ell_N^p\to\ell_N^q$
with $p>0$ and $q>2$. Assume $b_1\ge\cdots\ge b_N\ge0$ and
extend $\mathbf b$ by zeros. Let $\varepsilon>0$ and let $R$ satisfy
$\rank R<k$ and $\|D_{\mathbf b}-R\|\le\varepsilon$.
Testing on the coordinate vectors shows that each column of
$D_{\mathbf b}-R$ has $\ell_N^q$ norm at most $\varepsilon$.
Set $I:=\{j:b_j>2\varepsilon\}$ and
\[
S:=\sum_{\substack{j\in I\\j>64k}}b_j^q,\qquad
c_j:=\min\{b_j^q,S/(64k)\}\quad(j\in I),\qquad
B:=\sum_{j\in I}c_j.
\]
Suppose $S>0$. Then $64k\in I$. If $b_{64k}^q\ge S/(64k)$,
the first $64k$ terms give $B\ge S$. Otherwise, $c_j=b_j^q$
for $j\in I$ with $j>64k$, so again $B\ge S$.

Put $E:=(D_{\mathbf b}-R)|_{I\times I}$ and define
\[
E'_{ij}:=\frac{\sqrt{c_ic_j}}{b_i}E_{ij},\qquad
R'_{ij}:=\frac{\sqrt{c_ic_j}}{b_i}R_{ij},\qquad i,j\in I.
\]
Writing $\mathbf c:=(c_j)_{j\in I}$, we have $R'=D_{\mathbf c}-E'$
and $\rank R'<k$. Since $|E_{ii}|\le\varepsilon<b_i/2$ on $I$,
\[
\operatorname{Re}\operatorname{tr}R'
=\sum_{i\in I}c_i\operatorname{Re}(1-E_{ii}/b_i)\ge B/2.
\]
As $c_i\le b_i^q$, we have $(c_i/b_i^2)^{q/(q-2)}\le c_i$.
H\"older's inequality applied to each column gives
\[
\|E'\|_{\mathrm{HS}}^2
=\sum_{j\in I}c_j\sum_{i\in I}\frac{c_i}{b_i^2}|E_{ij}|^2
\le\varepsilon^2 B^{2-2/q}.
\]
Also, $\|D_{\mathbf c}\|_{\mathrm{HS}}^2\le SB/(64k)\le B^2/(64k)$.
Thus the trace bound
$|\operatorname{tr}R'|\le\sqrt{\rank R'}\,\|R'\|_{\mathrm{HS}}$
gives $B/2\le B/8+\sqrt{k}\,\varepsilon B^{1-1/q}$.
Since $S\le B$, we obtain
\begin{equation}\label{eq:critical-truncated-tail}
	S=\sum_{\substack{j\in I\\j>64k}}b_j^q
	\le Ck^{q/2}\varepsilon^q.
\end{equation}
which also holds when $S=0$.
For $\varepsilon>a_k(D_{\mathbf b})$, choose $R$ as above.
The terms with $b_j\le2\varepsilon$ and $j\le k^{q/2}$ contribute
at most $2^qk^{q/2}\varepsilon^q$. Letting
$\varepsilon\downarrow a_k(D_{\mathbf b})$ yields
\begin{equation}\label{eq:critical-tail-window}
	\left(\sum_{64k<j\le k^{q/2}}b_j^q\right)^{2/q}
	\le Ck\,a_k(D_{\mathbf b})^2.
\end{equation}

Set $F(t):=b_{\lfloor t\rfloor}$ for $t\ge1$ and fix $1<a<q/2$.
Choose an integer $K$ so large that, for $t\ge K$ and
$k:=\lfloor t/128\rfloor$, we have $k\ge1$,
$64k<\lfloor t\rfloor$, and $\lceil t^a\rceil\le k^{q/2}$.
By \eqref{eq:critical-tail-window},
$(\int_t^{t^a}F(x)^q\,dx)^{2/q}\le Ck\,a_k(D_{\mathbf b})^2$.
Integrating over $[128k,128(k+1))\cap[K,\infty)$ and using
$\int_{128k}^{128(k+1)}dt/t\le1/k$, we obtain
\begin{equation}\label{eq:critical-tail-window-integral}
	\int_K^\infty
	\left(\int_t^{t^a}F(x)^q\,dx\right)^{2/q}\frac{dt}{t}
	\le C\sum_{k\ge1}a_k(D_{\mathbf b})^2.
\end{equation}
Since $2/q<1$, splitting $(t,\infty)$ into the intervals $(t^{a^j},t^{a^{j+1}})$
for $j\in\mathbb N\cup\{0\}$
and substituting $y=t^{a^j}$ in each term gives
\begin{equation}\label{eq:critical-full-tail-integral}
	\int_K^\infty
	\left(\int_t^\infty F(x)^q\,dx\right)^{2/q}\frac{dt}{t}
	\le\sum_{j\ge0}a^{-j}
	\int_{K^{a^j}}^\infty
	\left(\int_y^{y^a}F(x)^q\,dx\right)^{2/q}\frac{dy}{y}
	\le C\sum_{k\ge1}a_k(D_{\mathbf b})^2.
\end{equation}
For each $k>K$, the contribution of $[k-1,k)$ to the integral
on the left of \eqref{eq:critical-full-tail-integral} is at least
$k^{-1}(\sum_{j\ge k}b_j^q)^{2/q}$. Moreover, \[
\sum_{k=1}^{K}\frac1k\left(\sum_{j\ge k}b_j^q\right)^{2/q}
\le \left(\sum_{k=1}^{K}\frac1k\right)\|\mathbf b\|_{\ell^q}^2
\le C\|\mathbf b\|_{\ell^q}^2,
\]
where $C$ depends only on $K$.
Splitting the $\ell^q$ sum at $2K$ and using $b_1\le a_1(D_{\mathbf b})$
gives
\[
\|\mathbf b\|_{\ell^q}^2
\le C\left[
a_1(D_{\mathbf b})^2+
\int_K^{2K}\left(\int_t^\infty F(x)^q\,dx\right)^{2/q}\frac{dt}{t}
\right].
\]
Consequently,
\begin{equation}\label{eq:critical-tail-necessity}
	\mathcal T_q(\mathbf b)\le C\|D_{\mathbf b}\|_{\mathfrak A_2},
\end{equation}
with $C$ independent of $N$ and $\mathbf b$.

If $p^\dagger\ge q$, then $\tau=q$ and
\eqref{eq:critical-tail-necessity} applies directly. If $p^\dagger<q$,
then $\tau=p^\dagger$ and we apply it to
$D_{\mathbf b}^*:\ell_N^{q'}\to\ell_N^{p^\dagger}$, using
$a_k(D_{\mathbf b}^*)=a_k(D_{\mathbf b})$. Thus
\begin{equation}\label{eq:critical-tail-lower}
	\mathcal T_\tau(\mathbf b)\le C\|D_{\mathbf b}\|_{\mathfrak A_2}.
\end{equation}

Suppose $p^\dagger\ne q$. Write
$D_{\mathbf b}=D_{(b_1,\ldots,b_m)}\oplus D_{(b_j)_{j>m}}$.
Since $\rank D_{(b_1,\ldots,b_m)}\le m$,
Lemma~\ref{lem:critical-weak-approximation} gives, for a fixed integer $L\ge2$,
\begin{equation}\label{eq:critical-tail-upper}
	a_{Lm+1}(D_{\mathbf b})
	\le a_{(L-1)m+1}(D_{(b_j)_{j>m}})
	\le Cm^{-1/2}\left(\sum_{j>m}b_j^\tau\right)^{1/\tau},
	\qquad m\ge1.
\end{equation}
For $Lm<k\le L(m+1)$, we have $a_k(D_{\mathbf b})\le a_{Lm+1}(D_{\mathbf b})$,
and for $k\le L$ we have $a_k(D_{\mathbf b})\le b_1$. Summing their squares
using \eqref{eq:critical-tail-upper} and enlarging each tail from $j>m$ to $j\ge m$ gives
\[
\|D_{\mathbf b}\|_{\mathfrak A_2}^2
\le C\left[b_1^2+
\sum_{m\ge1}\frac1m\left(\sum_{j>m}b_j^\tau\right)^{2/\tau}\right]
\le C\mathcal T_\tau(\mathbf b)^2.
\]
Together with \eqref{eq:critical-tail-lower}, this proves the finite
sequence estimate when $p^\dagger\ne q$.

Suppose $p^\dagger=q$, so $p=q'$, and put $\eta:=1/2-1/q$.
For $j\ge0$, put
$n_j:=2^j$, $I_j:=\{n_j,\ldots,2n_j-1\}$, and $c_j:=b_{2n_j-1}$.
Set $B:=\bigoplus_{j\ge0}c_jI_{n_j}$ and work on its nonzero blocks.
Since $c_j\le b_i$ for $i\in I_j$, diagonal domination gives
$\|B\|_{\mathfrak A_2}\le\|D_{\mathbf b}\|_{\mathfrak A_2}$.
For $\varepsilon>0$ and $\|B-R\|\le\varepsilon$, put $E:=B-R$.
Testing $E$ and $E^*$ on the coordinate vectors shows that every row
and column of $E$ has $\ell^q$ norm at most $\varepsilon$.
Hence the block $E_{jk}:\C^{I_k}\to\C^{I_j}$ satisfies
\begin{equation}\label{eq:critical-selfdual-block-error}
	\|E_{jk}\|_{\mathrm{HS}}^2
	\le\varepsilon^2
	\min\{n_kn_j^{1-2/q},n_jn_k^{1-2/q}\}.
\end{equation}
Set $J_\varepsilon:=\{j:c_jn_j^{-\eta}\le\varepsilon<c_j/2\}$
and $v_j:=c_jn_j^{1/q}\chi_{J_\varepsilon}(j)$.
Define the diagonal matrix $W$ by
$W|_{\C^{I_j}}:=v_j^{1/2}n_j^{-1/2+1/(2q)}I_{n_j}$.
By \eqref{eq:critical-selfdual-block-error},
\begin{equation}\label{eq:critical-weighted-error}
	\|WEW\|_{\mathrm{HS}}^2
	\le\varepsilon^2\sum_{j,k}v_jv_k2^{-|j-k|/q}
	\le C_q\varepsilon^2\sum_jv_j^2.
\end{equation}
The last inequality follows from $2v_jv_k\le v_j^2+v_k^2$
and $\sum_{\ell\in\mathbb Z}2^{-|\ell|/q}<\infty$.
Also, $\|WBW\|_{\mathrm{HS}}^2
=\sum_jv_j^2c_j^2n_j^{-2\eta}\le\varepsilon^2\sum_jv_j^2$.
Since $|E_{ii}|\le\varepsilon<c_j/2$ on the selected blocks,
$\operatorname{Re}\operatorname{tr}(WRW)\ge\frac12\sum_jv_j^2$.
Using $|\operatorname{tr}A|\le\sqrt{\rank A}\|A\|_{\mathrm{HS}}$
and \eqref{eq:critical-weighted-error} gives
\[
\frac12\sum_jv_j^2
\le\operatorname{Re}\operatorname{tr}(WRW)
\le C_q\varepsilon\sqrt{\rank R}
\left(\sum_jv_j^2\right)^{1/2}.
\]
Therefore,
\begin{equation}\label{eq:critical-selfdual-rank}
	\rank R\ge c_q\varepsilon^{-2}
	\sum_{j\in J_\varepsilon}c_j^2n_j^{2/q}.
\end{equation}
The same argument applies to any finite collection of dyadic blocks
with nonnegative coefficients, with the constant depending only on $q$.
For $\varepsilon>0$, finite-dimensional compactness gives
$N_B(\varepsilon)=\min\{\rank R:\|B-R\|\le\varepsilon\}$.
Integrating \eqref{eq:critical-selfdual-rank} in
\eqref{eq:approximation-layer-cake} with $s=2$ yields
\[
\|B\|_{\mathfrak A_2}^2
\gtrsim_q\sum_{j\ge0}c_j^2n_j^{2/q}
(\eta\log n_j-\log2)_+,
\]
since each nonempty active interval has logarithmic length
$\eta\log n_j-\log2$.
Since $n_j=2^j$, the last factor is comparable to $j+1$ for
$j\ge2/\eta$. The remaining finitely many terms are bounded
by $C_qb_1^2\le C_q\|D_{\mathbf b}\|_{\mathfrak A_2}^2$.
Using $b_{2^j}\le c_{j-1}$ for $j\ge1$, dyadic summation gives
\begin{equation}\label{eq:critical-selfdual-lower}
	\|\mathbf b\|_{\ell^{q,2,1/2}}^2
	\le C\left[b_1^2+\sum_{j\ge0}(j+1)2^{2j/q}c_j^2\right]
	\le C\|D_{\mathbf b}\|_{\mathfrak A_2}^2.
\end{equation}
Conversely, \eqref{eq:finite-crossing-sum} and
\eqref{eq:finite-block-summability} with $s=2$ give
\begin{equation}\label{eq:critical-finite-block-upper}
	\|D_{\mathbf b}\|_{\mathfrak A_2}^2
	\le C\sum_{j\ge0}(j+1)2^{2j/q}b_{2^j}^2
	\le C\|\mathbf b\|_{\ell^{q,2,1/2}}^2.
\end{equation}
The last inequality follows by summing over $2^{j-1}\le k<2^j$
for $j\ge1$. The $j=0$ term is $b_1^2$.

Finally, if $\|D_{\mathbf b}\|_{\mathfrak A_2}<\infty$, then
$\mathbf b\in c_0$ by \eqref{eq:diagonal-compact-coordinates}.
Finiteness of either expression on the right of
\eqref{eq:critical-diagonal-sufficiency} also implies $\mathbf b\in c_0$.
Thus, whenever one side is finite, we may arrange the positive entries
in non-increasing order, appending zeros if necessary.

Set $\mathbf b^{(N)}:=(b_1,\ldots,b_N,0,\ldots)$ and
$D_{\mathbf b}^{(N)}:=D_{\mathbf b^{(N)}}$. By the ideal property,
taking the first $N$ coordinates and extending by zero preserve the
approximation numbers of the finite diagonal block. Thus the
finite-dimensional estimates apply to $D_{\mathbf b}^{(N)}$ with
constants independent of $N$.

Let $Q_N$ be the projection that keeps the first $N$ coordinates and
sets the rest to zero. Since $D_{\mathbf b}^{(N)}=Q_ND_{\mathbf b}^{(N+1)}$
and $\|Q_N\|\le1$, the numbers $a_k(D_{\mathbf b}^{(N)})$ increase
with $N$. Moreover,
\[
\bigl|a_k(D_{\mathbf b}^{(N)})-a_k(D_{\mathbf b})\bigr|
\le \|D_{\mathbf b}^{(N)}-D_{\mathbf b}\|
=b_{N+1}\longrightarrow0.
\]
The squared sequence expressions in
\eqref{eq:critical-diagonal-sufficiency} also increase to their
values for $\mathbf b$, since each truncation adds a nonnegative
term. Passing to the limit by monotone convergence gives both
estimates for $D_{\mathbf b}$.
This completes the proof.
\end{proof}

\begin{remark}\label{rem:critical-separation}
	The critical conditions
	$\mathcal T_\tau(\mathbf b)<\infty$ and
	$\mathbf b\in\ell^{\tau,2,1/2}$ are not equivalent.
	Indeed, fix $\tau>2$ and $\alpha>0$, and let
	\[
	\mathbf b=(b_j)_{j\in \mathbb N},
	\qquad
	b_j:=j^{-1/\tau}[\log(ej)]^{-\alpha}.
	\]
	Since $b_j^\tau=j^{-1}[\log(ej)]^{-\alpha\tau}$, the integral test
	reduces the convergence of $\sum_{j\ge k}b_j^\tau$ to that of
	\[
	\int_k^\infty\frac{dx}{x[\log(ex)]^{\alpha\tau}}
	=\int_{\log(ek)}^\infty u^{-\alpha\tau}\,du.
	\]
	If $\alpha\tau>1$, integral comparison gives
	$\sum_{j\ge k}b_j^\tau\simeq[\log(ek)]^{1-\alpha\tau}$.
	If $\alpha\tau\le1$, the integral diverges for every $k\ge1$,
	so $\sum_{j\ge k}b_j^\tau=\infty$ for every $k\ge1$. Hence
	\[
	\mathcal T_\tau(\mathbf b)<\infty
	\quad\Longleftrightarrow\quad
	\alpha>\frac12+\frac1\tau,
	\]
	whereas
	\[
	\mathbf b\in\ell^{\tau,2,1/2}
	\quad\Longleftrightarrow\quad
	\alpha>1.
	\]
	In particular, for $\tau=4$ and $\alpha=7/8$,
	$\mathcal T_4(\mathbf b)<\infty$ but
	$\mathbf b\notin\ell^{4,2,1/2}$.
	This distinguishes the two conditions at the sequence level.
\end{remark}

Let $\Psi:[1,\infty)\to(0,\infty)$ be nondecreasing and satisfy
\begin{equation}\label{eq:rates-growth}
	\Psi(v)\le C_0(v/u)^\gamma\Psi(u),
	\qquad 1\le u\le v,
\end{equation}
for some $C_0\ge1$ and $\gamma>0$. The assumption in
Theorem~\ref{thm:rates-rearranged-decay} that $\Psi(t)/t^\gamma$
is nonincreasing gives \eqref{eq:rates-growth} with $C_0=1$.

Let $X$ be a quasi-Banach space and $Y$ a Banach space. For
$T\in\mathcal L(X,Y)$, set
\begin{equation}\label{eq:rates-ideal-definition}
	\|T\|_{\mathfrak A_\Psi}:=\sup_{k\ge1}\Psi(k)a_k(T),
	\qquad
	\mathfrak A_\Psi(X,Y)
	:=\{T\in\mathcal L(X,Y):\|T\|_{\mathfrak A_\Psi}<\infty\}.
\end{equation}
If $\Psi(t)\to\infty$, then $a_k(T)\to0$ for every
$T\in\mathfrak A_\Psi(X,Y)$, so $T$ is compact as a norm limit of
finite-rank operators. For $\Psi(t)=t^{1/s}$ with $s>0$, we write
$\mathfrak A_{s,\infty}(X,Y)$.

\begin{lemma}\label{lem:rates-ideal-properties}
	Let $\Psi:[1,\infty)\to(0,\infty)$ be nondecreasing and satisfy
	\eqref{eq:rates-growth} with constants $C_0$ and $\gamma$.
	Let $X$ be a quasi-Banach space and $Y$ a Banach space.
	Then the following hold.
	\begin{enumerate}
		\item[\textup{(i)}]
		If $T_1,\ldots,T_N\in\mathfrak A_\Psi(X,Y)$, then
		\begin{equation}\label{eq:rates-finite-sum}
			\left\|\sum_{\nu=1}^N T_\nu\right\|_{\mathfrak A_\Psi}
			\le C_0N^\gamma\sum_{\nu=1}^N\|T_\nu\|_{\mathfrak A_\Psi}.
		\end{equation}
		
		\item[\textup{(ii)}]
		If $T_j\in\mathfrak A_\Psi(X,Y)$ and $T_j\to T$ in operator norm, then
		\begin{equation}\label{eq:rates-Fatou}
			\|T\|_{\mathfrak A_\Psi}
			\le\liminf_{j\to\infty}\|T_j\|_{\mathfrak A_\Psi}.
		\end{equation}
		In particular, $T\in\mathfrak A_\Psi(X,Y)$ whenever the right-hand side is finite.
	\end{enumerate}
\end{lemma}

\begin{proof}
	For \textup{(i)}, fix $m\ge1$ and set $k:=\lceil m/N\rceil$.
	Then $N(k-1)+1\le m\le Nk$. Iterating
	Lemma~\ref{lem:ideal-properties}\textup{(ii)} gives
	$a_m(\sum_{\nu=1}^N T_\nu)\le\sum_{\nu=1}^N a_k(T_\nu)$.
	Since $\Psi(m)\le\Psi(Nk)\le C_0N^\gamma\Psi(k)$ by
	\eqref{eq:rates-growth},
	\[
	\Psi(m)a_m\left(\sum_{\nu=1}^N T_\nu\right)
	\le C_0N^\gamma\sum_{\nu=1}^N\Psi(k)a_k(T_\nu)
	\le C_0N^\gamma\sum_{\nu=1}^N\|T_\nu\|_{\mathfrak A_\Psi}.
	\]
	Taking the supremum over $m$ proves \textup{(i)}.
	
	For \textup{(ii)}, Lemma~\ref{lem:ideal-properties}\textup{(iv)}
	gives $a_k(T_j)\to a_k(T)$ for every fixed $k$. Hence
	\[
	\Psi(k)a_k(T)
	=\lim_{j\to\infty}\Psi(k)a_k(T_j)
	\le\liminf_{j\to\infty}\|T_j\|_{\mathfrak A_\Psi}.
	\]
	Taking the supremum over $k$ proves \textup{(ii)}.
	This completes the proof.
\end{proof}

Let $0<s,t<\infty$. For a nonnegative sequence
$\mathbf b=(b_j)_{j\in\mathbb N}$, define
\begin{equation}\label{eq:rates-window}
	W_{s,t}(\mathbf b;k)
	:=\frac1{\sqrt{k}}
	\left(
	\sum_{j=k}^{\lceil k^{s/2}\rceil}
	j^{t/s-1}(b_j^*)^t
	\right)^{1/t},
	\qquad k\in\mathbb N.
\end{equation}
	
\begin{lemma}\label{lem:rates-crossing-window}
	Suppose $0<p<2<q<\infty$ and
	$\mathbf b=(b_j)_{j\in\mathbb N}\in c_0$ is nonnegative.
	There exist $A\in\mathbb N$ and $C>0$, independent of
	$\mathbf b$ and $k$, such that
	\begin{equation}\label{eq:rates-window-upper}
		a_{Ak}(D_{\mathbf b})
		\le C W_{\tau,\varsigma}(\mathbf b;k),
		\qquad k\in \mathbb N,
	\end{equation}
	and
	\begin{equation}\label{eq:rates-window-lower}
		W_{\tau,\varsigma}(\mathbf b;k)
		\le C a_{\max\{1,\lfloor k/A\rfloor\}}(D_{\mathbf b}),
		\qquad k\in \mathbb N.
	\end{equation}
	Consequently,
	\begin{equation}\label{eq:rates-window-ideal}
		\|D_{\mathbf b}\|_{\mathfrak A_\Psi}
		\simeq\sup_{k\ge1}
		\Psi(k)W_{\tau,\varsigma}(\mathbf b;k).
	\end{equation}
\end{lemma}

\begin{proof}
	By \eqref{eq:diagonal-active-coordinates} and invariance under
	coordinate permutations, we may assume $b_j=b_j^*$.
	Put $M:=\lceil k^{\tau/2}\rceil$. We first take $k\ge32$
	with $k^{\tau/2}>64k$.
	
	Suppose $p^\dagger\ne q$, so $\varsigma=\tau$.
	By Lemma~\ref{lem:critical-weak-approximation} with $m=k$,
	there is an operator $\widetilde R$ of rank at most $Ck$ such that
	\[
	\|D_{(b_j)_{k<j\le M}}-\widetilde R\|
	\le Ck^{-1/2}\left(\sum_{j=k+1}^M b_j^\tau\right)^{1/\tau}.
	\]
	Set $R:=D_{(b_j)_{j\le k}}\oplus\widetilde R\oplus0$.
Then
\[
D_{\mathbf b}-R
=0\oplus\bigl(D_{(b_j)_{k<j\le M}}-\widetilde R\bigr)
\oplus D_{(b_j)_{j>M}}.
\]
Since $\operatorname{rank}R\le k+\operatorname{rank}\widetilde R$
and $\|D_{(b_j)_{j>M}}\|=b_{M+1}$, the bound for $\widetilde R$ gives
\begin{equation}\label{eq:rates-window-nonself-approximation}
	\operatorname{rank}R\le Ck,\qquad
	\|D_{\mathbf b}-R\|
	\le Ck^{-1/2}\left(\sum_{j=k+1}^M b_j^\tau\right)^{1/\tau}
	+b_{M+1}.
\end{equation}
	Since $M-k+1\simeq k^{\tau/2}$ and $b_j\ge b_{M+1}$ for $j\le M$,
	\[
	\sqrt{k}\,b_{M+1}
	\le C(M-k+1)^{1/\tau}b_{M+1}
	\le C\left(\sum_{j=k}^M b_j^\tau\right)^{1/\tau}.
	\]
	Hence \eqref{eq:rates-window-nonself-approximation} and
	\eqref{eq:rates-window} give
	$\|D_{\mathbf b}-R\|\le C W_{\tau,\tau}(\mathbf b;k)$.
	
	For the lower estimate, set
	$D^{(M)}:=D_{(b_1,\ldots,b_M)}:\ell_M^p\to\ell_M^q$.
	Apply \eqref{eq:critical-tail-window} to $D^{(M)}$ if $\tau=q$,
	and to $(D^{(M)})^*:\ell_M^{q'}\to\ell_M^{p^\dagger}$
	if $\tau=p^\dagger<q$, in which case $p>1$.
	Finite-dimensional duality and Lemma~\ref{lem:ideal-properties}\textup{(i)}
	give
	\begin{equation}\label{eq:rates-window-nonself-tail}
		\left(\sum_{64k<j\le\lfloor k^{\tau/2}\rfloor}b_j^\tau\right)^{1/\tau}
		\le C\sqrt{k}\,a_k(D^{(M)})
		\le C\sqrt{k}\,a_k(D_{\mathbf b}).
	\end{equation}
	By monotonicity and \eqref{eq:rates-window-nonself-tail},
	\[
	\sum_{j=k}^M b_j^\tau
	\le (63k+2)b_k^\tau
	+\sum_{64k<j\le\lfloor k^{\tau/2}\rfloor}b_j^\tau
	\le Ck b_k^\tau+Ck^{\tau/2}a_k(D_{\mathbf b})^\tau.
	\]
	Using \eqref{eq:rates-window} and then \eqref{eq:prefix-crossing-lower}
	with $m=k$ and approximation index $\lfloor k/2\rfloor$, we get
	\[
	W_{\tau,\tau}(\mathbf b;k)
	\le C\bigl(k^{-\eta_{p,q}}b_k+a_k(D_{\mathbf b})\bigr)
	\le C a_{\lfloor k/2\rfloor}(D_{\mathbf b}).
	\]
	
	Suppose $p^\dagger=q$, so $\tau=q$ and $\varsigma=2$.
	For the upper estimate, put $\varepsilon:=W_{q,2}(\mathbf b;k)$.
	If $\varepsilon=0$, then $b_k=0$ and
	$\operatorname{rank}D_{\mathbf b}<k$.
	Assume $\varepsilon>0$. For $i\ge0$, set
	$n_i:=2^i$ and $I_i:=\{n_i,\ldots,2n_i-1\}$.
	Choose $i_0<i_1$ such that
	\begin{equation}\label{eq:rates-window-dyadic-endpoints}
		2k\le n_{i_0}<4k,\qquad
		\tfrac12 k^{q/2}<n_{i_1}\le k^{q/2}.
	\end{equation}
	For $i_0\le i<i_1$, put $D_i:=D_{\mathbf b|_{I_i}}$.
	By Lemma~\ref{lem:ideal-properties}\textup{(i)} and
	\eqref{eq:finite-crossing-upper},
	\begin{align}\label{rates-window-block-width1}
	a_m(D_i)\le b_{n_i}a_m(I_{n_i})
	\le Cb_{n_i}n_i^{1/q}m^{-1/2},
	\qquad 1\le m\le n_i.
	\end{align}
	Put $r_i:=\min\{n_i,\lceil\varepsilon^{-2}n_i^{2/q}b_{n_i}^2\rceil\}$.
	For $1\le r_i<n_i$, \eqref{rates-window-block-width1} gives
	\[
	a_{r_i}(D_i)\le Cb_{n_i}n_i^{1/q}r_i^{-1/2}
	\le C\varepsilon,
	\]
	since $r_i\ge\varepsilon^{-2}n_i^{2/q}b_{n_i}^2$.
	By the definition of $a_{r_i}(D_i)$, choose $R_i$ with
	$\operatorname{rank}R_i<r_i$ and $\|D_i-R_i\|\le C\varepsilon$.
	For $r_i=0$ or $r_i=n_i$, take $R_i:=D_i$, which is zero when $r_i=0$.
	In all cases,
	\begin{equation}\label{eq:rates-window-block-approximation}
		\operatorname{rank}R_i\le r_i,\qquad
		\|D_i-R_i\|\le C\varepsilon.
	\end{equation}
	By monotonicity,
	\begin{equation}\label{eq:rates-dyadic-comparison}
		n_i^{2/q}b_{n_i}^2
		\le C\sum_{j=n_i/2}^{n_i-1}j^{2/q-1}b_j^2,
		\qquad i\ge1.
	\end{equation}
	For $i_0\le i<i_1$, these summation intervals are disjoint and
	contained in $\{k,\ldots,M\}$. By \eqref{eq:rates-window},
	\begin{equation}\label{eq:rates-window-block-cost}
		\sum_{i=i_0}^{i_1-1}n_i^{2/q}b_{n_i}^2
		\le C\sum_{j=k}^M j^{2/q-1}b_j^2
		=Ck\varepsilon^2.
	\end{equation}
	Set $R:=D_{(b_j)_{j<n_{i_0}}}\oplus
	\bigoplus_{i=i_0}^{i_1-1}R_i\oplus0$.
	Since $i_1-i_0\le C\log k$, the definition of $r_i$ and
	\eqref{eq:rates-window-block-cost} yield
	\begin{equation}\label{eq:rates-window-selfdual-rank-upper}
		\operatorname{rank}R
		\le n_{i_0}-1+\sum_{i=i_0}^{i_1-1}r_i
		\le 4k+C\log k+
		\varepsilon^{-2}\sum_{i=i_0}^{i_1-1}n_i^{2/q}b_{n_i}^2
		\le Ck.
	\end{equation}
	Also, $\{n_{i_1}/2,\ldots,n_{i_1}-1\}\subset\{k,\ldots,M\}$.
	By \eqref{eq:rates-dyadic-comparison} with $i=i_1$ and
	\eqref{eq:rates-window-dyadic-endpoints},
	\[
	k\varepsilon^2
	\ge \sum_{j=n_{i_1}/2}^{n_{i_1}-1}j^{2/q-1}b_j^2
	\ge c n_{i_1}^{2/q}b_{n_{i_1}}^2
	\ge ck b_{n_{i_1}}^2.
	\]
	Since $D_{\mathbf b}-R=0\oplus\bigoplus_{i=i_0}^{i_1-1}(D_i-R_i)
	\oplus D_{(b_j)_{j\ge n_{i_1}}}$, the block estimate
	\eqref{eq:block-error}, Lemma~\ref{lem:diag-facts}\textup{(i)},
	and \eqref{eq:rates-window-block-approximation} give
	\[
	\|D_{\mathbf b}-R\|
	\le\max\left\{\max_{i_0\le i<i_1}\|D_i-R_i\|,\,
	b_{n_{i_1}}\right\}\le C\varepsilon.
	\]
	
	For the lower estimate, put $l:=\lfloor k/16\rfloor$ and
	$c_i:=b_{2n_i-1}$ for $i\ge0$.
	On $I_i$ with $i\ge1$, we have $b_j\le c_{i-1}$, so
	$\sum_{j\in I_i}j^{2/q-1}b_j^2\le Cc_{i-1}^2n_{i-1}^{2/q}$.
	If $I_i\cap\{k,\ldots,M\}\ne\varnothing$, then $k/4<n_{i-1}\le M$.
	Summing over these blocks and reindexing gives
	\begin{equation}\label{eq:rates-window-dyadic-lower}
		\sum_{j=k}^M j^{2/q-1}b_j^2
		\le C\sum_{k/4<n_i\le M}c_i^2n_i^{2/q}.
	\end{equation}
	For $k/4<n_i\le M$, we have $n_i\ge2l$ and
	$b_{n_i}\ge c_i$. By \eqref{eq:prefix-crossing-lower},
	\[
	a_l(D_{\mathbf b})
	\ge cb_{n_i}\min\{1,n_i^{1/q}l^{-1/2}\}
	\ge cc_i n_i^{-\eta_{p,q}}.
	\]
	Choose $\varepsilon>a_l(D_{\mathbf b})$ and fix $C_1\ge1$ such that
	\begin{equation}\label{eq:rates-window-selfdual-threshold}
		c_i n_i^{-\eta_{p,q}}\le C_1\varepsilon,
		\qquad k/4<n_i\le M.
	\end{equation}
	Set $B:=\bigoplus_{k/4<n_i\le M}c_iI_{n_i}$.
	Since $c_i\le b_j$ on $I_i$, coordinate restriction and diagonal
	domination give $a_l(B)\le a_l(D_{\mathbf b})<\varepsilon$.
	Choose $R$ with $\operatorname{rank}R<l$ and
	$\|B-R\|<\varepsilon$. By \eqref{eq:rates-window-selfdual-threshold},
	\[
	k/4<n_i\le M,\quad c_i>2C_1\varepsilon
	\quad\Longrightarrow\quad
	c_i n_i^{-\eta_{p,q}}\le C_1\varepsilon<c_i/2.
	\]
	Thus \eqref{eq:critical-selfdual-rank} applies to these blocks
	with error parameter $C_1\varepsilon$.
	Therefore
	\[
	\sum_{\substack{k/4<n_i\le M\\c_i>2C_1\varepsilon}}
	c_i^2n_i^{2/q}\le Cl\varepsilon^2\le Ck\varepsilon^2.
	\]
	The blocks with $c_i\le2C_1\varepsilon$ contribute at most
	$C\varepsilon^2\sum_{n_i\le M}n_i^{2/q}\le Ck\varepsilon^2$.
	By \eqref{eq:rates-window-dyadic-lower} and \eqref{eq:rates-window},
	$W_{q,2}(\mathbf b;k)^2\le C\varepsilon^2$.
	Letting $\varepsilon\downarrow a_l(D_{\mathbf b})$ gives the lower estimate.
	
	Choose $k_0$ so that the preceding estimates apply for $k\ge k_0$,
	and choose a fixed integer $A\ge\max\{16,k_0\}$ larger than the
	constants in the rank bounds. The constructed approximants have
	rank less than $Ak$, and $\lfloor k/A\rfloor\le\lfloor k/16\rfloor$.
	This gives \eqref{eq:rates-window-upper} and
	\eqref{eq:rates-window-lower} for $k\ge k_0$.
	For $1\le k<k_0$, \eqref{eq:rates-window} gives
	$cb_k\le W_{\tau,\varsigma}(\mathbf b;k)\le Cb_1$.
	Now \eqref{eq:diagonal-truncation} gives $a_{Ak}(D_{\mathbf b})\le b_k$,
	and the lower estimate follows from $a_1(D_{\mathbf b})=b_1$
	and $\max\{1,\lfloor k/A\rfloor\}=1$.
	
	Finally, for $m\ge1$, put $k:=\max\{1,\lfloor m/A\rfloor\}$.
	Since $k\le m<2Ak$, \eqref{eq:rates-growth} and
	\eqref{eq:rates-window-upper} give
	\[
	\Psi(m)a_m(D_{\mathbf b})
	\le C\Psi(k)W_{\tau,\varsigma}(\mathbf b;k).
	\]
	Here for $m<A$ we use $a_m(D_{\mathbf b})\le b_1
	=W_{\tau,\varsigma}(\mathbf b;1)$.
	Conversely, for $l:=\max\{1,\lfloor k/A\rfloor\}$ we have
	$l\le k<2Al$, so \eqref{eq:rates-growth} and
	\eqref{eq:rates-window-lower} give
	\[
	\Psi(k)W_{\tau,\varsigma}(\mathbf b;k)
	\le C\Psi(l)a_l(D_{\mathbf b})
	\le C\|D_{\mathbf b}\|_{\mathfrak A_\Psi}.
	\]
	Taking the respective suprema proves \eqref{eq:rates-window-ideal}.
	This completes the proof.
\end{proof}

\begin{lemma}\label{lem:tail}
	Let $1\le q<p<\infty$ and $1/r=1/q-1/p$.
	For a nonnegative sequence $\mathbf b=(b_j)_{j\in \mathbb N  }\in\ell^r$,
	let $D_{\mathbf b}:\ell^p\to\ell^q$ be the diagonal operator.
	Then, for every $k\in\mathbb N$,
	\begin{equation}\label{eq:tail-formula}
		a_k(D_{\mathbf b})
		=\left(\sum_{j=k}^\infty(b_j^*)^r\right)^{1/r}.
	\end{equation}
\end{lemma}

\begin{proof}
	The case $\mathbf b=0$ is immediate.
	Put $J:=\{j:b_j>0\}$.
	By \eqref{Pi} and Lemma~\ref{lem:ideal-properties}\textup{(i)},
	\[
	D_{\mathbf b}=\iota_JD_{\mathbf b|_J}P_J,\qquad
	D_{\mathbf b|_J}=P_JD_{\mathbf b}\iota_J,\qquad
	a_k(D_{\mathbf b})=a_k(D_{\mathbf b|_J}).
	\]
	Let $J^*:=\{m:b_m^*>0\}$.
	Since $\mathbf b\in\ell^r\subset c_0$, there exists
	a bijection $\pi:J^*\to J$ such that
	$b_{\pi(m)}=b_m^*$.
	For $u\in\{p,q\}$, define the isometric isomorphism
	\[
	U_u:\ell_J^u\to\ell_{J^*}^u,\qquad
	U_u\mathbf x:=(x_{\pi(m)})_{m\in J^*}.
	\]
	Then
	\[
	D_{\mathbf b^*}
	=\iota_{J^*}U_qD_{\mathbf b|_J}U_p^{-1}P_{J^*},
	\qquad
	D_{\mathbf b|_J}
	=U_q^{-1}P_{J^*}D_{\mathbf b^*}\iota_{J^*}U_p.
	\]
	The ideal property therefore gives
	\[
	a_k(D_{\mathbf b})
	=a_k(D_{\mathbf b|_J})
	=a_k(D_{\mathbf b^*}).
	\]
	
	Put $M:=|J|\in\mathbb N\cup\{\infty\}$.
	If $M<k$, then
	\[
	a_k(D_{\mathbf b})=0
	=\left(\sum_{j=k}^\infty(b_j^*)^r\right)^{1/r}.
	\]
	Suppose $k\le M$.
	For $m\ge1$, put $J_m:=\{1,\ldots,m\}$ and
	$\Pi_{J_m}:=\iota_{J_m}P_{J_m}$, with $\Pi_{J_0}:=0$.
	The operator $R_{k-1}:=D_{\mathbf b^*}\Pi_{J_{k-1}}$
	has rank less than $k$.
	By Lemma~\ref{lem:diag-facts}\textup{(ii)},
	\begin{equation}\label{eq:tail-upper}
		a_k(D_{\mathbf b})
		\le\|D_{\mathbf b^*}-R_{k-1}\|
		=\left(\sum_{j=k}^\infty(b_j^*)^r\right)^{1/r}.
	\end{equation}
	
	Fix $N\in\mathbb N$ with $k\le N\le M$, and let
	$B_N:=P_{J_N}D_{\mathbf b^*}\iota_{J_N}:
	\ell_N^p\to\ell_N^q$.
	For $S:\ell_N^p\to\ell_N^q$ with $\rank S<k$, put
	\begin{equation}\label{eq:extreme-space}
		E:=\ker S,\qquad
		d:=\dim E\ge N-k+1,\qquad
		V:=D_{((b_j^*)^{-r/p})_{j=1}^N}E.
	\end{equation}
	Since $b_j^*>0$ for $1\le j\le N$, we have $\dim V=d$.
	The set
	\[
	\mathcal K:=
	\{\mathbf v\in V:\|\mathbf v\|_{\ell_N^\infty}\le1\}
	\]
	is compact.
	Choose $\mathbf w=(w_j)_{j=1}^N\in\mathcal K$
	maximizing $\|\mathbf w\|_{\ell_N^2}^2$, and put
	$K:=\{j\in J_N:|w_j|=1\}$.
	If $|K|<d$, then
	\[
	\dim\{\mathbf v\in V:v_j=0\text{ for }j\in K\}
	\ge d-|K|>0.
	\]
	Choose a nonzero $\mathbf h=(h_j)_{j=1}^N$ in this
	subspace and take
	\[
	0<\varepsilon<
	\min_{\substack{j\notin K\\h_j\ne0}}
	\frac{1-|w_j|}{|h_j|}.
	\]
	Then $\mathbf w\pm\varepsilon\mathbf h\in\mathcal K$,
	whereas maximality implies
	\[
	2\|\mathbf w\|_{\ell_N^2}^2
	\ge\|\mathbf w+\varepsilon\mathbf h\|_{\ell_N^2}^2
	+\|\mathbf w-\varepsilon\mathbf h\|_{\ell_N^2}^2
	=2\|\mathbf w\|_{\ell_N^2}^2
	+2\varepsilon^2\|\mathbf h\|_{\ell_N^2}^2,
	\]
	a contradiction.
	Thus $|K|\ge d\ge N-k+1$.
	
	Define
	\begin{equation}\label{eq:tail-vector}
		\mathbf x:=((b_j^*)^{r/p}w_j)_{j=1}^N\in E,\qquad
		A:=\sum_{j=1}^N(b_j^*)^r|w_j|^q,\qquad
		B:=\sum_{j=1}^N(b_j^*)^r|w_j|^p.
	\end{equation}
	Since $q<p$ and $|w_j|\le1$, we have $0<B\le A$.
	Using $S\mathbf x=0$, $q+rq/p=r$, and
	$1/q-1/p=1/r$, we obtain
	\begin{equation}\label{eq:tail-lower}
		\|B_N-S\|
		\ge\frac{\|(B_N-S)\mathbf x\|_{\ell_N^q}}
		{\|\mathbf x\|_{\ell_N^p}}
		=\frac{A^{1/q}}{B^{1/p}}
		\ge A^{1/r}
		\ge\left(\sum_{j\in K}(b_j^*)^r\right)^{1/r}
		\ge\left(\sum_{j=k}^N(b_j^*)^r\right)^{1/r},
	\end{equation}
	where the last inequality follows from
	$b_1^*\ge\cdots\ge b_N^*$ and $|K|\ge N-k+1$.
	Taking the infimum over $S$ and applying
	Lemma~\ref{lem:ideal-properties}\textup{(i)} gives
	\[
	a_k(D_{\mathbf b})
	\ge a_k(B_N)
	\ge\left(\sum_{j=k}^N(b_j^*)^r\right)^{1/r}.
	\]
	Consequently,
	\[
	a_k(D_{\mathbf b})
	\ge\sup_{\substack{N\in\mathbb N\\k\le N\le M}}
	\left(\sum_{j=k}^N(b_j^*)^r\right)^{1/r}
	=\left(\sum_{j=k}^\infty(b_j^*)^r\right)^{1/r}.
	\]
	Together with \eqref{eq:tail-upper}, this proves
	\eqref{eq:tail-formula}.
	This completes the proof.
\end{proof}

\section{Sampling and local multiplication estimates}\label{sec:sampling}

In this section, we establish sampling comparisons for the IDA function
and approximation estimates for local multiplication operators.
These estimates will be used to compare Hankel operators with
diagonal operators in the next section.
\subsection{Decomposition and radius comparisons}
For $1\le q<\infty$, $r>0$, and
$u\in L^q_{\rm loc}(\C^n)$, define
\begin{equation}\label{eq:local-mean}
	M_{q,r}(u)(z)
	:=\left(\frac1{v(B(z,r))}
	\int_{B(z,r)}|u(w)|^q\,dv(w)\right)^{1/q}.
\end{equation}
For a differential form $\omega$, set
$M_{q,r}(\omega):=M_{q,r}(|\omega|)$.

The following lemma follows from
\cite[Lemma 3.6 and Theorem 3.8]{HV23} and
\cite[Corollary 3.10]{HV23}.

\begin{lemma}\label{lem:known-decomposition}
	Let $1\le q<\infty$, $t>0$, and
	$f\in L^q_{\rm loc}(\C^n)$.
	Then there exist $f_1\in C^2(\C^n)$ and
	$f_2\in L^q_{\rm loc}(\C^n)$ such that $f=f_1+f_2$ and
	\begin{equation}\label{eq:decomposition}
		|\bar\partial f_1(z)|+M_{q,t/2}(\bar\partial f_1)(z)
		+M_{q,t/2}(f_2)(z)
		\le C_tG_{q,2t}(f)(z),\qquad z\in\C^n.
	\end{equation}
	Moreover, for every $0<s<\infty$ and $r,R>0$,
	\begin{equation}\label{eq:radius}
		\|G_{q,r}(f)\|_{L^s}\simeq\|G_{q,R}(f)\|_{L^s}.
	\end{equation}
\end{lemma}

The following lemma shows that the Lorentz and
Lorentz--Zygmund conditions on $G_{q,r}(f)$
are independent of the choice of $r>0$.

\begin{lemma}\label{lem:rearranged-radius}
	\label{lem:logarithmic-rearrangements}
	Let $1\le q<\infty$, $0<\sigma,s<\infty$,
	$\beta\in\mathbb R$, and $f\in L^q_{\rm loc}(\C^n)$.
	Then, for every $r_1,r_2>0$,
	\begin{equation}\label{eq:logarithmic-radius}
		\|G_{q,r_1}(f)\|_{\mathcal L^{\sigma,s,\beta}}
		\simeq\|G_{q,r_2}(f)\|_{\mathcal L^{\sigma,s,\beta}}.
	\end{equation}
	In particular, taking $\beta=0$ gives
	\begin{equation}\label{eq:rearranged-radius}
		\|G_{q,r_1}(f)\|_{\mathcal L^{\sigma,s}}
		\simeq\|G_{q,r_2}(f)\|_{\mathcal L^{\sigma,s}}.
	\end{equation}
	The implicit constants are independent of $f$.
\end{lemma}

\begin{proof}
	Fix $r_1,r_2>0$. By Lemma~\ref{lem:known-decomposition} with
	parameter $r_1/2$, write $f=f_1+f_2$, where
	$f_1\in C^2(\C^n)$ and $f_2\in L^q_{\rm loc}(\C^n)$ satisfy
	\begin{equation}\label{eq:radius-local-decomposition}
		M_{q,r_1/4}(\bar\partial f_1)(z)+M_{q,r_1/4}(f_2)(z)
		\le C G_{q,r_1}(f)(z).
	\end{equation}
	The Henkin estimate \cite[(3-16)]{HV23}, applied on $B(z,2r_2)$
	as in the proof of \cite[Theorem~3.8]{HV23}, gives
	\begin{equation}\label{eq:radius-local-distance}
		G_{q,r_2}(f)(z)
		\le C\bigl(M_{q,2r_2}(\bar\partial f_1)(z)
		+M_{q,2r_2}(f_2)(z)\bigr).
	\end{equation}
	Choose $v_1,\ldots,v_N\in\C^n$ such that
	\begin{equation}\label{eq:radius-cover}
		B(0,2r_2)\subseteq\bigcup_{\nu=1}^N B(v_\nu,r_1/4).
	\end{equation}
	Translating the covering by $z$ gives, for
	$u=\bar\partial f_1$ or $u=f_2$,
	\begin{equation}\label{eq:radius-average-cover}
		M_{q,2r_2}(u)(z)^q
		\le \frac{v(B(0,r_1/4))}{v(B(0,2r_2))}
		\sum_{\nu=1}^N M_{q,r_1/4}(u)(z+v_\nu)^q.
	\end{equation}
	Taking $q$th roots in \eqref{eq:radius-average-cover} and using $q\ge1$ gives
	\[
	M_{q,2r_2}(u)(z)
	\le C\left(\sum_{\nu=1}^N M_{q,r_1/4}(u)(z+v_\nu)^q\right)^{1/q}
	\le C\sum_{\nu=1}^N M_{q,r_1/4}(u)(z+v_\nu),
	\qquad u=\bar\partial f_1,\ f_2.
	\]
	Apply this estimate to both terms on the right of
	\eqref{eq:radius-local-distance}, and then use
	\eqref{eq:radius-local-decomposition} at each $z+v_\nu$ to obtain
	\begin{equation}\label{eq:radius-translation-comparison}
		\begin{aligned}
			G_{q,r_2}(f)(z)
			\le C\sum_{\nu=1}^N
			\bigl[M_{q,r_1/4}(\bar\partial f_1)(z+v_\nu)
			+M_{q,r_1/4}(f_2)(z+v_\nu)\bigr]
			\le C\sum_{\nu=1}^N G_{q,r_1}(f)(z+v_\nu).
		\end{aligned}
	\end{equation}
	
	For nonnegative measurable functions $F_1,\ldots,F_N$ and
	$a_1,\ldots,a_N\ge0$,
	\begin{equation}\label{eq:rearrangement-sum-distribution}
		\left\{\sum_{\nu=1}^N F_\nu>\sum_{\nu=1}^N a_\nu\right\}
		\subseteq\bigcup_{\nu=1}^N\{F_\nu>a_\nu\}.
	\end{equation}
	If all $F_\nu^*(t)$ are finite, take $a_\nu=F_\nu^*(t)+\varepsilon$.
	By \eqref{eq:rearrangement-sum-distribution},
	\[
	d_{\sum_{\nu=1}^N F_\nu}
	\left(\sum_{\nu=1}^N a_\nu\right)
	\le \sum_{\nu=1}^N d_{F_\nu}(a_\nu)
	\le Nt.
	\] By the definition of rearrangement
	and $\varepsilon\downarrow0$,
	\begin{equation}\label{eq:rearrangement-sum}
		\left(\sum_{\nu=1}^N F_\nu\right)^*(Nt)
		\le\sum_{\nu=1}^N F_\nu^*(t),\qquad t>0.
	\end{equation}
	The inequality also holds if the right-hand side is infinite.
	
	Recall that $\lambda(t):=1+\max\{0,\log t\}$. Since
	$\lambda(ct)\le\lambda(t)+|\log c|$ and $\lambda(t)\ge1$,
	applying the same bound to the dilation by $c^{-1}$ gives
	\begin{equation}\label{eq:logarithmic-dilation}
		(1+|\log c|)^{-1}\lambda(t)
		\le\lambda(ct)\le(1+|\log c|)\lambda(t),
		\qquad c,t>0.
	\end{equation}
	In particular, for $\beta\in\mathbb R$,
	\begin{equation}\label{eq:logarithmic-weight-dilation}
		\lambda(ct)^{\beta s}
		\le (1+|\log c|)^{|\beta|s}\lambda(t)^{\beta s}.
	\end{equation}
	Rescaling the integral in the definition of the norm gives
	\begin{equation}\label{eq:logarithmic-sum-proof}
		\left\|\sum_{\nu=1}^N F_\nu\right\|_{\mathcal L^{\sigma,s,\beta}}^s
		=N^{s/\sigma}\int_0^\infty t^{s/\sigma-1}
		\lambda(Nt)^{\beta s}
		\left[\left(\sum_{\nu=1}^N F_\nu\right)^*(Nt)\right]^s\,dt.
	\end{equation}
	By \eqref{eq:rearrangement-sum},
	\eqref{eq:logarithmic-weight-dilation}, and
	$(\sum_{\nu=1}^N a_\nu)^s\le N^{(s-1)_+}\sum_{\nu=1}^N a_\nu^s$,
	\begin{equation}\label{eq:logarithmic-sum}
		\left\|\sum_{\nu=1}^N F_\nu\right\|_{\mathcal L^{\sigma,s,\beta}}^s
		\le N^{s/\sigma+(s-1)_+}(1+\log N)^{|\beta|s}
		\sum_{\nu=1}^N\|F_\nu\|_{\mathcal L^{\sigma,s,\beta}}^s.
	\end{equation}
	For $\beta=0$, this gives
	\begin{equation}\label{eq:rearrangement-sum-norm}
		\left\|\sum_{\nu=1}^N F_\nu\right\|_{\mathcal L^{\sigma,s}}^s
		\le N^{s/\sigma+(s-1)_+}
		\sum_{\nu=1}^N\|F_\nu\|_{\mathcal L^{\sigma,s}}^s.
	\end{equation}
	Translations preserve distribution functions, so
	\begin{equation}\label{eq:logarithmic-translation}
		\|F(\cdot+a)\|_{\mathcal L^{\sigma,s,\beta}}
		=\|F\|_{\mathcal L^{\sigma,s,\beta}},\qquad a\in\C^n.
	\end{equation}
	
	Applying \eqref{eq:logarithmic-sum} to the translates in
	\eqref{eq:radius-translation-comparison} and using
	\eqref{eq:logarithmic-translation}, we obtain
	\begin{equation}\label{eq:logarithmic-radius-upper}
		\|G_{q,r_2}(f)\|_{\mathcal L^{\sigma,s,\beta}}^s
		\le C\|G_{q,r_1}(f)\|_{\mathcal L^{\sigma,s,\beta}}^s.
	\end{equation}
	Interchanging $r_1$ and $r_2$ proves \eqref{eq:logarithmic-radius}.
	Taking $\beta=0$ gives \eqref{eq:rearranged-radius}.
	This completes the proof.
\end{proof}

\subsection{Sampling on cubic grids}

Let $r_0$ be given by Lemma~\ref{lem:known-kernel}, and define
\[
\mathscr G:=
\left\{
\xi+h\mathbb Z^{2n}:
\xi\in\mathbb R^{2n},\
0<h<\frac{r_0}{10\sqrt{2n}}
\right\}.
\]
For $\Lambda=\xi+h\mathbb Z^{2n}\in\mathscr G$, choose any
enumeration $\Lambda=\{a_j\}_{j\in\mathbb N}$ and put
\[
Q_j:=a_j+[-h/2,h/2)^{2n},\qquad
\delta:=\frac{\sqrt{2n}}2h.
\]
Then $\{Q_j\}_{j\in\mathbb N}$ is a partition of
$\mathbb R^{2n}$ into half-open cubes with centers $a_j$.
Set
\begin{equation}\label{eq:parameters}
	t:=8\delta,\qquad R:=20\delta<r_0.
\end{equation}
For every $j\in\mathbb N$ and $z\in Q_j$,
\begin{equation}\label{eq:sampling-ball-inclusions}
	B(z,R/2)\subset B(a_j,R)\subset B(z,2R).
\end{equation}

\begin{lemma}\label{lem:power-rearrangement-sampling}
	\label{lem:logarithmic-sampling}
	Let $1\le q<\infty$, $0<\sigma,s<\infty$,
	$\beta\in\mathbb R$, and $f\in L^q_{\rm loc}(\C^n)$.
	For $\Lambda=\{a_j\}_{j\in\mathbb N}\in\mathscr G$,
	let $R$ be as in \eqref{eq:parameters} and put
	$\mathbf g:=(G_{q,R}(f)(a_j))_{j\in\mathbb N}$.
	Then
	\begin{equation}\label{eq:logarithmic-sampling}
		\|G_{q,1}(f)\|_{\mathcal L^{\sigma,s,\beta}}
		\simeq\|\mathbf g\|_{\ell^{\sigma,s,\beta}}.
	\end{equation}
	In particular, taking $\beta=0$ gives
	\begin{equation}\label{eq:power-sampling}
		\|G_{q,1}(f)\|_{\mathcal L^{\sigma,s}}
		\simeq\|\mathbf g\|_{\ell^{\sigma,s}}.
	\end{equation}
	For $\sigma>2$, one also has
	\begin{equation}\label{eq:critical-Copson-sampling}
		\mathcal T_\sigma(\mathbf g)^2
		\simeq
		\|G_{q,1}(f)\|_{\mathcal L^{\sigma,\sigma}}^2
		+\|G_{q,1}(f)\|_{\mathcal{CL}^{\sigma,2}}^2.
	\end{equation}
	The comparison constants are independent of $f$.
\end{lemma}

\begin{proof}
	Set $g_j:=G_{q,R}(f)(a_j)$, $v_h:=h^{2n}$, and
	$F_h:=\sum_{j\ge1}g_j\chi_{Q_j}$. By
	\eqref{eq:sampling-ball-inclusions},
	\begin{equation}\label{eq:power-sampling-sandwich}
		2^{-2n/q}G_{q,R/2}(f)\le F_h\le2^{2n/q}G_{q,2R}(f).
	\end{equation}
	By monotonicity and \eqref{eq:logarithmic-radius},
	\begin{equation}\label{eq:logarithmic-sampling-radius}
		\|G_{q,1}(f)\|_{\mathcal L^{\sigma,s,\beta}}
		\simeq\|F_h\|_{\mathcal L^{\sigma,s,\beta}}.
	\end{equation}
	Since the cubes are disjoint and $v(Q_j)=v_h$,
	$d_{F_h}(u)=v_h\#\{j\ge1:g_j>u\}$. Hence
	\begin{equation}\label{eq:power-sampling-step}
		F_h^*(t)=g_k^*,
		\qquad (k-1)v_h<t<kv_h,\quad k\ge1.
	\end{equation}
	Substituting \eqref{eq:power-sampling-step} into the norm gives
	\begin{equation}\label{eq:logarithmic-step-function}
		\|F_h\|_{\mathcal L^{\sigma,s,\beta}}^s
		=\sum_{k\ge1}(g_k^*)^s
		\int_{(k-1)v_h}^{kv_h}
		t^{s/\sigma-1}\lambda(t)^{\beta s}\,dt.
	\end{equation}
	By the substitution $t=v_hu$ and \eqref{eq:logarithmic-dilation},
	\begin{equation}\label{eq:logarithmic-cell-weight}
		\int_{(k-1)v_h}^{kv_h}t^{s/\sigma-1}\lambda(t)^{\beta s}\,dt
		\simeq v_h^{s/\sigma}
		k^{s/\sigma-1}(1+\log k)^{\beta s}.
	\end{equation}
	Indeed, for $u\in(k-1,k)$ with $k\ge2$, we have $u\simeq k$ and
	$\lambda(u)\simeq1+\log k$. For $k=1$, the rescaled integral
	with weight $\lambda(u)^{\beta s}$ equals $\sigma/s$.
	Thus \eqref{eq:logarithmic-step-function} gives
	\begin{equation}\label{eq:logarithmic-sampling-sum}
		\|F_h\|_{\mathcal L^{\sigma,s,\beta}}^s
		\simeq v_h^{s/\sigma}\|\mathbf g\|_{\ell^{\sigma,s,\beta}}^s.
	\end{equation}
	Since $h$ is fixed, \eqref{eq:logarithmic-sampling-radius} proves
	\eqref{eq:logarithmic-sampling}. Taking $\beta=0$ gives
	\eqref{eq:power-sampling}.
	
	For \eqref{eq:critical-Copson-sampling}, let $\sigma>2$ and write
	$A_F(t):=\int_t^\infty F^*(u)^\sigma\,du$.
	If $F:=\sum_{\nu=1}^N F_\nu$ with $F_\nu\ge0$, then
	\eqref{eq:rearrangement-sum} and H\"older's inequality give
	\[
	A_F(Nt)
	=N\int_t^\infty F^*(Nu)^\sigma\,du
	\le N^\sigma\sum_{\nu=1}^N A_{F_\nu}(t).
	\]
	Since $2/\sigma<1$ and $N/(1+Nt)\le N/(1+t)$,
	\[
	\|F\|_{\mathcal{CL}^{\sigma,2}}^2
	=\int_0^\infty A_F(Nt)^{2/\sigma}\frac{N\,dt}{1+Nt}
	\le N^3\sum_{\nu=1}^N\|F_\nu\|_{\mathcal{CL}^{\sigma,2}}^2.
	\]
	Together with $\|F\|_{L^\sigma}^2\le N\sum_{\nu=1}^N
	\|F_\nu\|_{L^\sigma}^2$, this gives
	\begin{equation}\label{eq:critical-Copson-sum}
		\|F\|_{L^\sigma}^2+\|F\|_{\mathcal{CL}^{\sigma,2}}^2
		\le N^3\sum_{\nu=1}^N
		\left(\|F_\nu\|_{L^\sigma}^2+
		\|F_\nu\|_{\mathcal{CL}^{\sigma,2}}^2\right).
	\end{equation}
In \eqref{eq:critical-Copson-sum}, take
$F_\nu:=G_{q,r_1}(f)(\cdot+v_\nu)$ for $1\le\nu\le N$ and
$F:=\sum_{\nu=1}^N F_\nu$. Then
\eqref{eq:radius-translation-comparison} gives $G_{q,r_2}(f)\le CF$.
By translation invariance,
\[
\|F_\nu\|_{L^\sigma}=\|G_{q,r_1}(f)\|_{L^\sigma},
\qquad
\|F_\nu\|_{\mathcal{CL}^{\sigma,2}}
=\|G_{q,r_1}(f)\|_{\mathcal{CL}^{\sigma,2}}.
\]
By monotonicity, \eqref{eq:critical-Copson-sum}, and translation
invariance,
\[
\|G_{q,r_2}(f)\|_{L^\sigma}
+\|G_{q,r_2}(f)\|_{\mathcal{CL}^{\sigma,2}}
\le C\left(
\|G_{q,r_1}(f)\|_{L^\sigma}
+\|G_{q,r_1}(f)\|_{\mathcal{CL}^{\sigma,2}}
\right).
\] Interchanging $r_1$ and $r_2$ gives
\begin{equation}\label{eq:critical-Copson-radius}
	\|G_{q,r_1}(f)\|_{L^\sigma}
	+\|G_{q,r_1}(f)\|_{\mathcal{CL}^{\sigma,2}}
	\simeq
	\|G_{q,r_2}(f)\|_{L^\sigma}
	+\|G_{q,r_2}(f)\|_{\mathcal{CL}^{\sigma,2}}.
\end{equation}
	
	Set $S_k:=\sum_{j\ge k}(g_j^*)^\sigma$. By
	\eqref{eq:power-sampling-step}, $\|F_h\|_{L^\sigma}^2
	=v_h^{2/\sigma}S_1^{2/\sigma}$ and
	\[
	v_hS_{k+1}\le A_{F_h}(v_ht)\le v_hS_k,
	\qquad k-1<t<k,\quad k\ge1.
	\]
	Changing variables in the Copson integral and using
	$v_h/(1+v_ht)\simeq(1+t)^{-1}\simeq k^{-1}$ for  $t\in(k-1,k)$,
	we obtain
	\[
	v_h^{2/\sigma}\sum_{k\ge1}\frac{S_{k+1}^{2/\sigma}}k
	\lesssim\|F_h\|_{\mathcal{CL}^{\sigma,2}}^2
	\lesssim v_h^{2/\sigma}\sum_{k\ge1}\frac{S_k^{2/\sigma}}k.
	\]
	Reindexing the lower sum and adding the $L^\sigma$ term gives
	\[
	\|F_h\|_{L^\sigma}^2+\|F_h\|_{\mathcal{CL}^{\sigma,2}}^2
	\simeq v_h^{2/\sigma}\mathcal T_\sigma(\mathbf g)^2.
	\]
	Now \eqref{eq:power-sampling-sandwich} and
	\eqref{eq:critical-Copson-radius} yield
	\eqref{eq:critical-Copson-sampling}.
	This completes the proof.
\end{proof}

\begin{lemma}\label{lem:indexed-sampling}
	Let $1\le q<\infty$, $f\in L^q_{\rm loc}(\C^n)$, and
	$\Lambda=\{a_j\}_{j\in\mathbb N}\in\mathscr G$, with $R$ as in
	\eqref{eq:parameters}.
	Let $\mathbf g=(g_j)_{j\in\mathbb N}$, where
	$g_j:=G_{q,R}(f)(a_j)$.
	There exists $C\ge1$, independent of $f$, $k$, and $t$,
	such that, for every $k\in\mathbb N$ and every $t>0$
	satisfying $k-1\le t<k$,
	\begin{equation}\label{eq:rates-indexed-sampling}
		C^{-1}\bigl(G_{q,1}(f)\bigr)^*(Ct)
		\le g_k^*
		\le C\bigl(G_{q,1}(f)\bigr)^*(t/C).
	\end{equation}
\end{lemma}

\begin{proof}
	Fix $k\in\mathbb N$ and $t>0$ with $k-1\le t<k$.
	Set $F_h:=\sum_{j\in\mathbb N}g_j\chi_{Q_j}$.
	By \eqref{eq:power-sampling-step} and the definition of rearrangement,
	we have $
	F_h^*(h^{2n}t)=g_k^*.
	$
 Applying
	\eqref{eq:radius-translation-comparison} and
	\eqref{eq:rearrangement-sum} with
	$(r_1,r_2)=(R/2,1)$ and $(1,2R)$ gives a common fixed
	integer $N\ge1$ such that
	\[
	\bigl(G_{q,1}(f)\bigr)^*(Nu)
	\lesssim \bigl(G_{q,R/2}(f)\bigr)^*(u),
	\]
	and
	\[
	\bigl(G_{q,2R}(f)\bigr)^*(Nu)
	\lesssim \bigl(G_{q,1}(f)\bigr)^*(u),
	\qquad u>0.
	\]
	Together with \eqref{eq:power-sampling-sandwich}, these imply
	\[
	C_0^{-1}\bigl(G_{q,1}(f)\bigr)^*(C_0u)
	\le F_h^*(u)
	\le C_0\bigl(G_{q,1}(f)\bigr)^*(u/C_0),
	\qquad u>0,
	\]
	where $C_0\ge1$ is independent of $f$ and $u$.
	Taking $u=h^{2n}t$ yields
	\[
	C_0^{-1}\bigl(G_{q,1}(f)\bigr)^*(C_0h^{2n}t)
	\le g_k^*
	\le C_0\bigl(G_{q,1}(f)\bigr)^*(h^{2n}t/C_0).
	\]
Choose $C:=C_0\max\{1,h^{2n},h^{-2n}\}$. Monotonicity then gives \eqref{eq:rates-indexed-sampling}. This completes the proof.
\end{proof}

\begin{lemma}\label{lem:continuous-tail}
	Let $1\le q<p<\infty$, $0<s<\infty$, and
	$f\in L^q_{\rm loc}(\C^n)$.
	For $\Lambda=\{a_j\}_{j\in\mathbb N}\in\mathscr G$,
	let $R$ be as in \eqref{eq:parameters} and put
	$\mathbf g:=(G_{q,R}(f)(a_j))_{j\in\mathbb N}$.
	Then $D_{\mathbf g}\in\mathfrak A_s(\ell^p,\ell^q)$ if and only if
	$G_{q,1}(f)\in\mathcal L^{\rho,s}(\C^n,dv)$. Moreover, for $1/r:=1/q-1/p$, we have
	\begin{equation}\label{eq:continuous-tail-bridge}
		\|G_{q,1}(f)\|_{\mathcal L^{\rho,s}}^s
		\simeq\sum_{k\ge1}
		\left(\sum_{j=k}^\infty(g_j^*)^r\right)^{s/r}
		=\As{D_{\mathbf g}}^s.
	\end{equation}
	The comparison constants are independent of $f$.
\end{lemma}

\begin{proof}
	Set $u:=s/r$ and, for $j\ge1$, put
	$g_j:=G_{q,R}(f)(a_j)$ and $y_j:=(g_j^*)^r$.
	Fix $N\ge1$ and work with the first $N$ terms.
	For $2\le k\le N$, monotonicity gives
	\[
	\sum_{j=\lfloor k/2\rfloor}^{N}y_j
	\ge \sum_{j=\lfloor k/2\rfloor}^{k}y_j
	\ge \frac{k}{2}y_k.
	\]
	Each value of $\lfloor k/2\rfloor$ occurs at most twice for
	$2\le k\le N$. Raising the preceding inequality to the power $u$,
	summing over $k$, and including $y_1\le\sum_{j=1}^N y_j$ gives
	\begin{equation}\label{eq:hardy-tail-lower}
		\sum_{k=1}^N k^uy_k^u
		\le (1+2^{u+1})\sum_{k=1}^N
		\left(\sum_{j=k}^N y_j\right)^u.
	\end{equation}
	
	For $1\le m\le N$ and $1\le k\le\lfloor N/m\rfloor$, set
	$A_{m,k}:=\{m(k-1)+1,\ldots,mk\}$.
	If $k\ge2$ and $j\in A_{m,k}$, then $j\ge mk/2$ and $y_j\ge y_{mk}$.
	Since $A_{m,k}$ contains $m$ indices,
	\begin{equation}\label{eq:hardy-block}
		\sum_{j\in A_{m,k}}j^uy_j^u
		\ge 2^{-u}m^{u+1}k^uy_{mk}^u.
	\end{equation}
	For the first block, $y_j\ge y_m$ and
	$\int_{j-1}^j t^u\,dt\le j^u$ for $1\le j\le m$. Hence
	\begin{equation}\label{eq:hardy-first-block}
		\sum_{j\in A_{m,1}}j^uy_j^u
		\ge y_m^u\int_0^m t^u\,dt
		=\frac{m^{u+1}}{u+1}y_m^u.
	\end{equation}
	Summing \eqref{eq:hardy-block} and \eqref{eq:hardy-first-block}
	over the disjoint blocks contained in $\{1,\ldots,N\}$ gives
	\begin{equation}\label{eq:Hardy-dilation}
		\sum_{k=1}^{\lfloor N/m\rfloor}k^uy_{mk}^u
		\le Cm^{-(u+1)}\sum_{j=1}^N j^uy_j^u,
	\end{equation}
	where $C$ depends only on $u$.
	
	For $1\le k\le N$, split the tail into dyadic blocks, ending the
	last block at $N$. By monotonicity,
	\begin{equation}\label{eq:hardy-dyadic-tail}
		\sum_{j=k}^N y_j
		=\sum_{\nu=0}^{\lfloor\log_2(N/k)\rfloor}
		\sum_{j=2^\nu k}^{\min\{2^{\nu+1}k-1,N\}}y_j
		\le\sum_{\nu=0}^{\lfloor\log_2(N/k)\rfloor}
		2^\nu k\,y_{2^\nu k}.
	\end{equation}
	If $0<u\le1$, subadditivity in \eqref{eq:hardy-dyadic-tail}
	and \eqref{eq:Hardy-dilation} with $m=2^\nu$ give
	\begin{equation}\label{eq:hardy-small-exponent}
		\sum_{k=1}^N\left(\sum_{j=k}^N y_j\right)^u
		\le\sum_{\nu=0}^{\lfloor\log_2 N\rfloor}2^{\nu u}
		\sum_{k=1}^{\lfloor N/2^\nu\rfloor}k^uy_{2^\nu k}^u
		\le C\sum_{j=1}^N j^uy_j^u,
	\end{equation}
	where the last bound uses $\sum_{\nu\ge0}2^{-\nu}<\infty$.
	If $u>1$, Minkowski's inequality in \eqref{eq:hardy-dyadic-tail}
	and \eqref{eq:Hardy-dilation} with $m=2^\nu$ give
	\begin{equation}\label{eq:hardy-large-exponent}
		\left[\sum_{k=1}^N\left(\sum_{j=k}^N y_j\right)^u\right]^{1/u}
		\le\sum_{\nu=0}^{\lfloor\log_2 N\rfloor}2^\nu
		\left(\sum_{k=1}^{\lfloor N/2^\nu\rfloor}
		k^uy_{2^\nu k}^u\right)^{1/u}
		\le C\left(\sum_{j=1}^N j^uy_j^u\right)^{1/u},
	\end{equation}
	using $\sum_{\nu\ge0}2^{-\nu/u}<\infty$.
	
	The constants in these bounds are independent of $N$.
	As $N\to\infty$, monotone convergence gives
	\[
	\sum_{k=1}^N k^uy_k^u\uparrow\sum_{k\ge1}k^uy_k^u,
	\qquad
	\sum_{k=1}^N\left(\sum_{j=k}^N y_j\right)^u
	\uparrow\sum_{k\ge1}\left(\sum_{j\ge k}y_j\right)^u.
	\]
	Passing to the limit in \eqref{eq:hardy-tail-lower} and the
	corresponding upper bound, and substituting $u=s/r$ and
	$y_j=(g_j^*)^r$, we obtain
	\begin{equation}\label{eq:hardy-tail}
		\sum_{k\ge1}\left(\sum_{j\ge k}(g_j^*)^r\right)^{s/r}
		\simeq\sum_{k\ge1}k^{s/r}(g_k^*)^s.
	\end{equation}
	
	Since $s/\rho-1=s/r$, Lemma~\ref{lem:power-rearrangement-sampling}
	and \eqref{eq:hardy-tail} give
	\begin{equation}\label{eq:continuous-tail-comparison}
		\|G_{q,1}(f)\|_{\mathcal L^{\rho,s}}^s
		\simeq\sum_{k\ge1}k^{s/r}(g_k^*)^s
		\simeq\sum_{k\ge1}\left(\sum_{j\ge k}(g_j^*)^r\right)^{s/r}.
	\end{equation}
	If $G_{q,1}(f)\in\mathcal L^{\rho,s}$, the $k=1$ term in the
	last sum gives
	\begin{equation}\label{eq:continuous-tail-boundedness}
		\|\mathbf g\|_{\ell^r}^s
		\lesssim\|G_{q,1}(f)\|_{\mathcal L^{\rho,s}}^s.
	\end{equation}
	If $D_{\mathbf g}\in\mathfrak A_s(\ell^p,\ell^q)$, then
	$\mathbf g\in\ell^r$ by Lemma~\ref{lem:diag-facts}\textup{(ii)}.
	In either case, Lemma~\ref{lem:tail} and
	\eqref{eq:continuous-tail-comparison} yield
	\begin{equation}\label{eq:continuous-tail-ideal}
		\As{D_{\mathbf g}}^s
		=\sum_{k\ge1}\left(\sum_{j\ge k}(g_j^*)^r\right)^{s/r}
		\simeq\|G_{q,1}(f)\|_{\mathcal L^{\rho,s}}^s.
	\end{equation}
	This completes the proof.
\end{proof}

\subsection{Local approximation of multiplication operators}
For $1\le q<\infty$ and $d\in\mathbb N$, let
$L^q_{\rm loc}(\C^n;\C^d)$ be the space of measurable maps
$\vartheta:\C^n\to\C^d$ given by
\[
z\mapsto\vartheta(z)
:=(\vartheta_1(z),\ldots,\vartheta_d(z)),
\]
such that $|\vartheta|\in L^q_{\rm loc}(\C^n)$, where
\[
|\vartheta(z)|
:=\left(\sum_{\ell=1}^d|\vartheta_\ell(z)|^2\right)^{1/2}.
\]

	\begin{lemma}\label{lem:transfer}\label{lem:weighted-transfer}
		Suppose that $\operatorname{Hess}_{\mathbb R}\varphi\simeq E$.
		Let $0<p<\infty$, $1\le q<\infty$, $0<s<\infty$,
		$d\in\mathbb N$, and
		$\vartheta\in L^q_{\rm loc}(\C^n;\C^d)$.
		For $\Lambda=\{a_j\}_{j\in \mathbb N}\in\mathscr G$, let $Q_j$
		be the corresponding cube centered at $a_j$, and put
		$\mathbf b:=(b_j)_{j\in \mathbb N}$, where
		\[
		b_j:=\left(\int_{Q_j}|\vartheta(z)|^q\,dv(z)\right)^{1/q}.
		\]
		If $D_{\mathbf b}\in\mathfrak A_s(\ell^p,\ell^q)$,
		then $J_\vartheta u:=\vartheta u$ belongs to
		$\mathfrak A_s(F^p(\varphi),L^q(\varphi;\C^d))$ and
		\begin{equation}\label{eq:transfer}
			\As{J_\vartheta}\lesssim \As{D_{\mathbf b}}.
		\end{equation}
		If $\Psi$ satisfies \eqref{eq:rates-growth} and
		$D_{\mathbf b}\in\mathfrak A_\Psi(\ell^p,\ell^q)$, then
		$J_\vartheta\in
		\mathfrak A_\Psi(F^p(\varphi),L^q(\varphi;\C^d))$ and
		\begin{equation}\label{eq:rates-multiplication}
			\|J_\vartheta\|_{\mathfrak A_\Psi}
			\lesssim\|D_{\mathbf b}\|_{\mathfrak A_\Psi}.
		\end{equation}
	The implicit constants are independent of $\vartheta$.
	\end{lemma}
	
	\begin{proof}
		Let $\delta$ be the cube radius defined before \eqref{eq:parameters}
		and choose $R_B>4\delta$. For $j\in\mathbb N$, put
		$B_j:=B(a_j,R_B)$ and define
		\begin{equation}\label{eq:local-spaces}
			X_j:=\mathcal H(B_j)\cap L^p(B_j,e^{-p\varphi}dv),
			\qquad
			\|u\|_{X_j}:=\left(\int_{B_j}|u|^pe^{-p\varphi}\,dv\right)^{1/p},
		\end{equation}
and
	\[
	Y_j:=L^q(Q_j,e^{-q\varphi}dv;\C^d),
	\qquad
	\|w\|_{Y_j}:=
	\left(\int_{Q_j}|w|^{q}e^{-q\varphi}\,dv\right)^{1/q}.
	\]
		Set $X:=(\bigoplus_{j\in\mathbb N}X_j)_{\ell^p}$ and
		$Y:=(\bigoplus_{j\in\mathbb N}Y_j)_{\ell^q}$, with
		\[
		\|\mathbf u\|_X^p:=\sum_{j\in\mathbb N}\|u_j\|_{X_j}^p,
		\qquad
		\|\mathbf v\|_Y^q:=\sum_{j\in\mathbb N}\|v_j\|_{Y_j}^q.
		\]
		Define $R_0u:=(u|_{B_j})_{j\in\mathbb N}$ and
		$\mathcal I\mathbf v:=\sum_{j\in\mathbb N}\chi_{Q_j}v_j$,
		where $v_j$ is extended by zero outside $Q_j$.
		The bounded overlap of the balls and the disjointness of the
		cubes give, respectively,
		\begin{equation}\label{eq:restriction-bound}
			\|R_0u\|_X^p
			=\int_{\C^n}|u|^pe^{-p\varphi}
			\sum_{j\in\mathbb N}\chi_{B_j}\,dv
			\le C\|u\|_{p,\varphi}^p
		\end{equation}
		and
		\begin{equation}\label{eq:assembly-isometry}
			\|\mathcal I\mathbf v\|_{q,\varphi}^q
			=\sum_{j\in\mathbb N}\|v_j\|_{Y_j}^q
			=\|\mathbf v\|_Y^q.
		\end{equation}
		
		For $j\in\mathbb N$, define $E_j:X_j\to Y_j$ by
		$E_ju:=(\vartheta u)|_{Q_j}$. Choose $M>0$ such that
		$\|\operatorname{Hess}_{\mathbb R}\varphi\|_{\rm op}\le M$, and put
		\[
		\gamma_j(z):=\exp\left(
		\varphi(a_j)+2\sum_{\nu=1}^n
		\partial_{z_\nu}\varphi(a_j)(z_\nu-a_{j,\nu})\right).
		\]
		With $\ell_j:=\log|\gamma_j|$, Taylor's theorem gives
		\begin{equation}\label{eq:gauge-remainder}
			|\varphi(z)-\ell_j(z)|
			\le \tfrac M2|z-a_j|^2\le \tfrac M2R_B^2,
			\qquad z\in B_j.
		\end{equation}
		Thus
		\begin{equation}\label{eq:gauge}
			e^{-MR_B^2/2}\le|\gamma_j(z)|e^{-\varphi(z)}
			\le e^{MR_B^2/2},\qquad z\in B_j.
		\end{equation}
		Choose $\delta<R_1<R_2<R_B$ and set $\theta:=\delta/R_1$.
		For $u\in X_j$, write $v:=u/\gamma_j$ and
		$v(a_j+w)=\sum_{k\ge0}H_{j,k}(w)$ for $|w|<R_B$,
		where $H_{j,k}$ is homogeneous of degree $k$.
		The submean inequality and \eqref{eq:gauge} give
		\begin{equation}\label{eq:taylor-submean}
			\sup_{B(a_j,R_1)}|v|
			\le C\|v\|_{L^p(B(a_j,R_2))}
			\le C\|u\|_{X_j}.
		\end{equation}
		By Cauchy's formula,
		\begin{equation}\label{eq:taylor-coefficients}
			|H_{j,k}(w)|\le C\|u\|_{X_j}(|w|/R_1)^k,
			\qquad |w|\le\delta,\quad k\ge0.
		\end{equation}
		For $L\ge0$, put
		$(\mathcal T_{j,L}u)(z):=\gamma_j(z)\sum_{k=0}^LH_{j,k}(z-a_j)$.
		Since $Q_j\subset\overline{B(a_j,\delta)}\subset B(a_j,R_1)$,
		equations \eqref{eq:gauge} and \eqref{eq:taylor-coefficients} give
		\[
		|u(z)-\mathcal T_{j,L}u(z)|e^{-\varphi(z)}
		\le C\|u\|_{X_j}\sum_{k>L}\theta^k
		\le C\theta^L\|u\|_{X_j},\qquad z\in Q_j.
		\]
		Also, \eqref{eq:gauge} and \eqref{eq:taylor-submean} give
		$\sup_{Q_j}|u|e^{-\varphi}\le C\|u\|_{X_j}$.
		Multiplying by $|\vartheta|$ and taking the $L^q(Q_j)$ norm yields
		\begin{equation}\label{eq:taylor}
			\|E_j-E_j\mathcal T_{j,L}\|\le Cb_j\theta^L,
			\qquad \|E_j\|\le Cb_j.
		\end{equation}
		Let $\mathcal P_L$ denote the space of polynomials in $n$ variables
		of total degree at most $L$. Then
		\[
		\operatorname{Ran}\mathcal T_{j,L}
		=\{\gamma_j P:P\in\mathcal P_L\}.
		\]
		Consequently,
		\begin{equation}\label{rankDL}
			\rank(E_j\mathcal T_{j,L})
			\le\dim\mathcal P_L=\binom{n+L}{n}.
		\end{equation}
		
		Set $E\mathbf u:=(E_ju_j)_{j\in\mathbb N}$. By \eqref{eq:taylor},
		\begin{equation}\label{eq:local-multiplication-bound}
			\begin{aligned}
				\|E\mathbf u\|_Y
				\le C\left(\sum_{j\in\mathbb N}b_j^q\|u_j\|_{X_j}^q\right)^{1/q}
				=C\left\|D_{\mathbf b}\bigl((\|u_j\|_{X_j})_{j\in\mathbb N}\bigr)\right\|_{\ell^q}
				\le C\|D_{\mathbf b}\|\|\mathbf u\|_X.
			\end{aligned}
		\end{equation}
		Thus $E:X\to Y$ is bounded, and the definitions give
		\begin{equation}\label{eq:factor}
			J_\vartheta=\mathcal IER_0.
		\end{equation}
		
	Let $C_s$ be the constant in Lemma~\ref{lem:ideal-properties}
	and put $C_\Psi:=C_0 2^\gamma$, where $C_0$ and $\gamma$
	are the constants in the growth condition on $\Psi$.
		Fix $*\in\{s,\Psi\}$ with $D_{\mathbf b}\in\mathfrak A_*$ and
		choose $K\in\mathbb N$ such that $C_*\theta^K<1/2$.
		For $j\in\mathbb N$ and $L\ge0$, set
		\[
		R_{j,L}:=E_j\mathcal T_{j,KL},\qquad
		R_{j,-1}:=0,\qquad S_{j,L}:=R_{j,L}-R_{j,L-1}.
		\]
		By the definition of the Taylor truncations and \eqref{rankDL},
		\begin{equation}\label{eq:pieces-rank}
			\rank S_{j,L}\le d_L:=\binom{n+KL}{n}\le C(L+1)^n,
		\end{equation}
		where $C$ is independent of $L$.
		For $L\ge1$, \eqref{eq:taylor} gives
		\[
		\|S_{j,L}\|\le\|E_j-R_{j,L}\|+\|E_j-R_{j,L-1}\|
		\le Cb_j\theta^{K(L-1)}.
		\]
		For $L=0$, it gives $\|S_{j,0}\|=\|R_{j,0}\|\le Cb_j$.
	By the definition of $S_{j,\ell}$,
	$
	\sum_{\ell=0}^L S_{j,\ell}
	=\sum_{\ell=0}^L(R_{j,\ell}-R_{j,\ell-1})
	=R_{j,L}.
	$
	Since $\|E_j-R_{j,L}\|\le Cb_j\theta^{KL}\to0$ by \eqref{eq:taylor},
	the series converges to $E_j$ in operator norm. Thus
	\begin{equation}\label{eq:pieces}
		\|S_{j,L}\|\le Cb_j\theta^{K(L-1)_+},
		\qquad E_j=\sum_{L\ge0}S_{j,L},
	\end{equation}
	where $(L-1)_+:=\max\{L-1,0\}$.
		
		Fix $L\ge0$ and, for $j\in\mathbb N$, put
		$V_{j,L}:=\operatorname{Ran}S_{j,L}$ and $r_{j,L}:=\dim V_{j,L}$.
		If $r_{j,L}>0$, fix coordinates on $V_{j,L}$ and choose
		$y_{j,1},\ldots,y_{j,r_{j,L}}\in V_{j,L}$ with
		$\|y_{j,\nu}\|_{Y_j}\le1$ such that
		\[
		|\det(y_{j,1},\ldots,y_{j,r_{j,L}})|
		=\max_{\substack{w_1,\ldots,w_{r_{j,L}}\in V_{j,L}\\
				\|w_\nu\|_{Y_j}\le1,\ 1\le\nu\le r_{j,L}}}
		|\det(w_1,\ldots,w_{r_{j,L}})|.
		\]
		Here $\det(y_{j,1},\ldots,y_{j,r_{j,L}})$ denotes the determinant
		of the matrix whose $\nu$th column is the coordinate vector of
		$y_{j,\nu}$ in the chosen basis of $V_{j,L}$.
		For its coordinate functionals
		$y_{j,\nu}^*$, Cramer's rule gives
		\begin{equation}\label{eq:auerbach-functional}
			|y_{j,\nu}^*(y)|
			=\frac{|\det(y_{j,1},\ldots,y,\ldots,y_{j,r_{j,L}})|}
			{|\det(y_{j,1},\ldots,y_{j,r_{j,L}})|}
			\le1,\qquad \|y\|_{Y_j}\le1,
		\end{equation}
		where $y$ replaces the $\nu$th column. Put
		$\phi_{j,\nu}:=y_{j,\nu}^*\circ S_{j,L}$ for $\nu\le r_{j,L}$,
		and set $y_{j,\nu}:=0$ and $\phi_{j,\nu}:=0$ for
		$r_{j,L}<\nu\le d_L$. Thus
		\begin{equation}\label{eq:auerbach}
			S_{j,L}u=\sum_{\nu=1}^{d_L}\phi_{j,\nu}(u)y_{j,\nu},
			\qquad u\in X_j,
		\end{equation}
		and \eqref{eq:auerbach-functional} and \eqref{eq:pieces} give
		\begin{equation}\label{eq:auerbach-bounds}
			\|\phi_{j,\nu}\|_{X_j^*}\le Cb_j\theta^{K(L-1)_+}.
		\end{equation}
		For $1\le\nu\le d_L$, define $F_{L,\nu}:X\to\ell^p$ and
		$Y_{L,\nu}:\ell^q\to Y$ by
		\[
		(F_{L,\nu}\mathbf u)_j:=
		\begin{cases}
			\phi_{j,\nu}(u_j)/b_j,&b_j>0,\\
			0,&b_j=0,
		\end{cases}
		\qquad Y_{L,\nu}\mathbf c:=(c_jy_{j,\nu})_{j\in\mathbb N}.
		\]
		By \eqref{eq:auerbach-bounds} and $\|y_{j,\nu}\|_{Y_j}\le1$ for all $j$,
		\begin{equation}\label{eq:coefficient-map-bound}
			\|F_{L,\nu}\|\le C\theta^{K(L-1)_+},\qquad 	\|Y_{L,\nu}\|\le1.
		\end{equation}
		Set $S_L\mathbf u:=(S_{j,L}u_j)_{j\in\mathbb N}$.
		Since $S_{j,L}=0$ when $b_j=0$, \eqref{eq:auerbach} gives
		\begin{equation}\label{eq:finite-factor}
			S_L=\sum_{\nu=1}^{d_L}Y_{L,\nu}D_{\mathbf b}F_{L,\nu}.
		\end{equation}
		
		For $k\in\mathbb N$ and $\varepsilon>0$, choose $A:\ell^p\to\ell^q$
		with $\rank A<k$ and
		$\|D_{\mathbf b}-A\|\le a_k(D_{\mathbf b})+\varepsilon$.
		Put $A_L:=\sum_{\nu=1}^{d_L}Y_{L,\nu}AF_{L,\nu}$. Then
		\begin{equation}\label{eq:piece-rank}
			\rank A_L\le d_L(k-1).
		\end{equation}
		By \eqref{eq:coefficient-map-bound} and \eqref{eq:finite-factor},
		\begin{equation}\label{eq:piece-error}
			\|S_L-A_L\|
			\le Cd_L\theta^{K(L-1)_+}\bigl(a_k(D_{\mathbf b})+\varepsilon\bigr).
		\end{equation}
		Letting $\varepsilon\downarrow0$ yields
		\begin{equation}\label{eq:piece-ideal}
			a_{d_L(k-1)+1}(S_L)
			\le Cd_L\theta^{K(L-1)_+}a_k(D_{\mathbf b}).
		\end{equation}
		For $m\in\bigl(d_L(k-1),d_Lk\bigr]\cap\mathbb N$, monotonicity gives
		$a_m(S_L)\le a_{d_L(k-1)+1}(S_L)$. Summing over these blocks,
		\begin{equation}\label{eq:piece-small}
			\As{S_L}^s
			\le d_L\sum_{k\ge1}a_{d_L(k-1)+1}(S_L)^s
			\le C^sd_L^{s+1}\theta^{sK(L-1)_+}\As{D_{\mathbf b}}^s.
		\end{equation}
		On the same blocks, \eqref{eq:rates-growth} gives
		$\Psi(m)\le\Psi(d_Lk)\le C_0d_L^\gamma\Psi(k)$.
		By \eqref{eq:piece-ideal},
		\[
		\|S_L\|_{\mathfrak A_\Psi}
		=\sup_{m\ge1}\Psi(m)a_m(S_L)
		\le C d_L^{1+\gamma}\theta^{K(L-1)_+}
		\sup_{k\ge1}\Psi(k)a_k(D_{\mathbf b})
		= C d_L^{1+\gamma}\theta^{K(L-1)_+}
		\|D_{\mathbf b}\|_{\mathfrak A_\Psi}.
		\] Together with
		\eqref{eq:piece-small}, this yields
		\begin{equation}\label{eq:piece-ideal-norm}
			\|S_L\|_{\mathfrak A_*}
			\le Cd_L^{1+\max\{1/s,\gamma\}}
			\theta^{K(L-1)_+}\|D_{\mathbf b}\|_{\mathfrak A_*}.
		\end{equation}
		
		For $N\ge0$, set $S^{(N)}:=\sum_{L=0}^N S_L$. Applying the
		quasi-triangle inequality from the right and using
		\eqref{eq:pieces-rank} and \eqref{eq:piece-ideal-norm}, we obtain
		\begin{equation}\label{eq:geometric}
			\|S^{(N)}\|_{\mathfrak A_*}
			\le C\|D_{\mathbf b}\|_{\mathfrak A_*}
			\sum_{L\ge0}C_*^{L+1}(L+1)^{n(1+\max\{1/s,\gamma\})}
			\theta^{K(L-1)_+}
			\le C\|D_{\mathbf b}\|_{\mathfrak A_*},
		\end{equation}
		since $C_*\theta^K<1/2$.
		 Since $R_{j,-1}=0$ for
		$j\in\mathbb N$, we have
		\[
		S^{(N)}\mathbf u
		=\left(\sum_{L=0}^N S_{j,L}u_j\right)_{j\in\mathbb N}
		=\left(\sum_{L=0}^N(R_{j,L}-R_{j,L-1})u_j\right)_{j\in\mathbb N}
		=(R_{j,N}u_j)_{j\in\mathbb N}.
		\]
		Using \eqref{eq:taylor} as in \eqref{eq:local-multiplication-bound},
		\begin{equation}\label{eq:transfer-convergence}
		\|(E-S^{(N)})\mathbf u\|_Y
		=\left(\sum_{j\in\mathbb N}
		\|(E_j-R_{j,N})u_j\|_{Y_j}^q\right)^{1/q}
		\le C\theta^{KN}\left(\sum_{j\in\mathbb N}
		b_j^q\|u_j\|_{X_j}^q\right)^{1/q}
		\le C\theta^{KN}\|D_{\mathbf b}\|\|\mathbf u\|_X.
		\end{equation}
		In particular,
		\begin{equation}\label{eq:transfer-operator-limit}
			\|E-S^{(N)}\|\le C\theta^{KN}\|D_{\mathbf b}\|\longrightarrow0,
\quad(N\to\infty).
		\end{equation}
		By Lemma~\ref{lem:ideal-properties}\textup{(iv)}, 
		\begin{equation}\label{eq:transfer-number-limit}
			|a_k(E)-a_k(S^{(N)})|\le\|E-S^{(N)}\|\longrightarrow0,
			\qquad(N\to\infty).
		\end{equation}
		Fatou's lemma for $*=s$ and \eqref{eq:rates-Fatou} for $*=\Psi$,
		together with \eqref{eq:geometric}, give
		\begin{equation}\label{eq:transfer-ideal-limit}
			\|E\|_{\mathfrak A_*}
			\le\liminf_{N\to\infty}\|S^{(N)}\|_{\mathfrak A_*}
			\le C\|D_{\mathbf b}\|_{\mathfrak A_*}.
		\end{equation}
		Finally, the factorization \eqref{eq:factor}, the bounds
		\eqref{eq:restriction-bound} and \eqref{eq:assembly-isometry},
		and the ideal property give
		\[
		\|J_\vartheta\|_{\mathfrak A_*}
		\le\|\mathcal I\|\|R_0\|\|E\|_{\mathfrak A_*}
		\le C\|D_{\mathbf b}\|_{\mathfrak A_*}.
		\]
		This completes the proof.
	\end{proof}

\section{Approximation numbers and ideals of Hankel operators}\label{sec:Hankel}
In this section, we establish the comparison between Hankel operators
and diagonal operators and use it to prove the main results on
approximation-number decay and ideal membership.
We also apply the local recovery argument to weighted Bergman spaces.

\subsection{Decay estimates}

\begin{proposition}\label{prop:bounded-diagonal-recovery}
	\label{thm:rates-diagonal-transfer}
	Suppose $\operatorname{Hess}_{\mathbb R}\varphi\simeq E$.
	Let $0<p<\infty$, $1\le q<\infty$, and
	$f\in\mathcal S\cap L^q_{\rm loc}(\C^n)$.
	Let $\Psi$ be nondecreasing and satisfy \eqref{eq:rates-growth}.
	For every $\Lambda=\{a_j\}_{j\in\mathbb N}\in\mathscr G$,
	let $R$ be defined by \eqref{eq:parameters} and put
	$\mathbf g:=(G_{q,R}(f)(a_j))_{j\in\mathbb N}$.
	Then $H_f\in\mathfrak A_\Psi(F^p(\varphi),L^q(\varphi))$
	if and only if
	$D_{\mathbf g}\in\mathfrak A_\Psi(\ell^p,\ell^q)$.
	Moreover,
	\begin{equation}\label{eq:rates-diagonal-equivalence}
		\|H_f\|_{\mathfrak A_\Psi}
		\simeq\|D_{\mathbf g}\|_{\mathfrak A_\Psi}.
	\end{equation}
	The comparison constants are independent of $f$.
\end{proposition}

\begin{proof}
	First suppose $D_{\mathbf g}\in\mathfrak A_\Psi(\ell^p,\ell^q)$.
	Set $g_j:=G_{q,R}(f)(a_j)$ for all $j\in\mathbb N$.
	Let $f=f_1+f_2$ be the decomposition in
	Lemma~\ref{lem:known-decomposition}.
	By \eqref{eq:parameters}, $Q_j\subset B(a_j,t/2)$ and
	$B(z,2t)\subset B(a_j,R)$ for $z\in Q_j$. Consequently,
	\begin{equation}\label{eq:local-G}
		G_{q,2t}(f)(z)
		\le\left(\frac{R}{2t}\right)^{2n/q}g_j,
		\qquad z\in Q_j.
	\end{equation}
	For a locally $q$-integrable function $u$, put
	$b_j(u):=(\int_{Q_j}|u|^q\,dv)^{1/q}$.
	By \eqref{eq:decomposition} and \eqref{eq:local-G},
	\[
	b_j(f_2)
	\le v(B(0,t/2))^{1/q}M_{q,t/2}(f_2)(a_j)
	\le Cg_j,
	\]
	and
	\begin{equation}\label{eq:local-b}
		\sup_{z\in Q_j}|\bar\partial f_1(z)|\le Cg_j.
	\end{equation}
	Thus $b_j(\bar\partial f_1)\le Cg_j$.
	Lemma~\ref{lem:diag-facts} and
	Lemma~\ref{lem:transfer} give
	\begin{equation}\label{eq:J}
		\|J_{f_2}\|_{\mathfrak A_\Psi}
		+\|J_{\bar\partial f_1}\|_{\mathfrak A_\Psi}
		\le C\|D_{\mathbf g}\|_{\mathfrak A_\Psi}.
	\end{equation}
	Since $\|\mathbf g\|_{\ell^\infty}\le\|D_{\mathbf g}\|$,
	\eqref{eq:local-b} implies
	$\bar\partial f_1\in L^\infty$.
	If $p>q$, put $r:=pq/(p-q)>1$.
	Lemma~\ref{lem:diag-facts}\textup{(ii)} gives
	\begin{equation}\label{eq:derivative-class}
		\|\bar\partial f_1\|_{L^r}^r
		=\sum_{j\in\mathbb N}\int_{Q_j}|\bar\partial f_1|^r\,dv
		\le Ch^{2n}\sum_{j\in\mathbb N}g_j^r
		=Ch^{2n}\|D_{\mathbf g}\|^r<\infty.
	\end{equation}
	Moreover, \eqref{eq:kernel-upper},
	\eqref{eq:normalized-kernel}, and H\"older's inequality give
	\begin{align}\label{eq:symbol-check}
		\int_{\C^n}|f_2k_a|e^{-\varphi}\,dv
		&\le C\sum_{j\ge1}e^{-c|a_j-a|}
		\int_{Q_j}|f_2|\,dv\notag\\
		&\le C\sum_{j\ge1}e^{-c|a_j-a|}
		v(Q_j)^{1-1/q}b_j(f_2)\notag\\
		&\le C\|\mathbf g\|_{\ell^\infty}
		\sum_{j\ge1}e^{-c|a_j-a|}<\infty,
		\qquad a\in\C^n.
	\end{align}
	Hence $f_2,f_1=f-f_2\in\mathcal S$.
	By \eqref{eq:Aidentity},
	\begin{equation}\label{eq:hankel-factor}
		H_f=(I-P)\bigl(J_{f_2}
		+A_\varphi J_{\bar\partial f_1}\bigr)
		\quad\text{on }\Gamma.
	\end{equation}
	Equations \eqref{eq:A}, \eqref{eq:J}, and
	Lemma~\ref{lem:rates-ideal-properties} yield
	\begin{equation}\label{eq:rates-upper-transfer}
		\|H_f\|_{\mathfrak A_\Psi}
		\le C\bigl(
		\|J_{f_2}\|_{\mathfrak A_\Psi}
		+\|J_{\bar\partial f_1}\|_{\mathfrak A_\Psi}\bigr)
		\le C\|D_{\mathbf g}\|_{\mathfrak A_\Psi}.
	\end{equation}
	
	Conversely, suppose
	$H_f\in\mathfrak A_\Psi(F^p(\varphi),L^q(\varphi))$.
	Fix a grid point $a_*$ and write
	\[
	m_j:=\frac{a_j-a_*}{h}\in\mathbb Z^{2n},
	\qquad j\in\mathbb N.
	\]
	For $L_0\in\mathbb N$ to be chosen below, partition $\mathbb N$ into
	\[
	\Lambda_\nu:=\{j\in\mathbb N:m_j-\nu\in L_0\mathbb Z^{2n}\},
	\qquad \nu\in\{0,\ldots,L_0-1\}^{2n}.
	\]
	Set $B_j:=B(a_j,R)$. For distinct $j,k\in\Lambda_\nu$,
	\[
	m_j-m_k\in L_0\bigl(\mathbb Z^{2n}\setminus\{0\}\bigr),
	\qquad |a_j-a_k|=h|m_j-m_k|\ge hL_0.
	\]
	Thus $B_j\cap B_k=\varnothing$ whenever $hL_0>2R$.
	By \eqref{eq:kernel-upper}, \eqref{eq:kernel-lower}, and
	\eqref{eq:normalized-kernel},
	\begin{equation}\label{eq:local-distance}
		d_j:=\inf_{u\in\mathcal H(B_j)}
		\left(\int_{B_j}|f-u|^q|k_{a_j}|^qe^{-q\varphi}\,dv\right)^{1/q}
		\simeq g_j,\qquad j\in\mathbb N.
	\end{equation}
	Moreover, since $R<r_0$ by \eqref{eq:parameters},
	the same kernel estimates imply that $k_{a_j}$ is
	nonvanishing on $B_j$ and
	\begin{equation}\label{eq:local-kernel-ratio}
		\sup_{w\in B_j}
		\left|\frac{k_{a_k}(w)}{k_{a_j}(w)}\right|
		\lesssim\sup_{w\in B_j}e^{-\theta|w-a_k|}
		\le e^{\theta R}e^{-\theta|a_j-a_k|}.
	\end{equation}
	Multiplication by $k_{a_j}$
	identifies the holomorphic approximants in \eqref{eq:local-distance}
	with those to $fk_{a_j}$ in
	$\mathcal H_j:=\mathcal H(B_j)\cap L^q(B_j,e^{-q\varphi}dv)$.
	Thus $d_j=\operatorname{dist}(fk_{a_j},\mathcal H_j)$ in the weighted
	$L^q$ norm. Hahn--Banach gives
	$\lambda_j\in L^q(B_j,e^{-q\varphi}dv)^*$ with
	\begin{equation}\label{eq:lambda}
		\|\lambda_j\|\le1,\qquad
		\lambda_j|_{\mathcal H_j}=0,\qquad
		\lambda_j(fk_{a_j})=d_j.
	\end{equation}
	Take $\lambda_j:=0$ when $d_j=0$.
	
	For each class $\Lambda_\nu$,
	define $U_\nu:\ell^p(\Lambda_\nu)\to F^p(\varphi)$
	and $V_\nu:L^q(\varphi)\to\ell^q(\Lambda_\nu)$ by
	\begin{equation}\label{eq:UV}
		U_\nu\mathbf c:=\sum_{j\in\Lambda_\nu}c_jk_{a_j},
		\qquad (V_\nu u)_j:=\lambda_j(u|_{B_j}).
	\end{equation}
	For finitely supported $\mathbf c$, \eqref{eq:UV},
	\eqref{eq:kernel-upper}, and \eqref{eq:normalized-kernel} give
	\[
	|U_\nu\mathbf c(z)|e^{-\varphi(z)}
	\le C\sum_{j\in\Lambda_\nu}|c_j|e^{-c|z-a_j|}.
	\]
	Set $\alpha:=\min\{1,p\}$.
	Using the uniform lattice sum bound from \eqref{eq:symbol-check},
	H\"older's inequality for $p>1$ and subadditivity for $0<p\le1$ give
	\[
	|U_\nu\mathbf c(z)|^pe^{-p\varphi(z)}
	\le C\sum_{j\in\Lambda_\nu}|c_j|^pe^{-c\alpha|z-a_j|}.
	\]
	Integrating this estimate and using
	\eqref{eq:lambda} and disjointness for $V_\nu$, we get
	\begin{equation}\label{eq:UVnorm}
		\|U_\nu\mathbf c\|_{p,\varphi}^p\le C\|\mathbf c\|_{\ell^p}^p,
		\qquad
		\|V_\nu u\|_{\ell^q}^q
		\le\sum_{j\in\Lambda_\nu}\int_{B_j}|u|^qe^{-q\varphi}\,dv
		\le\|u\|_{q,\varphi}^q.
	\end{equation}
	By \eqref{eq:UVnorm} and the density of finitely supported sequences
	in $\ell^p(\Lambda_\nu)$, $U_\nu$ extends uniquely to a bounded
	operator on $\ell^p(\Lambda_\nu)$. Since the finite truncations of
	$\mathbf c$ converge to $\mathbf c$ in $\ell^p(\Lambda_\nu)$,
	\[
	U_\nu\mathbf c
	=\lim_{N\to\infty}
	\sum_{\substack{j\in\Lambda_\nu\\j\le N}}c_jk_{a_j}
	\qquad\text{in }F^p(\varphi).
	\]
	Set $M_\nu:=V_\nu H_fU_\nu$. Since
	$(fk_{a_k}-H_fk_{a_k})|_{B_j}\in\mathcal H_j$, we have
	\begin{equation}\label{eq:matrix}
		(M_\nu)_{jk}=\lambda_j(fk_{a_k}),
		\qquad (M_\nu)_{jj}=d_j.
	\end{equation}
	If $u\in\mathcal H(B_j)$ and
	$(f-u)k_{a_j}\in L^q(B_j,e^{-q\varphi}dv)$,
	then \eqref{eq:local-kernel-ratio} gives
	$uk_{a_k}=(uk_{a_j})(k_{a_k}/k_{a_j})\in\mathcal H_j$.
	Hence \eqref{eq:lambda} implies
	\[
	|(M_\nu)_{jk}|
	=|\lambda_j((f-u)k_{a_k})|
	\le C_Re^{-\theta|a_j-a_k|}
	\|(f-u)k_{a_j}\|_{L^q(B_j,e^{-q\varphi}dv)}.
	\]
	Taking the infimum in \eqref{eq:local-distance} yields
	\begin{equation}\label{eq:offdiag}
		|(M_\nu)_{jk}|
		\le C_Rd_je^{-\theta|a_j-a_k|},
		\qquad j\ne k.
	\end{equation}
	
	Define $B^{(\nu)}_{jk}:=(M_\nu)_{jk}/d_j$
	if $j\ne k$ and $d_j>0$,
	and $B^{(\nu)}_{jk}:=0$ otherwise.
	For each $m\in\mathbb Z^{2n}\setminus\{0\}$,
	the entries with $a_j-a_k=L_0hm$
	form a weighted shift on $\ell^p(\Lambda_\nu)$.
	Thus the $\alpha$-triangle inequality and \eqref{eq:offdiag} give
	\begin{equation}\label{eq:Schur}
		\|B^{(\nu)}\|_{\ell^p\to\ell^p}^{\alpha}
		\le\sum_{m\in\mathbb Z^{2n}\setminus\{0\}}
		\sup_{\substack{j,k\in\Lambda_\nu\\a_j-a_k=L_0hm}}
		|B^{(\nu)}_{jk}|^\alpha
		\le C\sum_{m\in\mathbb Z^{2n}\setminus\{0\}}
		e^{-c\alpha L_0h|m|}.
	\end{equation}
	Since $|m|\ge1$ for $m\in\mathbb Z^{2n}\setminus\{0\}$,
	the right-hand side of \eqref{eq:Schur} is bounded by
	\[
	Ce^{-c\alpha h(L_0-1)}
	\sum_{m\in\mathbb Z^{2n}\setminus\{0\}}e^{-c\alpha h|m|}
	\longrightarrow0\qquad(L_0\to\infty).
	\]
	Choose $L_0$ so that $L_0h>2R$ and
	$\|B^{(\nu)}\|_{\ell^p\to\ell^p}^{\alpha}\le\tfrac12$
	for every $\nu$.
	The Neumann series
	$(I+B^{(\nu)})^{-1}
	=\sum_{\ell=0}^{\infty}(-B^{(\nu)})^\ell$
	converges in operator norm, and the $\alpha$-triangle inequality gives
	\[
	\|(I+B^{(\nu)})^{-1}\|_{\ell^p\to\ell^p}^{\alpha}
	\le\sum_{\ell=0}^{\infty}
	\|B^{(\nu)}\|_{\ell^p\to\ell^p}^{\ell\alpha}
	\le\sum_{\ell=0}^{\infty}2^{-\ell}=2.
	\]
	For finitely supported $\mathbf c$,
	\eqref{eq:matrix} and the definition of $B^{(\nu)}$ give
	\begin{equation}\label{eq:MB}
		(M_\nu\mathbf c)_j
		=d_jc_j+d_j\sum_{\substack{k\in\Lambda_\nu\\k\ne j}}
		B^{(\nu)}_{jk}c_k
		=d_j\bigl((I+B^{(\nu)})\mathbf c\bigr)_j,
		\qquad j\in\Lambda_\nu.
	\end{equation}
	By density, \eqref{eq:MB} holds coordinatewise for all
	$\mathbf c\in\ell^p(\Lambda_\nu)$.
	Put $\mathbf d:=(d_j)_{j\in\mathbb N}$ and write
	$D_{\mathbf d|_{\Lambda_\nu}}\mathbf c
	:=(d_jc_j)_{j\in\Lambda_\nu}$.
	Applying it to $(I+B^{(\nu)})^{-1}\mathbf c$ shows that
	$D_{\mathbf d|_{\Lambda_\nu}}\mathbf c\in\ell^q(\Lambda_\nu)$
	and gives
	\begin{equation}\label{eq:recover}
		D_{\mathbf d|_{\Lambda_\nu}}
		=V_\nu H_fU_\nu(I+B^{(\nu)})^{-1}.
	\end{equation}
	By \eqref{eq:local-distance}, \eqref{eq:UVnorm},
	and the ideal property,
	\[
	a_k(D_{\mathbf g|_{\Lambda_\nu}})
	\le Ca_k(D_{\mathbf d|_{\Lambda_\nu}})
	\le Ca_k(H_f),\qquad k\in\mathbb N.
	\]
	Set $N:=L_0^{2n}$. By \eqref{Pi},
	\begin{equation}\label{eq:colors}
		D_{\mathbf g}
		=\sum_\nu\iota_{\Lambda_\nu}
		D_{\mathbf g|_{\Lambda_\nu}}P_{\Lambda_\nu}.
	\end{equation}
	Lemma~\ref{lem:ideal-properties}\textup{(ii)} gives
	\begin{equation}\label{eq:rates-recovery-numbers}
		a_{N(k-1)+1}(D_{\mathbf g})\le Ca_k(H_f),
		\qquad k\in\mathbb N.
	\end{equation}
	For $m\in\bigl(N(k-1),Nk\bigr]\cap\mathbb N$, monotonicity,
	\eqref{eq:rates-growth}, and \eqref{eq:rates-recovery-numbers} give
	\[
	\Psi(m)a_m(D_{\mathbf g})
	\le C_0N^\gamma\Psi(k)a_{N(k-1)+1}(D_{\mathbf g})
	\le C\Psi(k)a_k(H_f).
	\]
	Taking the supremum over $m\in\mathbb N$ yields
	$\|D_{\mathbf g}\|_{\mathfrak A_\Psi}
	\le C\|H_f\|_{\mathfrak A_\Psi}$.
	Together with \eqref{eq:rates-upper-transfer}, this proves
	\eqref{eq:rates-diagonal-equivalence}.
	This completes the proof.
\end{proof}

\begin{proof}[Proof of Theorem~\ref{thm:rates-rearranged-decay}]
	Fix $\Lambda=\{a_j\}_{j\in\mathbb N}\in\mathscr G$ and put
	\[
	g_j:=G_{q,R}(f)(a_j),\qquad
	\mathbf g:=(g_j)_{j\in\mathbb N},\qquad
	F:=G_{q,1}(f).
	\]
Define $B(t):=g_j^*$ for $j-1\le t<j$, $t>0$,
and $j\in\mathbb N$.
	The assumption on the symbol $f$ ensures that $H_f$ is compact,
	so $a_k(H_f)\to0$ by the approximation property of $L^q(\varphi)$.
	By \eqref{eq:rates-recovery-numbers}, $D_{\mathbf g}$ is compact
	and hence $\mathbf g\in c_0$.
By Lemma~\ref{lem:indexed-sampling}, there exists a fixed $c\ge1$
such that
\[
F^*(ct)\lesssim B(t),\qquad
B(ct)\lesssim F^*(t),\qquad t>0,
\]
with implicit constants independent of $f$ and $t$.
Letting $t\downarrow0$ in the first inequality gives
$\|F\|_{L^\infty}\lesssim g_1^*$.
For each fixed $d>0$, monotonicity of $\Psi$ and
\eqref{eq:rates-growth}, applied to $\min\{t,dt\}$ and
$\max\{t,dt\}$, yield
\[
\Psi(dt)\simeq\Psi(t),
\qquad t\ge\max\{1,d^{-1}\}.
\]
	
	Fix $k_0\ge2$. By monotonicity of $\Psi$ and
	$a_1(H_f)=\|H_f\|$, the contribution of the indices $k<k_0$ satisfies
	\[
	\sup_{1\le k<k_0}\Psi(k)F^*(k)
	\le \Psi(k_0)\|F\|_\infty,
	\qquad
	\sup_{1\le k<k_0}\Psi(k)a_k(H_f)
	\le \Psi(k_0)\|H_f\|.
	\]
	Thus, to determine whether these weighted suprema are finite,
	it suffices to consider $k\ge k_0$.
	
	\textup{(i)}
	By \eqref{eq:same-side} and Proposition~\ref{thm:rates-diagonal-transfer},
	\[
	\|H_f\|_{\mathfrak A_\Psi}\simeq\sup_{k\ge1}\Psi(k)g_k^*.
	\]
	For $k\ge4c$, set $m:=\lfloor k/(2c)\rfloor$. Then
	$k/(4c)\le m$ and $cm\le k-\tfrac12$. By monotonicity and
	\eqref{eq:rates-indexed-sampling},
	\[
	F^*(k)\le F^*(cm)\lesssim B(m)\le g_m^*,
	\qquad
	g_k^*=B(k-\tfrac12)\le B(cm)\lesssim F^*(m).
	\]
	Since $\Psi(k)\lesssim\Psi(m)$ by \eqref{eq:rates-growth}, it follows that
	\[
	\Psi(k)F^*(k)\lesssim\Psi(m)g_m^*,
	\qquad
	\Psi(k)g_k^*\lesssim\Psi(m)F^*(m).
	\]
	Taking the supremum over integers $k\ge4c$ in both inequalities gives
	\[
	\sup_{k\ge4c}\Psi(k)F^*(k)
	\lesssim \sup_{m\ge1}\Psi(m)g_m^*,
	\qquad
	\sup_{k\ge4c}\Psi(k)g_k^*
	\lesssim \sup_{m\ge1}\Psi(m)F^*(m).
	\]
	By the preceding reduction to large indices and the norm equivalence
	for $H_f$, we obtain
	\[
	a_k(H_f)\lesssim\Psi(k)^{-1}
	\quad\Longleftrightarrow\quad
	F^*(k)\lesssim\Psi(k)^{-1}.
	\]
	This completes the proof of \textup{(i)}.
	
	\textup{(ii)}
	By Lemma~\ref{lem:rates-crossing-window} and
	Proposition~\ref{thm:rates-diagonal-transfer},
	\begin{equation}\label{eq:rates-crossing}
		\|H_f\|_{\mathfrak A_\Psi}
		\simeq\sup_{k\ge1}\Psi(k)
		W_{\tau,\varsigma}(\mathbf g;k).
	\end{equation}
	Put $a:=\tau/2>1$ and define
	\[
	I_\phi(t):=
	\int_t^{t^a}
	u^{\varsigma/\tau-1}\phi(u)^\varsigma\,du.
	\]
	Since $c\ge1$ and $a>1$, for $t\ge c^{1+1/a}$ we have
	\[
	[t/c,t^a/c]
	=[t/c,(t/c)^a]\cup[t/c^{1/a},(t/c^{1/a})^a].
	\]
	Thus \eqref{eq:rates-indexed-sampling} gives
	\[
	I_{F^*}(t)
	\le C\bigl(I_B(t/c)+I_B(t/c^{1/a})\bigr).
	\]
	The reverse comparison holds with $B$ and $F^*$
	interchanged.
	By \eqref{eq:rates-growth},
	\begin{equation}\label{eq:rates-window-continuous-transfer}
		I_{F^*}(t)\lesssim
		t^{\varsigma/2}\Psi(t)^{-\varsigma}
		\quad\Longleftrightarrow\quad
		I_B(t)\lesssim
		t^{\varsigma/2}\Psi(t)^{-\varsigma}.
	\end{equation}
	
For $t>1$, set $k:=\lfloor t\rfloor$. If
$(j-1,j)\cap(t,t^{\tau/2})\ne\varnothing$, then $j>t$ and
$j-1<t^{\tau/2}$. Since $k\le t<k+1$, it follows that
\[
k+1\le j\le\lceil t^{\tau/2}\rceil
\le\lceil(k+1)^{\tau/2}\rceil.
\]
For each such $j$ and every integration variable $u\in(j-1,j)$,
we have $B(u)=g_j^*$ and $j/2\le u<j$. Splitting the integral
defining $I_B(t)$ over the intersections
$(j-1,j)\cap(t,t^{\tau/2})$ and using
\eqref{eq:rates-window}, we obtain
\[
I_B(t)
\lesssim \sum_{j=k+1}^{\lceil(k+1)^{\tau/2}\rceil}
j^{\varsigma/\tau-1}(g_j^*)^\varsigma
=(k+1)^{\varsigma/2}W_{\tau,\varsigma}(\mathbf g;k+1)^\varsigma
\lesssim k^{\varsigma/2}W_{\tau,\varsigma}(\mathbf g;k+1)^\varsigma.
\]
	
	Conversely, for all sufficiently large $k$,
	\[
	[k-1,\lceil k^a\rceil]\subset
	[k-1,(k-1)^a]\cup[k+1,(k+1)^a],
	\]
and hence
	\[
	k^{\varsigma/2}
	W_{\tau,\varsigma}(\mathbf g;k)^\varsigma
	\le C\bigl(I_B(k-1)+I_B(k+1)\bigr).
	\]
The finitely many remaining terms are controlled by $g_1^*$.
	Together with \eqref{eq:rates-crossing},
	\eqref{eq:rates-window-continuous-transfer}, and
	\eqref{eq:rates-growth}, these estimates yield
	\begin{equation}\label{eq:rates-window-real-criterion}
		a_k(H_f)\lesssim\Psi(k)^{-1}
		\quad\Longleftrightarrow\quad
		I_{F^*}(t)\lesssim
		t^{\varsigma/2}\Psi(t)^{-\varsigma},
		\qquad t\ge2.
	\end{equation}
	
To finish the proof of \textup{(ii)}, we show that the bound
$I_{F^*}(t)\lesssim t^{\varsigma/2}\Psi(t)^{-\varsigma}$
for all sufficiently large real $t$ is equivalent to the same
bound at integers $t=k+1$. The forward implication follows by
setting $t=k+1$, so it suffices to prove the converse.

Assume the bound holds at all sufficiently large integers and set
$m:=\lfloor t\rfloor$. Since $m\le t<m+1$ and $m^a\ge m+1$ for
sufficiently large $m$, we have
\begin{equation}\label{eq:rates-integer-windows}
	[t,t^a]\subseteq[m,m^a]\cup[m+1,(m+1)^a].
\end{equation}
By \eqref{eq:rates-integer-windows} and the integer bounds,
\[
I_{F^*}(t)\le I_{F^*}(m)+I_{F^*}(m+1)
\lesssim m^{\varsigma/2}\Psi(m)^{-\varsigma}
+(m+1)^{\varsigma/2}\Psi(m+1)^{-\varsigma}
\lesssim t^{\varsigma/2}\Psi(t)^{-\varsigma}.
\]
The last step follows from monotonicity of $\Psi$ and
\eqref{eq:rates-growth}. Combining this equivalence with \eqref{eq:rates-window-real-criterion}, taking the $\varsigma$th root at $t=k+1$, and using
$\sqrt{k+1}/\Psi(k+1)\simeq\sqrt{k}/\Psi(k)$, we obtain
\[
a_k(H_f)\lesssim\Psi(k)^{-1}
\quad\Longleftrightarrow\quad
\left(
\int_{k+1}^{(k+1)^{\tau/2}}
u^{\varsigma/\tau-1}F^*(u)^\varsigma\,du
\right)^{1/\varsigma}
\lesssim\frac{\sqrt{k}}{\Psi(k)}.
\]
This completes the proof of \textup{(ii)}.
	
	\textup{(iii)}
First, we compare the discrete and continuous
$L^r$ tails. Since $D_{\mathbf g}$ is bounded,
Lemma~\ref{lem:diag-facts}\textup{(ii)} gives $\mathbf g\in\ell^r$.
Then \eqref{eq:power-sampling} gives $F\in L^r$.
By Lemma~\ref{lem:tail} and
Proposition~\ref{thm:rates-diagonal-transfer},
\begin{equation}\label{eq:rates-tail}
	\|H_f\|_{\mathfrak A_\Psi}
	\simeq \sup_{k\ge1}\Psi(k)
	\left(\sum_{j\ge k}(g_j^*)^r\right)^{1/r}.
\end{equation}
For $\phi=F^*$ or $B$, set
\[
J_\phi(t):=\left(\int_t^\infty\phi(u)^r\,du\right)^{1/r},
\qquad t\ge0.
\]
Integrating the $r$th powers in \eqref{eq:rates-indexed-sampling}
and using the change of variables $u=cv$ gives
\[
J_{F^*}(ct)\lesssim J_B(t),
\qquad
J_B(ct)\lesssim J_{F^*}(t),
\qquad t>0.
\]
Since $B(u)=g_j^*$ for $j-1<u<j$, we also have
\[
J_B(k-1)^r=\sum_{j\ge k}(g_j^*)^r,
\qquad k\in\mathbb N.
\]

For $k\ge4c$, set $m:=\lfloor k/(2c)\rfloor$. Then
$m\ge k/(4c)$ and $cm\le k-1$. By monotonicity of the tails,
\[
J_{F^*}(k)\le J_{F^*}(cm)\lesssim J_B(m)\le J_B(m-1),
\qquad
J_B(k-1)\le J_B(cm)\lesssim J_{F^*}(m).
\]
By \eqref{eq:rates-growth}, $\Psi(k)\lesssim\Psi(m)$.
Multiplying by $\Psi(k)$ and taking suprema therefore shows that
\[
\sup_{k\ge1}\Psi(k)J_{F^*}(k)<\infty
\quad\Longleftrightarrow\quad
\sup_{k\ge1}\Psi(k)J_B(k-1)<\infty,
\]
where the finitely many smaller indices give finite contributions
because $F\in L^r$ and $\mathbf g\in\ell^r$.
Combining this with \eqref{eq:rates-tail} yields
\[
a_k(H_f)\lesssim\Psi(k)^{-1}
\quad\Longleftrightarrow\quad
\left(\int_k^\infty F^*(u)^r\,du\right)^{1/r}
\lesssim\Psi(k)^{-1}.
\]
This completes the proof.
\end{proof}

The next corollary reduces the integral conditions in
Theorem~\ref{thm:rates-rearranged-decay} to pointwise
bounds for $(G_{q,1}(f))^*$.

\begin{corollary}\label{cor:rates-direct-rearrangement}
	Under the hypotheses of
	Theorem~\ref{thm:rates-rearranged-decay},
	the following equivalences hold for all integers $k\ge1$.
	\begin{enumerate}
		\item[\textup{(i)}]
		If $p<2<q$ and \eqref{eq:rates-growth} holds with
		$0<\gamma<1/2$, then
		\begin{equation}\label{eq:rates-direct-slow}
			a_k(H_f)\lesssim\Psi(k)^{-1}
			\quad\Longleftrightarrow\quad
			(G_{q,1}(f))^*(k)
			\lesssim\Psi(k^{2/\tau})^{-1}.
		\end{equation}
		\item[\textup{(ii)}]
		If $p<2<q$ and $\Psi(t)/t^\delta$ is nondecreasing
		for some $\delta>1/2$, then
		\begin{equation}\label{eq:rates-direct-fast}
			a_k(H_f)\lesssim\Psi(k)^{-1}
			\quad\Longleftrightarrow\quad
			(G_{q,1}(f))^*(k)
			\lesssim k^{\eta_{p,q}}\Psi(k)^{-1}.
		\end{equation}
		\item[\textup{(iii)}]
		If $p>q$ and $\Psi(t)/t^\delta$ is nondecreasing
		for some $\delta>0$, then
		\begin{equation}\label{eq:rates-direct-tail}
			a_k(H_f)\lesssim\Psi(k)^{-1}
			\quad\Longleftrightarrow\quad
			(G_{q,1}(f))^*(k)
			\lesssim k^{-1/r}\Psi(k)^{-1}.		\end{equation}
	\end{enumerate}
The implicit constants are independent of $k$.
\end{corollary}

\begin{proof}
	Put $F:=G_{q,1}(f)$. For \textup{(i)} and \textup{(ii)}, define
	\[
	I(t):=\int_t^{t^{\tau/2}}
	u^{\varsigma/\tau-1}F^*(u)^\varsigma\,du.
	\]
	By monotonicity and \eqref{eq:rates-growth},
	it suffices to prove the pointwise bounds for
	sufficiently large real $t$.
	The remaining indices are absorbed into the constants,
	since $F\in L^\infty$ and $H_f$ is bounded.
	
	\textup{(i)}
	The proof of Theorem~\ref{thm:rates-rearranged-decay}\textup{(ii)}
	gives the real-variable form of \eqref{eq:rates-rearranged-crossing},
	\begin{equation}\label{eq:direct-window-decay-criterion}
		a_k(H_f)\lesssim\Psi(k)^{-1}
		\quad\Longleftrightarrow\quad
		I(t)\lesssim t^{\varsigma/2}\Psi(t)^{-\varsigma},
	\end{equation}
	where the estimates hold for sufficiently large $k$ and $t$,
	respectively.
	
	Since $\tau>2$, we have $t^{\tau/2}/2\ge t$ for sufficiently
	large $t$. By monotonicity of $F^*$,
	\[
	I(t)\ge F^*(t^{\tau/2})^\varsigma
	\int_{t^{\tau/2}/2}^{t^{\tau/2}}u^{\varsigma/\tau-1}\,du
	\simeq t^{\varsigma/2}F^*(t^{\tau/2})^\varsigma.
	\]
	Thus \eqref{eq:direct-window-decay-criterion} gives
	$F^*(t^{\tau/2})\lesssim\Psi(t)^{-1}$. Taking $t=k^{2/\tau}$
	proves the pointwise estimate.
	
	Conversely, assume $F^*(k)\lesssim\Psi(k^{2/\tau})^{-1}$.
	For sufficiently large $u$, monotonicity of $F^*$ and
	\eqref{eq:rates-growth} give
	\[
	F^*(u)\le F^*(\lfloor u\rfloor)
	\lesssim\Psi(\lfloor u\rfloor^{2/\tau})^{-1}
	\lesssim\Psi(u^{2/\tau})^{-1}.
	\]
	Substituting $v=u^{2/\tau}$ in $I(t)$, we obtain
	\[
	I(t)\lesssim\int_{t^{2/\tau}}^t
	v^{\varsigma/2-1}\Psi(v)^{-\varsigma}\,dv.
	\]
	On this interval, \eqref{eq:rates-growth} with the exponent
	$0<\gamma<1/2$ gives $\Psi(t)\lesssim(t/v)^\gamma\Psi(v)$.
	Consequently,
	\[
	I(t)\lesssim
	\frac{t^{\varsigma\gamma}}{\Psi(t)^\varsigma}
	\int_{t^{2/\tau}}^t v^{\varsigma(1/2-\gamma)-1}\,dv
	\lesssim t^{\varsigma/2}\Psi(t)^{-\varsigma},
	\]
	where the last estimate uses $\varsigma(1/2-\gamma)>0$.
	The reverse implication follows from
	\eqref{eq:direct-window-decay-criterion}, completing the proof of
	\textup{(i)}. 
	
\textup{(ii)}
Since $\tau>2$, the interval $[t,2t]$ is contained in
$[t,t^{\tau/2}]$ for sufficiently large $t$. By monotonicity of $F^*$,
\[
I(t)\ge F^*(2t)^\varsigma
\int_t^{2t}u^{\varsigma/\tau-1}\,du
\simeq t^{\varsigma/\tau}F^*(2t)^\varsigma.
\]
Thus \eqref{eq:direct-window-decay-criterion} and
$\eta_{p,q}=1/2-1/\tau$ give
$F^*(2t)\lesssim t^{\eta_{p,q}}\Psi(t)^{-1}$.
Taking $t=k/2$ and using $\Psi(k/2)\simeq\Psi(k)$ by
\eqref{eq:rates-growth} proves necessity.

Conversely, assume $F^*(k)\lesssim k^{\eta_{p,q}}\Psi(k)^{-1}$.
As in \textup{(i)}, monotonicity and \eqref{eq:rates-growth}
extend this estimate to real $u\ge1$. Since $\Psi(u)/u^\delta$
is nondecreasing, $\Psi(u)\ge(u/t)^\delta\Psi(t)$ for $u\ge t$.
Therefore
\[
I(t)\lesssim
\frac{t^{\varsigma\delta}}{\Psi(t)^\varsigma}
\int_t^\infty u^{\varsigma(1/2-\delta)-1}\,du
\lesssim t^{\varsigma/2}\Psi(t)^{-\varsigma},
\]
where the last step uses $\delta>1/2$. The converse follows from
\eqref{eq:direct-window-decay-criterion}, completing the proof of
\textup{(ii)}.

\textup{(iii)}
By Theorem~\ref{thm:rates-rearranged-decay}\textup{(iii)},
the approximation number bound is equivalent to
\eqref{eq:rates-rearranged-tail}. For $k\ge2$, set
$m:=\lfloor k/2\rfloor$. By monotonicity of $F^*$ and
\eqref{eq:rates-rearranged-tail},
\[
(k-m)F^*(k)^r
\le \int_m^k F^*(u)^r\,du
\le \int_m^\infty F^*(u)^r\,du
\lesssim \Psi(m)^{-r}.
\]
Since $k-m\ge k/2$ and $\Psi(m)\simeq\Psi(k)$ by
\eqref{eq:rates-growth}, this proves necessity.

Conversely, assume $F^*(k)\lesssim k^{-1/r}\Psi(k)^{-1}$.
Monotonicity and \eqref{eq:rates-growth} again extend the estimate
to real $u\ge1$. Using $\Psi(u)\ge(u/k)^\delta\Psi(k)$ for $u\ge k$,
we obtain
\[
\int_k^\infty F^*(u)^r\,du
\lesssim \frac{k^{r\delta}}{\Psi(k)^r}
\int_k^\infty u^{-r\delta-1}\,du
=\frac{1}{r\delta\,\Psi(k)^r}.
\]
Taking the $r$th root and applying \eqref{eq:rates-rearranged-tail},
we obtain $a_k(H_f)\lesssim\Psi(k)^{-1}$.
This completes the proof.
\end{proof}

\subsection{A converse to the weighted Bergman decay criterion}\label{sec:Bergman-converse}
We next give an application of the local
recovery method to weighted Bergman spaces. In the Hilbert
space setting, a finite-dimensional trace estimate gives
a singular value lower bound that supplies the converse
to \cite[Theorem 1.4(2)]{FWZ26}.

 We use the notation and hypotheses of \cite{FWZ26}. Let $\Omega\subset\C$ be a domain and let $\mathcal W^*(\Omega)$
denote the class of weights defined in \cite[Definition~2.6]{FWZ26}.
For $\omega\in\mathcal W^*(\Omega)$, set $dA_\omega=\omega\,dA$,
where $dA$ denotes the Lebesgue measure. The corresponding
weighted Bergman space is
\[
A_\omega^2
:=\left\{
f\in\mathcal{H}(\Omega):
\|f\|_{A_\omega^2}
:=\left(\int_\Omega |f(z)|^2\,dA_\omega(z)\right)^{1/2}
<\infty
\right\}.
\] This family includes standard
Bergman spaces on the unit disc, standard Fock spaces, and weighted
Fock spaces. The reproducing kernel $K(\cdot,\cdot)$ of $A_\omega^2$
satisfies the hypotheses $(\mathrm{H})_1$--$(\mathrm{H})_4$ in
\cite[Definition~2.2]{FWZ26}.
In particular, set
\[
\tau_\omega(z):=
\frac1{\omega(z)^{1/2}K(z,z)^{1/2}},
\qquad
d\lambda_\omega(z):=\frac{dA(z)}{\tau_\omega(z)^2}.
\]

To ensure that the Hankel operator $H_f$ is densely defined from
$A_\omega^2$ into $L_\omega^2:=L^2(\Omega,dA_\omega)$, the authors
of \cite[Definition~2.7]{FWZ26} assume that $f$ belongs to the
symbol class
\begin{align}\label{BergS}
\mathcal S_\omega
:=\bigl\{f\text{ measurable on }\Omega:
          fg\in L_\omega^2\text{ for every }g\in\Gamma\bigr\},
\end{align}
where $\Gamma=\operatorname{span}\{K(\cdot,z):z\in\Omega\}$.

We next recall the IDA function $G_{\delta,\omega}(f)$ for
weights $\omega\in\mathcal W^*(\Omega)$, as introduced in
\cite[Definition~1.2]{FWZ26}. Let $B_\omega>1$ and $\delta_\omega>0$ be the constants
in the covering lemma [15, Lemma~2.3].
For $f\in L_{\mathrm{loc}}^2(\Omega)$ and
$0<\delta<B_\omega\delta_\omega$, define $G_{\delta,\omega}(f)$ by
\begin{align}\label{BergIDA}
G_{\delta,\omega}(f)(z)
=\inf_{h\in\mathcal{H}(D_{\delta,\omega}(z))}
\left(
\frac{1}{A(D_{\delta,\omega}(z))}
\int_{D_{\delta,\omega}(z)} |f(w)-h(w)|^2\,dA(w)
\right)^{1/2},
\qquad z\in\Omega,
\end{align}
where $A(E)=\int_E dA$ and $D_{\delta,\omega}(z)
=\{w\in\Omega:|w-z|<\delta\tau_\omega(z)\}.$

\begin{proposition}\label{prop:bergman-diagonal-recovery}
Let $\omega\in\mathcal W^*(\Omega)$, $0<\delta<\delta_\omega$,
and $f\in\mathcal S_\omega$. Suppose that
$H_f:A_\omega^2\to L_\omega^2$ is compact. Then there exist
constants $A\ge1$ and $C>0$, independent of $f$ and $k$, such that
\begin{equation}\label{eq:bergman-rearrangement-recovery}
(G_{\delta,\omega}(f))^*(Ak)\le C s_k(H_f),
\qquad k\ge1,
\end{equation}
where the rearrangement is taken with respect to $d\lambda_\omega$.
\end{proposition}

\begin{proof}
	We adapt the local compression argument for Fock spaces in the proof of
	Proposition~\ref{thm:rates-diagonal-transfer}, specifically
	\eqref{eq:local-distance} and \eqref{eq:UV}--\eqref{eq:matrix}.
	
	Let $B_\omega>1$, $\delta_\omega>0$, and the lattice
	$\{a_j\}_{j\in\mathbb N}\subset\Omega$ be as in
	\cite[Lemma~2.3]{FWZ26}. For $0<\delta<\delta_\omega$, set
	\begin{equation}\label{eq:bergman-covering-discs}
		D_j:=D(a_j,\delta\tau_\omega(a_j)),\qquad
		\widehat D_j:=D(a_j,B_\omega\delta\tau_\omega(a_j)),\qquad
		\widetilde D_j:=D(a_j,2B_\omega\delta\tau_\omega(a_j)).
	\end{equation}
	Take $\delta_\omega$ small enough that $\widetilde D_j\subset\Omega$
	and $(\mathrm H)_4$ applies on $\widetilde D_j$.
By \cite[Lemma~2.3(2), (4) and Definition~2.2, $(\mathrm H)_4$]{FWZ26},
there is $\varepsilon_\omega>0$ such that
$D(a_j,\varepsilon_\omega\delta\tau_\omega(a_j))$ are disjoint and
\[
z\in\widetilde D_j\ \Longrightarrow\
\tau_\omega(a_j)\simeq\tau_\omega(z),\qquad
D(a_j,\varepsilon_\omega\delta\tau_\omega(a_j))
\subset D(z,C\delta\tau_\omega(z)).
\]
Area comparison gives
	\begin{equation}\label{eq:bergman-overlap}
		\sum_{j\in\mathbb N}\chi_{\widetilde D_j}(z)\lesssim1,
		\qquad z\in\Omega.
	\end{equation}
	Let $k_{a_j}$ be the normalized reproducing kernel of $A_\omega^2$
	at $a_j$. By \cite[Definition~2.2, $(\mathrm H)_4$]{FWZ26},
	\begin{equation}\label{eq:bergman-local-kernel}
		|k_{a_j}(z)|^2\omega(z)\simeq\tau_\omega(a_j)^{-2},
		\qquad z\in\widetilde D_j.
	\end{equation}
 Set
	\begin{equation}\label{eq:bergman-local-distance}
		d_j:=\inf_{h\in\mathcal H(\widehat D_j)}
		\|(f-h)k_{a_j}\|_{L^2(\widehat D_j,dA_\omega)}.
	\end{equation}
	By \cite[proof of Theorem~4.2, Step~2]{FWZ26},
	\begin{equation}\label{eq:bergman-local-ida}
		d_j\simeq G_{B_\omega\delta,\omega}(f)(a_j).
	\end{equation}
	Let $P_j$ be the orthogonal projection onto the holomorphic subspace
	of $L^2(\widehat D_j,dA_\omega)$. Since $k_{a_j}$ has no zeros on
	$\widehat D_j$ by \eqref{eq:bergman-local-kernel},
	\begin{equation}\label{eq:bergman-local-projection}
		d_j=\|(I-P_j)(fk_{a_j})\|_{L^2(\widehat D_j,dA_\omega)}.
	\end{equation}
	For $d_j>0$, put $v_j:=d_j^{-1}(I-P_j)(fk_{a_j})$. Then
	\begin{equation}\label{eq:bergman-local-test}
		\|v_j\|_{L^2(\widehat D_j,dA_\omega)}=1,\qquad
		v_j\perp\operatorname{Ran}P_j,\qquad
		\langle fk_{a_j},v_j\rangle_{L^2(\widehat D_j,dA_\omega)}=d_j.
	\end{equation}
	
	Fix $t>0$ and a finite set $J\subset\{j\in\mathbb N:d_j>t\}$.
	For $\mathbf c\in\ell^2_J$ and $v\in L_\omega^2$, define
	\begin{equation}\label{eq:bergman-compression-maps}
		U_J\mathbf c:=\sum_{j\in J}c_jk_{a_j},\qquad
		(V_Jv)_j:=d_j^{-1}
		\langle v,v_j\rangle_{L^2(\widehat D_j,dA_\omega)}.
	\end{equation}
	For $u\in A_\omega^2$, applying the submean inequality
	to $u/k_{a_j}$ and using \eqref{eq:bergman-local-kernel},
	we obtain
	\begin{equation}\label{eq:bergman-local-submean}
		\begin{aligned}
			\sup_{\widehat D_j}|u/k_{a_j}|^2
			\lesssim
			\frac1{\delta^2\tau_\omega(a_j)^2}
			\int_{\widetilde D_j}|u/k_{a_j}|^2\,dA
			\lesssim
			\delta^{-2}\int_{\widetilde D_j}|u|^2\,dA_\omega.
		\end{aligned}
	\end{equation}
	Let $e_j$ be the $j$th coordinate vector in $\ell^2_J$.
	By \eqref{eq:bergman-compression-maps}, $U_Je_j=k_{a_j}$.
	The adjoint identity gives
	\begin{equation}\label{eq:bergman-compression-adjoint}
		(U_J^*u)_j
		=\langle U_J^*u,e_j\rangle_{\ell^2_J}
		=\langle u,U_Je_j\rangle_{A_\omega^2}
		=\langle u,k_{a_j}\rangle_{A_\omega^2}.
	\end{equation}
	By the reproducing property,
	\begin{equation}\label{eq:bergman-kernel-evaluation}
		\langle u,k_{a_j}\rangle_{A_\omega^2}
		=\frac{u(a_j)}{\sqrt{K(a_j,a_j)}}
		=(u/k_{a_j})(a_j),
	\end{equation}
	where $k_{a_j}(a_j)=\sqrt{K(a_j,a_j)}$ by normalization.
	Combining \eqref{eq:bergman-compression-adjoint},
	\eqref{eq:bergman-kernel-evaluation},
	\eqref{eq:bergman-local-submean}, and \eqref{eq:bergman-overlap},
	we obtain
	\begin{equation}\label{eq:bergman-synthesis-bound}
		\|U_J^*u\|_{\ell^2_J}^2
		\lesssim\sum_{j\in J}\int_{\widetilde D_j}|u|^2\,dA_\omega
		\lesssim\|u\|_{A_\omega^2}^2.
	\end{equation}
	By \eqref{eq:bergman-local-test},
	\eqref{eq:bergman-compression-maps}, and \eqref{eq:bergman-overlap},
	\begin{equation}\label{eq:bergman-analysis-bound}
		\|V_Jv\|_{\ell^2_J}^2
		\le t^{-2}\sum_{j\in J}\int_{\widehat D_j}|v|^2\,dA_\omega
		\lesssim t^{-2}\|v\|_{L_\omega^2}^2.
	\end{equation}
	Consequently,
	\begin{equation}\label{eq:bergman-compression-bounds}
		\|U_J\|\lesssim1,\qquad \|V_J\|\lesssim t^{-1}.
	\end{equation}
	
	Set $M_J:=V_JH_fU_J$ and let $u:=U_J\mathbf c$.
	Since $f\in\mathcal S_\omega$, we have $fu\in L_\omega^2$.
	The restriction of $P_\omega(fu)$ to $\widehat D_j$ belongs to
	$\operatorname{Ran}P_j$, so \eqref{eq:bergman-local-test} gives
	\begin{equation}\label{eq:bergman-compression-coordinate}
		(V_JH_fu)_j
		=d_j^{-1}\langle fu,v_j\rangle_{L^2(\widehat D_j,dA_\omega)}.
	\end{equation}
	Taking $u=k_{a_j}=U_Je_j$ and applying
	\eqref{eq:bergman-local-test}, we obtain
	\begin{equation}\label{eq:bergman-compression-diagonal}
		(M_Je_j)_j=1,\qquad j\in J.
	\end{equation}
	
	For $h\in\mathcal H(\widehat D_j)$ with
	$(f-h)k_{a_j}\in L^2(\widehat D_j,dA_\omega)$, we have
	$hk_{a_j}\in L^2(\widehat D_j,dA_\omega)$. By
	\eqref{eq:bergman-local-submean},
	$hu=(hk_{a_j})(u/k_{a_j})\in\operatorname{Ran}P_j$.
	Thus \eqref{eq:bergman-compression-coordinate} and
	\eqref{eq:bergman-local-test} give
	\begin{equation}\label{eq:bergman-compression-local-bound}
		|(V_JH_fu)_j|
		\le d_j^{-1}\|(f-h)k_{a_j}\|_{L^2(\widehat D_j,dA_\omega)}
		\sup_{\widehat D_j}|u/k_{a_j}|.
	\end{equation}
	Taking the infimum over $h$ in
	\eqref{eq:bergman-compression-local-bound} and using
	\eqref{eq:bergman-local-distance},
	\eqref{eq:bergman-local-submean}, and \eqref{eq:bergman-overlap},
	we obtain
	\begin{equation}\label{eq:bergman-normalized-analysis}
		\|V_JH_fu\|_{\ell^2_J}^2
		\le\sum_{j\in J}\sup_{\widehat D_j}|u/k_{a_j}|^2
		\lesssim\sum_{j\in J}\int_{\widetilde D_j}|u|^2\,dA_\omega
		\lesssim\|u\|_{A_\omega^2}^2.
	\end{equation}
	Substituting $u=U_J\mathbf c$ and applying
	\eqref{eq:bergman-compression-bounds}, we conclude that
	\begin{equation}\label{eq:bergman-normalized-compression}
		\|M_J\|\le C_*
	\end{equation}
	for a fixed $C_*\ge1$.
	By \eqref{eq:bergman-compression-diagonal} and
	\eqref{eq:bergman-normalized-compression},
	\begin{equation}\label{eq:bergman-compression-singular-count}
		\#J\le\sum_{\ell=1}^{\#J}s_\ell(M_J)^2
		\le C_*^2\#\{1\le\ell\le\#J:s_\ell(M_J)\ge\tfrac12\}
		+\tfrac14\#J.
	\end{equation}
	Set $c:=3/(4C_*^2)$. By
	\eqref{eq:bergman-compression-singular-count} and the decreasing
	order of the singular values,
	\begin{equation}\label{eq:bergman-compression-singular-lower}
		s_k(M_J)\ge\tfrac12,\qquad 1\le k\le c\,\#J.
	\end{equation}
	By Lemma~\ref{lem:ideal-properties}\textup{(i)} and
	\eqref{eq:bergman-compression-bounds},
	\begin{equation}\label{eq:bergman-compression-singular-upper}
		s_k(M_J)
		\le\|V_J\|\|U_J\|s_k(H_f)
		\le C_1t^{-1}s_k(H_f).
	\end{equation}
	If $\#\{j:d_j>t\}\ge\lceil k/c\rceil$, choose
	$J\subset\{j:d_j>t\}$ with $\#J=\lceil k/c\rceil$.
	Then $k\le c\,\#J$, so
	\eqref{eq:bergman-compression-singular-lower} and
	\eqref{eq:bergman-compression-singular-upper} give
	\begin{equation}\label{eq:bergman-compression-threshold}
		\tfrac12\le s_k(M_J)\le C_1t^{-1}s_k(H_f),
		\qquad t\le2C_1s_k(H_f).
	\end{equation}
	Consequently,
	\begin{equation}\label{eq:bergman-local-count}
		\#\{j:d_j>t\}<\lceil k/c\rceil\lesssim k,
		\qquad t>2C_1s_k(H_f).
	\end{equation}
	
	By \cite[proof of Theorem~4.2, Step~2]{FWZ26} and
	\eqref{eq:bergman-local-ida}, there is a fixed $C'>0$ such that
	\begin{equation}\label{eq:bergman-ida-covering}
		G_{\delta,\omega}(f)(z)\le C'd_j,\qquad z\in D_j.
	\end{equation}
	Since the discs $D_j$ cover $\Omega$,
	\eqref{eq:bergman-ida-covering} implies
	\begin{equation}\label{eq:bergman-superlevel-covering}
		\{G_{\delta,\omega}(f)>C't\}
		\subseteq\bigcup_{j:\,d_j>t}D_j.
	\end{equation}
	By $(\mathrm H)_4$ and the definition of $d\lambda_\omega$,
	$\lambda_\omega(D_j)\lesssim1$. Hence
	\eqref{eq:bergman-superlevel-covering} and
	\eqref{eq:bergman-local-count} give
	\begin{equation}\label{eq:bergman-distribution-bound}
		\lambda_\omega\bigl(\{G_{\delta,\omega}(f)>C't\}\bigr)
		\le\sum_{j:\,d_j>t}\lambda_\omega(D_j)
		\lesssim\#\{j:d_j>t\}
		\le Ak,\qquad t>2C_1s_k(H_f),
	\end{equation}
	for a fixed $A\ge1$. By \eqref{eq:bergman-distribution-bound}
	and the definition of decreasing rearrangement,
	\begin{equation}\label{eq:bergman-rearrangement-threshold}
		(G_{\delta,\omega}(f))^*(Ak)\le C't,
		\qquad t>2C_1s_k(H_f).
	\end{equation}
	Taking the infimum over $t>2C_1s_k(H_f)$ in
	\eqref{eq:bergman-rearrangement-threshold} proves
	\eqref{eq:bergman-rearrangement-recovery}.
	This completes the proof.
\end{proof}

\begin{theorem}\label{cor:bergman-prescribed-decay}
Under the hypotheses of Proposition~\ref{prop:bergman-diagonal-recovery},
let $\Psi:[1,\infty)\to(0,\infty)$ be nondecreasing and suppose
that $\Psi(t)/t^\gamma$ is nonincreasing for some $\gamma>0$.
Then
\begin{equation}\label{eq:bergman-prescribed-decay}
s_k(H_f)\lesssim\Psi(k)^{-1}\quad(k\ge1)
\quad\Longleftrightarrow\quad
(G_{\delta,\omega}(f))^*(k)\lesssim\Psi(k)^{-1}\quad(k\ge1).
\end{equation}
The implicit constants are independent of $k$.
\end{theorem}

\begin{proof}
	Suppose $s_k(H_f)\lesssim\Psi(k)^{-1}$. For $m\ge2A$, set
	$k:=\lfloor m/A\rfloor$. Then $Ak\le m$ and $m/(2A)\le k$.
	Applying \eqref{eq:rates-growth} with $u=k$ and $v=m$ gives
	\begin{equation}\label{eq:bergman-weight-dilation}
		\Psi(m)\le C_0(m/k)^\gamma\Psi(k)
		\le C_0(2A)^\gamma\Psi(k).
	\end{equation}
	By monotonicity of the rearrangement,
	\eqref{eq:bergman-rearrangement-recovery}, and
	\eqref{eq:bergman-weight-dilation},
	\begin{equation}\label{eq:bergman-decay-necessity}
		(G_{\delta,\omega}(f))^*(m)
		\le (G_{\delta,\omega}(f))^*(Ak)
		\lesssim s_k(H_f)
		\lesssim\Psi(k)^{-1}
		\lesssim\Psi(m)^{-1}.
	\end{equation}
	For the finitely many integers $m<2A$,
	\cite[Theorem~3.1]{FWZ26} gives
	\begin{equation}\label{eq:bergman-decay-small-indices}
		\Psi(m)(G_{\delta,\omega}(f))^*(m)
		\le\Psi(\lceil2A\rceil)\|G_{\delta,\omega}(f)\|_\infty
		\lesssim\Psi(\lceil2A\rceil)\|H_f\|<\infty.
	\end{equation}
	Together, \eqref{eq:bergman-decay-necessity} and
	\eqref{eq:bergman-decay-small-indices} prove necessity.
	
	Conversely, suppose
	$(G_{\delta,\omega}(f))^*(m)\lesssim\Psi(m)^{-1}$ for all
	$m\in\mathbb N$. Since $H_f$ is compact,
	\cite[Theorem~3.2]{FWZ26} gives $f\in\mathrm{VDA}_\omega(\Omega)$.
	Thus \cite[Theorem~1.4(2)]{FWZ26} applies and yields
	$s_m(H_f)\lesssim\Psi(m)^{-1}$, proving the reverse implication
	in \eqref{eq:bergman-prescribed-decay}.
	This completes the proof.
\end{proof}

	\subsection{Approximation ideals}
	
	\begin{proposition}\label{prop:diagonal-equivalence}
		\label{prop:upper-diagonal}\label{prop:recovery-diagonal}
		Under the hypotheses of
		Proposition~\ref{thm:rates-diagonal-transfer},
		for every $s>0$,
		\[
		H_f\in\mathfrak A_s(F^p(\varphi),L^q(\varphi))
		\quad\Longleftrightarrow\quad
		D_{\mathbf g}\in\mathfrak A_s(\ell^p,\ell^q).
		\]
		Moreover, $\As{H_f}\simeq\As{D_{\mathbf g}}$,
		with comparison constants independent of $f$.
	\end{proposition}
	
\begin{proof}
	Suppose $D_{\mathbf g}\in\mathfrak A_s$. If $D_{\mathbf g}=0$,
	\eqref{eq:rates-diagonal-equivalence} with $\Psi(t)=t$ gives $H_f=0$.
	Otherwise, define
	\begin{equation}\label{eq:diagonal-majorant-definition}
		b_k:=k^{-2/s}\max_{1\le j\le k}j^{2/s}a_j(D_{\mathbf g}),
		\qquad k\in\mathbb N.
	\end{equation}
	Then $b_k>0$ and $a_k(D_{\mathbf g})\le b_k$. By monotonicity of
	the approximation numbers,
	\begin{equation}\label{eq:diagonal-majorant-growth}
		b_{k+1}
		=\max\{(k/(k+1))^{2/s}b_k,a_{k+1}(D_{\mathbf g})\}\le b_k,
		\qquad k\in\mathbb N.
	\end{equation}
	Moreover, \eqref{eq:diagonal-majorant-definition} gives
	$b_k\ge(m/k)^{2/s}b_m$ for $1\le m\le k$.
	
	For $j\in\mathbb N$, monotonicity also gives
	\begin{equation}\label{eq:diagonal-majorant-local}
		j^2a_j(D_{\mathbf g})^s
		\lesssim\sum_{m=\lceil j/2\rceil}^{j}m\,a_m(D_{\mathbf g})^s.
	\end{equation}
	Using \eqref{eq:diagonal-majorant-definition} and
	\eqref{eq:diagonal-majorant-local}, then applying Tonelli's theorem,
	we obtain
	\begin{equation}\label{eq:diagonal-summable-majorant}
		\sum_{k=1}^{\infty}b_k^s
		\lesssim\sum_{k=1}^{\infty}k^{-2}
		\sum_{m=1}^{k}m\,a_m(D_{\mathbf g})^s
		=\sum_{m=1}^{\infty}m\,a_m(D_{\mathbf g})^s
		\sum_{k=m}^{\infty}k^{-2}
		\lesssim\sum_{m=1}^{\infty}a_m(D_{\mathbf g})^s.
	\end{equation}
	Thus $b_k\downarrow0$. Set $\Psi(t):=b_{\lfloor t\rfloor}^{-1}$
	for $t\ge1$. By \eqref{eq:diagonal-majorant-definition},
	\begin{equation}\label{eq:diagonal-majorant-weight}
		\frac{\Psi(v)}{\Psi(u)}
		=\frac{b_{\lfloor u\rfloor}}{b_{\lfloor v\rfloor}}
		\le\left(\frac{\lfloor v\rfloor}{\lfloor u\rfloor}\right)^{2/s}
		\le2^{2/s}(v/u)^{2/s},
		\qquad 1\le u\le v.
	\end{equation}
	Hence $\Psi$ is nondecreasing and satisfies
	\eqref{eq:rates-growth} with $C_0=2^{2/s}$ and $\gamma=2/s$.
	Since $a_k(D_{\mathbf g})\le b_k$ for every $k\in\mathbb N$,
	we have $\|D_{\mathbf g}\|_{\mathfrak A_\Psi}\le1$.
	Applying the proof of
	Proposition~\ref{thm:rates-diagonal-transfer}
	with these fixed growth constants gives
	\[
	a_k(H_f)\le Cb_k,\qquad k\in\mathbb N,
	\]
	where $C$ is independent of $f$ and $k$.
	Consequently, \eqref{eq:diagonal-summable-majorant} yields
	\begin{equation}\label{eq:summable-upper-transfer}
		\As{H_f}^s
		=\sum_{k=1}^{\infty}a_k(H_f)^s
		\le C^s\sum_{k=1}^{\infty}b_k^s
		\lesssim \As{D_{\mathbf g}}^s.
	\end{equation}
	
	Conversely, suppose $H_f\in\mathfrak A_s$.
	Let $N$ be the integer in \eqref{eq:rates-recovery-numbers}.
	Summing over $N(k-1)<m\le Nk$ for $k\in\mathbb N$,
	monotonicity and \eqref{eq:rates-recovery-numbers} give
	\begin{equation}\label{eq:summable-lower-transfer}
		\As{D_{\mathbf g}}^s
		\le N\sum_{k=1}^{\infty}a_{N(k-1)+1}(D_{\mathbf g})^s
		\lesssim\sum_{k=1}^{\infty}a_k(H_f)^s
		=\As{H_f}^s.
	\end{equation}
	This completes the proof.
\end{proof}

\begin{proof}[Proof of Theorem~\ref{thm:increasing}]
	First we verify the local integrability required by
	Proposition~\ref{prop:diagonal-equivalence}.
	If $H_f\in\mathfrak A_s$, then $H_fk_a\in L^q(\varphi)$
	for every $a\in\C^n$. By \eqref{eq:kernel-lower} and
	\eqref{eq:normalized-kernel}, $k_a$ has no zeros near $a$.
	Since $P(fk_a)$ is holomorphic, the identity
	\[
	f=\frac{H_fk_a+P(fk_a)}{k_a}
	\]
	gives $f\in L^q_{\rm loc}(\C^n)$.
	Conversely, each asserted function-space condition implies
	$G_{q,1}(f)<\infty$ almost everywhere. For every $z\in\C^n$,
	choose $a\in B(z,1/2)$ with $G_{q,1}(f)(a)<\infty$.
	The definition of $G_{q,1}$ gives $f\in L^q(B(z,1/2))$.
	Thus Proposition~\ref{prop:diagonal-equivalence} applies
	in both directions.
	Fix $\Lambda=\{a_j\}_{j\in\mathbb N}\in\mathscr G$ and put $
	\mathbf g:=(g_j)_{j\in\mathbb N}$, where $g_j:=G_{q,R}(f)(a_j)$.
	
	\textup{(i)}
	Suppose $(p,q,s)\notin\mathcal E$ and set $\rho:=\rho_{p,q,s}$.
	For $p\le q$, \eqref{eq:power-sampling} with $\sigma=\rho$ gives
	\begin{equation}\label{eq:main-power-sampling}
		G_{q,1}(f)\in\mathcal L^{\rho,s}
		\quad\Longleftrightarrow\quad
		\mathbf g\in\ell^{\rho,s}.
	\end{equation}
	Since $\ell^{\rho,s}\subset\ell^\infty$,
	Lemma~\ref{lem:diag-facts}\textup{(i)} ensures boundedness of
	$D_{\mathbf g}:\ell^p\to\ell^q$ under either condition in
	\eqref{eq:main-power-sampling}.
	Lemma~\ref{lem:same-side}\textup{(i)} and
	\eqref{eq:main-power-sampling} therefore give
	\begin{equation}\label{eq:main-power-criterion}
		D_{\mathbf g}\in\mathfrak A_s(\ell^p,\ell^q)
		\quad\Longleftrightarrow\quad
		G_{q,1}(f)\in\mathcal L^{\rho,s}.
	\end{equation}
	For $p>q$, set $r:=pq/(p-q)$, so $\rho=sr/(s+r)$.
	Lemma~\ref{lem:continuous-tail} gives
	\eqref{eq:main-power-criterion} in this case.
	
	By Proposition~\ref{prop:diagonal-equivalence},
	Lemma~\ref{lem:same-side}\textup{(i)}, \eqref{eq:power-sampling},
	and Lemma~\ref{lem:continuous-tail}, we obtain
	\begin{equation}\label{eq:main-power-diagonal}
		\As{H_f}\simeq\As{D_{\mathbf g}}
		\simeq\|G_{q,1}(f)\|_{\mathcal L^{\rho,s}}.
	\end{equation}
	Equations \eqref{eq:main-power-criterion} and
	\eqref{eq:main-power-diagonal} prove
	\eqref{eq:ascending-rearrangement-main}.
	
	\textup{(ii)}
	Suppose $(p,q,s)\in\mathcal E$. Then $0<p<2<q<\infty$,
	$s=2$, and $\tau:=\min\{p^\dagger,q\}>2$.
	
	If $p^\dagger=q$, then $\tau=q$. By
	\eqref{eq:logarithmic-sampling} with $\sigma=q$, $s=2$,
	and $\beta=1/2$,
	\begin{equation}\label{eq:main-critical-log-sampling}
		G_{q,1}(f)\in\mathcal L^{q,2,1/2}
		\quad\Longleftrightarrow\quad
		\mathbf g\in\ell^{q,2,1/2}\subset\ell^\infty.
	\end{equation}
	Thus Lemma~\ref{lem:diag-facts}\textup{(i)} gives boundedness
	of $D_{\mathbf g}$ whenever the function-space condition holds.
	By \eqref{eq:critical-diagonal-sufficiency},
	\eqref{eq:main-critical-log-sampling}, and
	Proposition~\ref{prop:diagonal-equivalence},
	\begin{equation}\label{eq:main-critical-selfdual-criterion}
		H_f\in\mathfrak A_2
		\quad\Longleftrightarrow\quad
		D_{\mathbf g}\in\mathfrak A_2
		\quad\Longleftrightarrow\quad
		G_{q,1}(f)\in\mathcal L^{q,2,1/2}.
	\end{equation}
	The norm estimates in
	Proposition~\ref{prop:diagonal-equivalence},
	\eqref{eq:critical-diagonal-sufficiency}, and
	\eqref{eq:logarithmic-sampling} give
	\begin{equation}\label{eq:critical-function-lower}
		\|H_f\|_{\mathfrak A_2}
		\simeq\|\mathbf g\|_{\ell^{q,2,1/2}}
		\simeq\|G_{q,1}(f)\|_{\mathcal L^{q,2,1/2}}.
	\end{equation}
	Equations \eqref{eq:main-critical-selfdual-criterion} and
	\eqref{eq:critical-function-lower} prove
	\eqref{eq:critical-selfdual-main}.
	
	If $p^\dagger\ne q$, \eqref{eq:critical-Copson-sampling}
	with $\sigma=\tau$ gives
	\begin{equation}\label{eq:main-critical-tail-sampling}
		G_{q,1}(f)\in
		\mathcal L^{\tau,\tau}\cap\mathcal{CL}^{\tau,2}
		\quad\Longleftrightarrow\quad
		\mathcal T_\tau(\mathbf g)<\infty.
	\end{equation}
	The $k=1$ term in the definition of $\mathcal T_\tau$ gives
	\begin{equation}\label{eq:main-critical-tail-boundedness}
		\|\mathbf g\|_{\ell^\infty}
		\le\|\mathbf g\|_{\ell^\tau}
		\le\mathcal T_\tau(\mathbf g).
	\end{equation}
	Hence \eqref{eq:main-critical-tail-sampling},
	\eqref{eq:main-critical-tail-boundedness}, and
	Lemma~\ref{lem:diag-facts}\textup{(i)} ensure boundedness of
	$D_{\mathbf g}$ under the function-space condition.
	By Proposition~\ref{prop:diagonal-equivalence},
	\eqref{eq:critical-diagonal-sufficiency}, and
	\eqref{eq:main-critical-tail-sampling},
	\begin{equation}\label{eq:main-critical-tail-criterion}
		H_f\in\mathfrak A_2
		\quad\Longleftrightarrow\quad
		D_{\mathbf g}\in\mathfrak A_2
		\quad\Longleftrightarrow\quad
		G_{q,1}(f)\in
		\mathcal L^{\tau,\tau}\cap\mathcal{CL}^{\tau,2}.
	\end{equation}
	The norm estimates in
	Proposition~\ref{prop:diagonal-equivalence},
	\eqref{eq:critical-diagonal-sufficiency}, and
	\eqref{eq:critical-Copson-sampling} yield
	\begin{equation}\label{eq:critical-function-upper}
		\|H_f\|_{\mathfrak A_2}^2
		\simeq\mathcal T_\tau(\mathbf g)^2
		\simeq
		\|G_{q,1}(f)\|_{\mathcal L^{\tau,\tau}}^2
		+\|G_{q,1}(f)\|_{\mathcal{CL}^{\tau,2}}^2.
	\end{equation}
	Taking square roots in \eqref{eq:critical-function-upper} and using
	$(x^2+y^2)^{1/2}\simeq x+y$ for $x,y\ge0$ proves
	\eqref{eq:critical-main-upper}.
	This completes the proof of Theorem~\ref{thm:increasing}.
\end{proof}

\begin{proof}[Proof of Theorem~\ref{cor:weak-approximation-ideals}]
	Take $\Psi(t):=t^{1/s}$.
	If $H_f\in\mathfrak A_{s,\infty}$, then $H_f$ is compact,
	and the argument in the proof of Theorem~\ref{thm:increasing}
	gives $f\in L^q_{\rm loc}(\C^n)$.
	If the corresponding function-space condition holds,
	the same local integrability argument applies.
	In either case, set
	$\mathbf g:=(G_{q,R}(f)(a_j))_{j\in\mathbb N}$.
	
	We first verify compactness under the function-space condition.
	For $p\le q$, \eqref{eq:rates-indexed-sampling} gives
	$\mathbf g\in c_0$. Equations
	\eqref{eq:power-sampling-sandwich} and
	\eqref{eq:radius-translation-comparison} then imply
	$f\in\mathrm{VDA}^q$.
	For $p>q$, put $r:=pq/(p-q)$. Since $\rho=sr/(s+r)<r$,
	\eqref{eq:rates-indexed-sampling} gives
	$\mathbf g\in\ell^{\rho,\infty}\subset\ell^r$.
	Equation~\eqref{eq:power-sampling} therefore yields
	$f\in\mathrm{IDA}^{r,q}$.
	Thus $H_f$ is compact by \cite[Theorem~1.1]{HV23}.
	
	If $(p,q,s)\notin\mathcal E$, we may apply
	Theorem~\ref{thm:rates-rearranged-decay} and
	Corollary~\ref{cor:rates-direct-rearrangement}
	in both directions to obtain the membership criteria.
	Their quantitative estimates, combined with
	\eqref{eq:rates-diagonal-equivalence} and
	Lemma~\ref{lem:indexed-sampling}, give the stated norm comparisons,
	with $a_1(D_{\mathbf g})$ controlling the initial terms.
	
	Suppose $(p,q,s)\in\mathcal E$, so $s=2$.
	If $H_f\in\mathfrak A_{2,\infty}$,
	Theorem~\ref{thm:rates-rearranged-decay}\textup{(ii)}
	gives the integral condition defining $\mathcal X^{\tau,\varsigma}$.
	Moreover, \eqref{eq:prefix-crossing-lower}, applied with
	dimension $2k$ and index $k$, and
	\eqref{eq:rates-diagonal-equivalence} give
	\[
	k^{1/\tau}g_{2k}^*
	\lesssim \sqrt{k}\,a_k(D_{\mathbf g})
	\lesssim \|H_f\|_{\mathfrak A_{2,\infty}}.
	\]
	Together with $g_1^*=a_1(D_{\mathbf g})$,
	monotonicity and \eqref{eq:rates-indexed-sampling} yield
	\[
	\|G_{q,1}(f)\|_{\mathcal L^{\tau,\infty}}
	\lesssim \|H_f\|_{\mathfrak A_{2,\infty}}.
	\]
	Hence $G_{q,1}(f)\in\mathcal X^{\tau,\varsigma}$.
	
Conversely, suppose
$G_{q,1}(f)\in\mathcal X^{\tau,\varsigma}$.
The compactness established above and the integral condition
in $\mathcal X^{\tau,\varsigma}$ imply
$H_f\in\mathfrak A_{2,\infty}$ by
Theorem~\ref{thm:rates-rearranged-decay}\textup{(ii)}.
Moreover, taking $k=1$ and $t=1/2$ in
\eqref{eq:rates-indexed-sampling}, we obtain
$
g_1^*
\lesssim \|G_{q,1}(f)\|_{\mathcal L^{\tau,\infty}}
\le \|G_{q,1}(f)\|_{\mathcal X^{\tau,\varsigma}}.
$
Together with this estimate,
\eqref{eq:rates-window-ideal},
\eqref{eq:rates-diagonal-equivalence}, and the window comparisons
in the proof of Theorem~\ref{thm:rates-rearranged-decay} give
\[
\|H_f\|_{\mathfrak A_{2,\infty}}
\simeq \|G_{q,1}(f)\|_{\mathcal X^{\tau,\varsigma}}.
\]
This completes the proof.
\end{proof}

\subsection{Paired Hankel operators}\label{subsec:paired}
We now apply the preceding results to the pair $H_f$ and $H_{\bar f}$.
Their simultaneous approximation-number decay and ideal membership
can be characterized in terms of the mean oscillation of $f$.
\begin{definition}\label{def:mean-oscillation}
	Let $1\le q<\infty$ and
	$f\in L^q_{\rm loc}(\C^n)$.
	For $r>0$, define
	\begin{equation}\label{eq:paired-MO}
		\operatorname{MO}_{q,r}(f)(z):=
		\left(
		\frac1{v(B(z,r))}
		\int_{B(z,r)}|f-f_{B(z,r)}|^q\,dv
		\right)^{1/q},
	\end{equation}
	where
	\[
	f_{B(z,r)}:=
	\frac1{v(B(z,r))}\int_{B(z,r)}f\,dv.
	\]
	For $f\in L^1_{\rm loc}(\C^n)$, we use the same definition,
	allowing the value $+\infty$.
\end{definition}

Recall that, for $0<\sigma,s<\infty$,
$f\in\mathcal L^{\sigma,\sigma}\cap\mathcal{CL}^{\sigma,s}$
if and only if
\[
\|f\|_{\mathcal L^{\sigma,\sigma}\cap\mathcal{CL}^{\sigma,s}}
:=
\|f\|_{\mathcal L^{\sigma,\sigma}}
+\|f\|_{\mathcal{CL}^{\sigma,s}}<\infty.
\]

\begin{lemma}\label{lem:paired-comparison}
	Let $1\le q<\infty$, $0<\beta<\infty$,
	$0<s\le\infty$, and $f\in L^q_{\rm loc}(\C^n)$.
	Then
	\begin{equation}\label{eq:paired-space-comparison}
		\|\operatorname{MO}_{q,1}(f)\|_{\mathcal L^{\beta,s}}
		\simeq
		\|G_{q,1}(f)\|_{\mathcal L^{\beta,s}}
		+\|G_{q,1}(\bar f)\|_{\mathcal L^{\beta,s}}.
	\end{equation}
	The same assertion holds with $\mathcal L^{\beta,s}$
	replaced by $\mathcal L^{\beta,s,\theta}$ when
	$s<\infty$ and $\theta\in\mathbb R$, by
	$\mathcal L^{\sigma,\sigma}\cap\mathcal{CL}^{\sigma,2}$
	when $\sigma>2$, and by $\mathcal X^{\tau,\varsigma}$
	when $\tau>2$ and $1\le\varsigma<\infty$.
	The comparison constants are independent of $f$.
\end{lemma}

\begin{proof}
	For $u\in\{f,\bar f\}$, constants are admissible
	in \eqref{eq:G}, so
	\[
	G_{q,1}(u)\le\operatorname{MO}_{q,1}(f).
	\]
Fix $z$ and $\varepsilon>0$. As in the proof of
\cite[Lemma~6.1]{HV22}, choose $h,k\in\mathcal H(B(z,2))$ with
\[
M_{q,2}(f-h)(z)+M_{q,2}(\bar f-k)(z)
\le G_{q,2}(f)(z)+G_{q,2}(\bar f)(z)+\varepsilon.
\]
Since $h-\bar k$ is harmonic, interior estimates give
\[
\sup_{B(z,1)}|\partial h|
\lesssim M_{q,2}(h-\bar k)(z)
\le M_{q,2}(f-h)(z)+M_{q,2}(\bar f-k)(z).
\]
Thus $\operatorname{MO}_{q,1}(f)(z)
\le2M_{q,1}(f-h)(z)+C\sup_{B(z,1)}|\partial h|$.
Letting $\varepsilon\downarrow0$ gives
\begin{equation}\label{eq:paired-local-two}
		\operatorname{MO}_{q,1}(f)
		\le C\bigl(G_{q,2}(f)+G_{q,2}(\bar f)\bigr).
	\end{equation}
	By \eqref{eq:radius-translation-comparison} and
	\eqref{eq:rearrangement-sum}, there is a fixed $c\ge1$
	such that
	\begin{equation}\label{eq:paired-rearrangement-comparison}
		\bigl(\operatorname{MO}_{q,1}(f)\bigr)^*(ct)
		\le C\bigl(
		(G_{q,1}(f))^*(t)+(G_{q,1}(\bar f))^*(t)
		\bigr).
	\end{equation}
	The finite-sum and dilation estimates used in
	\eqref{eq:logarithmic-sum},
	\eqref{eq:critical-Copson-sum},
	and \eqref{eq:rates-window-continuous-transfer}
	give the upper bound in
	\eqref{eq:paired-space-comparison}.
	For the weak Lorentz space, this follows directly
	by taking $\sup_{t>0}t^{1/\beta}$ times
	\eqref{eq:paired-rearrangement-comparison}.
	The lower bound follows from
	$G_{q,1}(u)\le\operatorname{MO}_{q,1}(f)$.
	This completes the proof.
\end{proof}

\begin{corollary}\label{cor:rates-prescribed-paired}
	Under the corresponding hypotheses of
	Theorem~\ref{thm:rates-rearranged-decay} or
	Corollary~\ref{cor:rates-direct-rearrangement},
	assume $f,\bar f\in\mathcal S$ and both
	$H_f,H_{\bar f}$ are compact.
	Each equivalence remains valid with
	\[
	a_k(H_f)+a_k(H_{\bar f})\lesssim\Psi(k)^{-1}
	\]
	on the left and with $(G_{q,1}(f))^*$ replaced by
	$(\operatorname{MO}_{q,1}(f))^*$ on the right.
\end{corollary}

\begin{proof}
	Use \eqref{eq:paired-rearrangement-comparison} and
	$G_{q,1}(u)\le\operatorname{MO}_{q,1}(f)$
	for $u=f,\bar f$.
	The pointwise and tail bounds are preserved by finite
	sums and fixed dilations, by \eqref{eq:rates-growth}.
	For the window bounds, apply the covering argument
	in \eqref{eq:rates-window-continuous-transfer}.
	The asserted equivalences therefore follow from
	Theorem~\ref{thm:rates-rearranged-decay} and
	Corollary~\ref{cor:rates-direct-rearrangement}.
	This completes the proof.
\end{proof}

\begin{corollary}\label{thm:paired}
	\label{thm:paired-increasing}\label{thm:paired-tail}
	Let $0<p<\infty$, $1\le q<\infty$, and $0<s<\infty$.
	Suppose $\operatorname{Hess}_{\mathbb R}\varphi\simeq E$
	and $f,\bar f\in\mathcal S$.
	Put
	\[
	Y:=
	\begin{cases}
		\mathcal L^{\rho,s},
		&(p,q,s)\notin\mathcal E,\\
		\mathcal L^{q,2,1/2},
		&(p,q,s)\in\mathcal E,\ p^\dagger=q,\\
		\mathcal L^{\tau,\tau}\cap\mathcal{CL}^{\tau,2},
		&(p,q,s)\in\mathcal E,\ p^\dagger\ne q.
	\end{cases}
	\]
	Then both $H_f$ and $H_{\bar f}$ belong to
	$\mathfrak A_s(\F p,\Lw q)$ if and only if
	$\operatorname{MO}_{q,1}(f)\in Y$.
	Moreover,
	\begin{equation}\label{eq:paired-increasing}
		\|H_f\|_{\mathfrak A_s}
		+\|H_{\bar f}\|_{\mathfrak A_s}
		\simeq\|\operatorname{MO}_{q,1}(f)\|_Y,
	\end{equation}
	with comparison constants independent of $f$.
\end{corollary}

\begin{proof}
	Theorem~\ref{thm:increasing}, applied to $f$ and
	$\bar f$, gives the corresponding membership criteria
	and
	\[
	\|H_f\|_{\mathfrak A_s}
	+\|H_{\bar f}\|_{\mathfrak A_s}
	\simeq
	\|G_{q,1}(f)\|_Y+\|G_{q,1}(\bar f)\|_Y.
	\]
	Equation~\eqref{eq:paired-space-comparison} proves
	the assertion.
	As in the proof of Theorem~\ref{thm:increasing},
	ideal membership ensures local $q$-integrability
	when needed.
	This completes the proof.
\end{proof}

\begin{corollary}\label{cor:rates-weak-paired}
	Let $0<p<\infty$, $1\le q<\infty$, and $0<s<\infty$.
	Suppose $\operatorname{Hess}_{\mathbb R}\varphi\simeq E$
	and $f,\bar f\in\mathcal S\cap L^q_{\rm loc}(\C^n)$.
	Put
	\[
	Y:=
	\begin{cases}
		\mathcal L^{\rho,\infty},
		&(p,q,s)\notin\mathcal E,\\
		\mathcal X^{\tau,\varsigma},
		&(p,q,s)\in\mathcal E.
	\end{cases}
	\]
	Then both $H_f$ and $H_{\bar f}$ belong to
	$\mathfrak A_{s,\infty}(\F p,\Lw q)$ if and only if
	$\operatorname{MO}_{q,1}(f)\in Y$.
	Moreover,
	\begin{equation}\label{eq:rates-paired-norm}
		\|H_f\|_{\mathfrak A_{s,\infty}}
		+\|H_{\bar f}\|_{\mathfrak A_{s,\infty}}
		\simeq\|\operatorname{MO}_{q,1}(f)\|_Y,
	\end{equation}
	with comparison constants independent of $f$.
\end{corollary}

\begin{proof}
	Apply Theorem~\ref{cor:weak-approximation-ideals}
	to $f$ and $\bar f$, and use
	\eqref{eq:paired-space-comparison}.
	This completes the proof.
\end{proof}

\section{The Berger--Coburn phenomenon}\label{sec:BC}
In this section, we prove the Berger--Coburn equivalence for bounded
symbols and establish the sharpness of the exponent range.
We begin with maximal operator and derivative estimates needed
for the proof.
\subsection{Maximal operators and derivative estimates}
All Lebesgue spaces below are over $\mathbb R^{2n}$, and
$dv$ denotes Lebesgue measure.
For $u\in L^1_{\rm loc}(\mathbb R^{2n})$, define
\[
\mathcal Mu(x)
:=\sup_{r>0}\frac1{v(B(x,r))}
\int_{B(x,r)}|u(y)|\,dv(y).
\]
For a Schwartz function $u$, write
\[
\widehat u(\xi)
:=\int_{\mathbb R^{2n}}u(x)e^{-\mathrm i x\cdot\xi}\,dv(x).
\]
For $\ell\in\{1,\ldots,2n\}$, define the Riesz transform by
\[
\widehat{\mathcal R_\ell u}(\xi)
:=\frac{\mathrm i\xi_\ell}{|\xi|}\widehat u(\xi),
\qquad \xi\ne0.
\]

For $1\le a_0<a_1\le\infty$, $0<\theta<1$, and
$0<s<\infty$, the real interpolation quasi-norm is
\[
\|u\|_{(L^{a_0},L^{a_1})_{\theta,s}}
:=
\left(
\int_0^\infty
\bigl[t^{-\theta}K(t,u;L^{a_0},L^{a_1})\bigr]^s
\frac{dt}{t}
\right)^{1/s},
\]
where
\[
K(t,u;L^{a_0},L^{a_1})
:=
\inf_{\substack{u=u_0+u_1\\u_i\in L^{a_i},\ i=0,1}}
\left(\|u_0\|_{L^{a_0}}+t\|u_1\|_{L^{a_1}}\right).
\]
For $s=\infty$, set
\[
\|u\|_{(L^{a_0},L^{a_1})_{\theta,\infty}}
:=
\sup_{t>0}t^{-\theta}K(t,u;L^{a_0},L^{a_1}).
\]
See \cite[Section 3.1]{BL76}.

\begin{lemma}\label{lem:known-Lorentz-calculus}
	Let $1<\sigma<\infty$ and $0<s\le\infty$. Then
	\begin{equation}\label{eq:known-Lorentz-calculus}
		\|\mathcal Mu\|_{\mathcal L^{\sigma,s}}
		+\sum_{\ell=1}^{2n}
		\|\mathcal R_\ell u\|_{\mathcal L^{\sigma,s}}
		\le C\|u\|_{\mathcal L^{\sigma,s}}.
	\end{equation}
	The same assertion holds with $\mathcal L^{\sigma,s}$
	replaced by $\mathcal L^{\sigma,2,1/2}$ or
	$\mathcal L^{\sigma,\sigma}\cap\mathcal{CL}^{\sigma,2}$,
	and also by $\mathcal X^{\sigma,s}$ when
	$\sigma>2$ and $1\le s<\infty$.
	The constants are independent of $u$.
\end{lemma}

\begin{proof}
	Fix $T\in\{\mathcal M,\mathcal R_1,\ldots,\mathcal R_{2n}\}$.
	Choose $1<a_0<\sigma<a_1<\infty$ and $\theta\in(0,1)$ such that
	\[
	\frac1\sigma=\frac{1-\theta}{a_0}+\frac{\theta}{a_1}.
	\]
	The classical bounds \cite[Chapters I and III]{Stein70}
	(see also \cite{G14})
	give
	\begin{equation}\label{eq:calculus-Lebesgue-bounds}
		\|Tu\|_{L^{a_i}}\le C\|u\|_{L^{a_i}},
		\qquad i=0,1.
	\end{equation}
	By \cite[Section 5.3]{BL76},
	\begin{equation}\label{eq:Lorentz-interpolation}
		(L^{a_0},L^{a_1})_{\theta,s}
		=\mathcal L^{\sigma,s},
		\qquad 0<s\le\infty,
	\end{equation}
	with equivalent quasi-norms.
	The linear and sublinear interpolation theorems
	\cite[Section 3.1]{BL76}, together with
	\eqref{eq:calculus-Lebesgue-bounds} and
	\eqref{eq:Lorentz-interpolation}, yield
	\[
	\|Tu\|_{\mathcal L^{\sigma,s}}
	\le C\|u\|_{\mathcal L^{\sigma,s}}.
	\]
	Summing over $T$ proves \eqref{eq:known-Lorentz-calculus}.
	
	For the remaining spaces, let
	\[
	Y\in\left\{
	\mathcal L^{\sigma,2,1/2},\
	\mathcal L^{\sigma,\sigma}\cap\mathcal{CL}^{\sigma,2},\
	\mathcal X^{\sigma,s}
	\right\},
	\]
	with $\sigma>2$ and $1\le s<\infty$ in the last case.
	Recall that $\lambda(t):=1+\log_+t$, and define
	\[
	\|\phi\|_{\widehat Y}:=
	\begin{cases}
		\displaystyle
		\left(
		\int_0^\infty
		t^{2/\sigma-1}\lambda(t)|\phi(t)|^2\,dt
		\right)^{1/2},
		&Y=\mathcal L^{\sigma,2,1/2},
		\\[6pt]
		\displaystyle
		\|\phi\|_{L^\sigma(0,\infty)}
		+\left(
		\int_1^\infty
		\left(\int_t^\infty|\phi(r)|^\sigma\,dr\right)^{2/\sigma}
		\frac{dt}{t}
		\right)^{1/2},
		&Y=\mathcal L^{\sigma,\sigma}\cap\mathcal{CL}^{\sigma,2},
		\\[6pt]
		\displaystyle
		\operatorname*{ess\,sup}_{t>0}t^{1/\sigma}|\phi(t)|
		+\sup_{t\ge2}
		\left(
		\int_t^{t^{\sigma/2}}
		r^{s/\sigma-1}|\phi(r)|^s\,dr
		\right)^{1/s},
		&Y=\mathcal X^{\sigma,s}.
	\end{cases}
	\]
By the definitions of these spaces,
\begin{equation}\label{eq:calculus-profile-norm}
	\|u\|_Y\simeq\|u^*\|_{\widehat Y}.
\end{equation}
The defining integrals and Minkowski's inequality yield
\[
|\phi|\le|\psi|
\ \Longrightarrow\
\|\phi\|_{\widehat Y}\le\|\psi\|_{\widehat Y},
\qquad
\|\phi+\psi\|_{\widehat Y}
\le\|\phi\|_{\widehat Y}+\|\psi\|_{\widehat Y}.
\]
	For $c>0$, put $D_c\phi(t):=\phi(ct)$. For the first two spaces, we claim that
	\begin{equation}\label{eq:critical-calculus-dilation}
		\|D_c\phi\|_{\widehat Y}
		\le Cc^{-1/\sigma}
		\bigl(1+\log_+(1/c)\bigr)^{1/2}
		\|\phi\|_{\widehat Y}.
	\end{equation}

	For $Y=\mathcal L^{\sigma,2,1/2}$, substitution gives
	\[
	\|D_c\phi\|_{\widehat Y}^2
	=c^{-2/\sigma}
	\int_0^\infty
	t^{2/\sigma-1}\lambda(t/c)|\phi(t)|^2\,dt.
	\]
Since $\lambda(t/c)\le
\bigl(1+\log_+(1/c)\bigr)\lambda(t)$,
\eqref{eq:critical-calculus-dilation} follows.
	
	For $Y=\mathcal L^{\sigma,\sigma}\cap\mathcal{CL}^{\sigma,2}$,
	the $L^\sigma$ term scales by $c^{-1/\sigma}$, while
	\[
	\begin{aligned}
		&\int_1^\infty
		\left(\int_t^\infty|\phi(cr)|^\sigma\,dr\right)^{2/\sigma}
		\frac{dt}{t}\\
		&\quad=c^{-2/\sigma}
		\int_c^\infty
		\left(\int_t^\infty|\phi(r)|^\sigma\,dr\right)^{2/\sigma}
		\frac{dt}{t}\\
		&\quad\le c^{-2/\sigma}
		\left[
		\int_1^\infty
		\left(\int_t^\infty|\phi(r)|^\sigma\,dr\right)^{2/\sigma}
		\frac{dt}{t}
		+\log_+(1/c)\|\phi\|_{L^\sigma(0,\infty)}^2
		\right].
	\end{aligned}
	\]
	Thus \eqref{eq:critical-calculus-dilation} holds for this space as well.
	
	For $Y=\mathcal X^{\sigma,s}$, put
	$A:=\operatorname*{ess\,sup}_{r>0}r^{1/\sigma}|\phi(r)|.$
	For $t\ge2$,
	\begin{equation}\label{eq:calculus-window-substitution}
		\int_t^{t^{\sigma/2}}
		r^{s/\sigma-1}|\phi(cr)|^s\,dr
		=
		c^{-s/\sigma}
		\int_{ct}^{ct^{\sigma/2}}
		r^{s/\sigma-1}|\phi(r)|^s\,dr.
	\end{equation}
	If $c\ge1$, then
	$[ct,ct^{\sigma/2}]\subseteq[ct,(ct)^{\sigma/2}]$.
	If $0<c<1$ and $ct\ge2$, splitting at $(ct)^{\sigma/2}$ gives
	\[
	\int_{ct}^{ct^{\sigma/2}}
	r^{s/\sigma-1}|\phi(r)|^s\,dr
	\le
	\|\phi\|_{\widehat Y}^s
	+A^s(\sigma/2-1)\log(1/c).
	\]
	If $ct<2$, then
	\[
	\int_{ct}^{ct^{\sigma/2}}
	r^{s/\sigma-1}|\phi(r)|^s\,dr
	\le A^s(\sigma/2-1)\log t
	\le CA^s\bigl(1+\log(1/c)\bigr).
	\]
	Together with \eqref{eq:calculus-window-substitution}, these bounds yield
	\begin{equation}\label{eq:rates-window-dilation}
		\|D_c\phi\|_{\widehat Y}
		\le Cc^{-1/\sigma}
		\bigl(1+\log_+(1/c)\bigr)^{1/s}
		\|\phi\|_{\widehat Y}.
	\end{equation}
	
	We now treat all three spaces simultaneously.
	For $Y=\mathcal L^{\sigma,2,1/2}$, monotonicity gives
	\[
	t^{2/\sigma}u^*(t)^2
	\lesssim
	\int_{t/2}^t r^{2/\sigma-1}u^*(r)^2\,dr
	\le\|u\|_Y^2.
	\]
	For $Y=\mathcal L^{\sigma,\sigma}\cap\mathcal{CL}^{\sigma,2}$,
	\[
	tu^*(t)^\sigma
	\le\int_0^t u^*(r)^\sigma\,dr
	\le\|u\|_Y^\sigma.
	\]
	For $\mathcal X^{\sigma,s}$, the weak Lorentz bound is part
	of its definition. Hence
	\begin{equation}\label{eq:critical-calculus-embedding}
		\|u\|_{\mathcal L^{\sigma,\infty}}\lesssim\|u\|_Y.
	\end{equation}
	
For $t>0$, choose a measurable set $E_t$ with
\[
v(E_t)=t,\qquad
\{|u|>u^*(t)\}\subseteq E_t\subseteq\{|u|\ge u^*(t)\},
\]
and set $u_0:=u\chi_{E_t}$ and $u_1:=u-u_0$.
By \eqref{eq:critical-calculus-embedding} and $a_0<\sigma<a_1$,
\begin{equation}\label{eq:calculus-level-decomposition}
	\|u_0\|_{L^{a_0}}
	=\left(\int_0^t u^*(r)^{a_0}\,dr\right)^{1/a_0}<\infty,
	\qquad
	\|u_1\|_{L^{a_1}}
	=\left(\int_t^\infty u^*(r)^{a_1}\,dr\right)^{1/a_1}<\infty.
\end{equation}
Thus $u\in L^{a_0}+L^{a_1}$, and we define
\[
\mathcal R_\ell u
:=\mathcal R_\ell u_0+\mathcal R_\ell u_1,
\qquad \ell\in\{1,\ldots,2n\},
\]
using the $L^{a_0}$ and $L^{a_1}$ extensions, respectively,
which agree on $L^{a_0}\cap L^{a_1}$.
By sublinearity and \eqref{eq:rearrangement-sum},
\[
(Tu)^*(2t)\le (Tu_0)^*(t)+(Tu_1)^*(t).
\]
Chebyshev's inequality and \eqref{eq:calculus-Lebesgue-bounds} give
\[
(Tu_i)^*(t)
\le t^{-1/a_i}\|Tu_i\|_{L^{a_i}}
\lesssim t^{-1/a_i}\|u_i\|_{L^{a_i}},
\qquad i=0,1.
\]
Substituting \eqref{eq:calculus-level-decomposition}, we obtain
\begin{equation}\label{eq:calculus-endpoint-rearrangement}
	(Tu)^*(2t)
	\lesssim
	t^{-1/a_0}
	\left(\int_0^t u^*(r)^{a_0}\,dr\right)^{1/a_0}
	+t^{-1/a_1}
	\left(\int_t^\infty u^*(r)^{a_1}\,dr\right)^{1/a_1}.
\end{equation}

Define
\[
P_{a_0}u^*(t):=\int_0^1 r^{1/a_0-1}u^*(tr)\,dr,
\qquad
Q_{a_1}u^*(t):=\int_1^\infty r^{1/a_1-1}u^*(tr)\,dr.
\]
By monotonicity of $u^*$,
\[
\begin{aligned}
	t^{-1/a_0}
	\left(\int_0^t u^*(r)^{a_0}\,dr\right)^{1/a_0}
	=
	\left(
	\frac1t\sum_{j\in\mathbb N}
	\int_{2^{-j}t}^{2^{1-j}t}u^*(r)^{a_0}\,dr
	\right)^{1/a_0}
	\le
	\sum_{j\in\mathbb N}2^{-j/a_0}u^*(2^{-j}t)
	\lesssim P_{a_0}u^*(t),
\end{aligned}
\]
and
\[
\begin{aligned}
	t^{-1/a_1}
	\left(\int_t^\infty u^*(r)^{a_1}\,dr\right)^{1/a_1}
	=
	\left(
	\frac1t\sum_{j\in\mathbb N}
	\int_{2^{j-1}t}^{2^jt}u^*(r)^{a_1}\,dr
	\right)^{1/a_1}
	\le
	\sum_{j\in\mathbb N}
	2^{(j-1)/a_1}u^*(2^{j-1}t)
	\lesssim u^*(t)+Q_{a_1}u^*(t).
\end{aligned}
\]
Since $P_{a_0}u^*(t)\ge a_0u^*(t)$,
\eqref{eq:calculus-endpoint-rearrangement} implies
\begin{equation}\label{eq:critical-calculus-rearrangement}
	(Tu)^*(2t)
	\lesssim P_{a_0}u^*(t)+Q_{a_1}u^*(t).
\end{equation}
Using Minkowski's inequality and
\eqref{eq:critical-calculus-dilation} or
\eqref{eq:rates-window-dilation}, with the logarithmic factor
bounded by $1+\log_+(1/r)$, we obtain
\begin{equation}\label{eq:calculus-Hardy-P-bound}
	\|P_{a_0}u^*\|_{\widehat Y}
	\lesssim
	\|u^*\|_{\widehat Y}
	\int_0^1
	r^{1/a_0-1/\sigma-1}\bigl(1+\log(1/r)\bigr)\,dr
	\lesssim\|u^*\|_{\widehat Y},
\end{equation}
and
\begin{equation}\label{eq:calculus-Hardy-Q-bound}
	\|Q_{a_1}u^*\|_{\widehat Y}
	\lesssim
	\|u^*\|_{\widehat Y}
	\int_1^\infty r^{1/a_1-1/\sigma-1}\,dr
	\lesssim\|u^*\|_{\widehat Y}.
\end{equation}
Both integrals converge because $a_0<\sigma<a_1$.
By \eqref{eq:calculus-profile-norm},
\eqref{eq:critical-calculus-rearrangement},
\eqref{eq:calculus-Hardy-P-bound},
\eqref{eq:calculus-Hardy-Q-bound}, and
\eqref{eq:critical-calculus-dilation} or
\eqref{eq:rates-window-dilation} at $c=1/2$,
\[
\|Tu\|_Y
\lesssim
\|D_{1/2}(P_{a_0}u^*+Q_{a_1}u^*)\|_{\widehat Y}
\lesssim\|u^*\|_{\widehat Y}
\simeq\|u\|_Y.
\]
Summing over $T\in\{\mathcal M,\mathcal R_1,\ldots,\mathcal R_{2n}\}$
gives
\begin{equation}\label{eq:critical-calculus-bound}
	\|\mathcal Mu\|_Y
	+\sum_{\ell=1}^{2n}\|\mathcal R_\ell u\|_Y
	\le C\|u\|_Y.
\end{equation}
This completes the proof.
\end{proof}

\begin{lemma}\label{lem:Lorentz-derivatives}
	Let $1<\sigma<\infty$, $0<s\le\infty$ and
	$u\in C^2(\C^n)\cap L^\infty(\C^n)$.
	If $|\bar\partial u|\in\mathcal L^{\sigma,s}$, then
	\begin{equation}\label{eq:Lorentz-derivatives}
		\|\partial u\|_{\mathcal L^{\sigma,s}}
		\le C\|\bar\partial u\|_{\mathcal L^{\sigma,s}}.
	\end{equation}
	The same assertion holds with $\mathcal L^{\sigma,s}$
	replaced by $\mathcal L^{\sigma,2,1/2}$ or
	$\mathcal L^{\sigma,\sigma}\cap\mathcal{CL}^{\sigma,2}$,
	and also by $\mathcal X^{\sigma,s}$ when
	$\sigma>2$ and $1\le s<\infty$.
	The constants are independent of $u$.
\end{lemma}

\begin{proof}
	Fix $1<\sigma<\infty$ and let $Y$ denote
	$\mathcal L^{\sigma,s}$, $\mathcal L^{\sigma,2,1/2}$,
	$\mathcal L^{\sigma,\sigma}\cap\mathcal{CL}^{\sigma,2}$,
	or $\mathcal X^{\sigma,s}$, with $0<s\le\infty$ in the first
	case and $\sigma>2$, $1\le s<\infty$ in the last case.
	For $0<s<\infty$, monotonicity gives
	\[
	t^{s/\sigma}h^*(t)^s
	\lesssim\int_{t/2}^t r^{s/\sigma-1}h^*(r)^s\,dr
	\le\|h\|_{\mathcal L^{\sigma,s}}^s.
	\]
	For $s=\infty$, the corresponding bound follows by definition.
	Together with \eqref{eq:critical-calculus-embedding}, this proves
	\begin{equation}\label{eq:elliptic-weak-embedding}
		\|h\|_{\mathcal L^{\sigma,\infty}}\lesssim\|h\|_Y.
	\end{equation}
	
	Write $z_j=x_{2j-1}+\mathrm i x_{2j}$ and define
	\begin{equation}\label{eq:elliptic-transform-definition}
		\mathcal T_{kj}
		:=-(\mathcal R_{2k-1}-\mathrm i\mathcal R_{2k})
		(\mathcal R_{2j-1}-\mathrm i\mathcal R_{2j}),
		\qquad k,j\in\{1,\ldots,n\}.
	\end{equation}
	For a Schwartz function $\eta$, the Fourier multipliers give
	\begin{equation}\label{eq:elliptic-multiplier}
		\widehat{\mathcal T_{kj}\eta}(\xi)
		=
		\frac{(\xi_{2k-1}-\mathrm i\xi_{2k})
			(\xi_{2j-1}-\mathrm i\xi_{2j})}{|\xi|^2}
		\widehat\eta(\xi),
		\qquad \xi\ne0.
	\end{equation}
	Since
	$\widehat{\partial_{z_j}\eta}(\xi)
	=\frac{\mathrm i}{2}(\xi_{2j-1}-\mathrm i\xi_{2j})\widehat\eta(\xi)$
	and $\widehat{\Delta\eta}(\xi)=-|\xi|^2\widehat\eta(\xi)$,
	\eqref{eq:elliptic-multiplier} yields
	\[
	\widehat{\Delta\mathcal T_{kj}\eta}(\xi)
	=-(\xi_{2k-1}-\mathrm i\xi_{2k})
	(\xi_{2j-1}-\mathrm i\xi_{2j})\widehat\eta(\xi)
	=\widehat{4\partial_{z_k}\partial_{z_j}\eta}(\xi).
	\]
Taking inverse Fourier transforms, we obtain
$\Delta\mathcal T_{kj}\eta
=4\partial_{z_k}\partial_{z_j}\eta$. Choose $1<a_0<\sigma<a_1<\infty$.
	By \eqref{eq:calculus-level-decomposition},
	$\mathcal L^{\sigma,\infty}\subset L^{a_0}+L^{a_1}$.
	The boundedness and consistency of the $L^{a_0}$ and $L^{a_1}$
	extensions of $\mathcal T_{kj}$, together with Schwartz density,
	therefore give
	\begin{equation}\label{eq:elliptic-operator-identity}
		\Delta\mathcal T_{kj}h
		=4\partial_{z_k}\partial_{z_j}h,
		\qquad h\in L^{a_0}+L^{a_1},
	\end{equation}
	in the sense of distributions.
	
	Put $\omega_j:=\partial_{\bar z_j}u$ and
	$V_k:=\sum_{j=1}^n\mathcal T_{kj}\omega_j$.
	Lemma~\ref{lem:known-Lorentz-calculus} and 
	\eqref{eq:elliptic-transform-definition} imply
	\begin{equation}\label{eq:elliptic-Lorentz-bound}
		\|V_k\|_Y
		\lesssim\sum_{j=1}^n\|\omega_j\|_Y
		\lesssim\|\bar\partial u\|_Y,
		\qquad k\in\{1,\ldots,n\}.
	\end{equation}
	By \eqref{eq:elliptic-weak-embedding},
	$\omega_j\in\mathcal L^{\sigma,\infty}\subset L^{a_0}+L^{a_1}$.
	Using \eqref{eq:elliptic-operator-identity} and
	$\Delta u=4\sum_{j=1}^n\partial_{z_j}\omega_j$, we obtain
	\begin{equation}\label{eq:elliptic-distribution}
		\Delta V_k
		=4\sum_{j=1}^n\partial_{z_k}\partial_{z_j}\omega_j
		=\Delta\partial_{z_k}u
	\end{equation}
	in the sense of distributions.
	Thus $w_k:=\partial_{z_k}u-V_k$ satisfies $\Delta w_k=0$,
	and Weyl's lemma gives a smooth harmonic representative of $w_k$.
	
	To show that $w_k=0$, first note that
	\eqref{eq:elliptic-weak-embedding} gives
	\begin{equation}\label{eq:Lorentz-local-integral}
		\int_E|h|\,dv
		\le\int_0^{v(E)}h^*(t)\,dt
		\le\|h\|_{\mathcal L^{\sigma,\infty}}
		\int_0^{v(E)}t^{-1/\sigma}\,dt
		\lesssim\|h\|_Yv(E)^{1-1/\sigma}
	\end{equation}
	for every $h\in Y$ and measurable set $E$ of finite measure.
	Fix $x_0\in\mathbb R^{2n}$ and choose a nonnegative radial
	$\chi\in C_c^\infty(B(0,1))$ with $\int_{B(0,1)}\chi\,dv=1$.
	For $A>0$, put
	$\chi_{A,x_0}(x):=A^{-2n}\chi((x-x_0)/A)$.
	The mean-value property and integration by parts on
	$\mathbb R^{2n}$ give
	\begin{equation}\label{eq:elliptic-mean-value}
		w_k(x_0)
		=\int_{\mathbb R^{2n}}w_k\chi_{A,x_0}\,dv
		=-\int_{\mathbb R^{2n}}
		u\,\partial_{z_k}\chi_{A,x_0}\,dv
		-\int_{\mathbb R^{2n}}V_k\chi_{A,x_0}\,dv.
	\end{equation}
	Since $\chi$ is smooth and compactly supported,
	$\operatorname{supp}\chi_{A,x_0}\subset B(x_0,A)$, and scaling gives
	\[
	\|\partial_{z_k}\chi_{A,x_0}\|_{L^1}
	=A^{-1}\|\partial_{z_k}\chi\|_{L^1}\lesssim A^{-1},
	\qquad
	\|\chi_{A,x_0}\|_\infty
	=A^{-2n}\|\chi\|_\infty\lesssim A^{-2n}.
	\]
	Thus \eqref{eq:elliptic-mean-value} and
	\eqref{eq:Lorentz-local-integral} imply
\begin{align}\label{eq:eliminate-harmonic-part}
			|w_k(x_0)|
			&\lesssim
			A^{-1}\|u\|_\infty
			+A^{-2n}\int_{B(x_0,A)}|V_k|\,dv\notag\\
			&\lesssim
			A^{-1}\|u\|_\infty
			+A^{-2n/\sigma}\|V_k\|_Y
			\longrightarrow0
			\qquad(A\to\infty).
		\end{align}
	Since $x_0$ is arbitrary, $w_k=0$, and hence
	\begin{equation}\label{eq:elliptic-reconstruction}
		\partial_{z_k}u
		=V_k
		=\sum_{j=1}^n\mathcal T_{kj}\partial_{\bar z_j}u,
		\qquad k\in\{1,\ldots,n\},
	\end{equation}
	in the sense of distributions and almost everywhere.
	Finally, \eqref{eq:elliptic-reconstruction},
	\eqref{eq:elliptic-Lorentz-bound}, and the finite-sum
	quasi-triangle inequality yield
	\[
	\|\partial u\|_Y
	\lesssim\sum_{k=1}^n\|V_k\|_Y
	\lesssim\|\bar\partial u\|_Y.
	\]
	This proves \eqref{eq:Lorentz-derivatives} and its stated variants.
	This completes the proof.
\end{proof}

\subsection{The Berger--Coburn equivalence}
\begin{proof}[Proof of Theorem~\ref{thm:BC-approx}]
	Recall from \eqref{eq:BC-range} that
	\[
	\begin{aligned}
		\mathcal B
		={}&\left\{(p,q,s):0<p\le q<\infty,\ q\ge1,\
		\frac1{1+\eta_{p,q}}<s<\infty\right\}\\
		&{}\cup
		\left\{(p,q,s):1\le q<p<\infty,\
		\frac r{r-1}<s<\infty\right\},
	\end{aligned}
	\]
	where $r:=pq/(p-q)$ when $p>q$.
	Put $\sigma:=\rho_{p,q,s}$ when $(p,q,s)\notin\mathcal E$
	and $\sigma:=\tau_{p,q}$ when $(p,q,s)\in\mathcal E$.
	For \textup{(i)}, take
	\[
	Y:=
	\begin{cases}
		\mathcal L^{\sigma,s},
		&(p,q,s)\notin\mathcal E,\\
		\mathcal L^{\sigma,2,1/2},
		&(p,q,s)\in\mathcal E,\ p^\dagger=q,\\
		\mathcal L^{\sigma,\sigma}\cap\mathcal{CL}^{\sigma,2},
		&(p,q,s)\in\mathcal E,\ p^\dagger\ne q.
	\end{cases}
	\]
	For \textup{(ii)}, take
	\[
	Y:=
	\begin{cases}
		\mathcal L^{\sigma,\infty},
		&(p,q,s)\notin\mathcal E,\\
		\mathcal X^{\sigma,\varsigma},
		&(p,q,s)\in\mathcal E.
	\end{cases}
	\]
	By \eqref{eq:lorentz-index}, the lower bounds on $s$ ensure
	$\sigma>1$, and $\sigma=\tau_{p,q}>2$ when
	$(p,q,s)\in\mathcal E$.
	These are the exponent conditions required in
	Lemmas~\ref{lem:known-Lorentz-calculus}
	and~\ref{lem:Lorentz-derivatives}.
	By Theorems~\ref{thm:increasing}
	and~\ref{cor:weak-approximation-ideals}, respectively,
	it suffices to prove
	\begin{equation}\label{eq:BC-Lorentz-transfer}
		\|G_{q,1}(\bar f)\|_Y
		\lesssim\|G_{q,1}(f)\|_Y
		\qquad\text{whenever }G_{q,1}(f)\in Y.
	\end{equation}
	
	The radius comparisons \eqref{eq:rearranged-radius},
	\eqref{eq:logarithmic-radius}, and
	\eqref{eq:critical-Copson-radius} apply to the first three spaces.
	For $\mathcal X^{\sigma,\varsigma}$, use
	\eqref{eq:radius-translation-comparison},
	\eqref{eq:rearrangement-sum}, and
	\eqref{eq:rates-window-dilation}.
	Consequently, for each fixed $r>0$,
	\begin{equation}\label{eq:BC-radius-comparison}
		\|G_{q,r}(h)\|_Y\simeq\|G_{q,1}(h)\|_Y,
		\qquad h\in\{f,\bar f\}.
	\end{equation}
	
	Put $\varepsilon:=1/16$ and $A:=G_{q,2}(f)$.
	By Lemma~\ref{lem:known-decomposition}, write $f=f_1+f_2$ with
	$f_1=\sum_{j\in\mathbb N}\psi_jh_j$, where
	$h_j\in\mathcal H(B(a_j,1))$ and
	$\operatorname{supp}\psi_j\subset B(a_j,1/2)$.
	The partition has uniformly finite overlap and satisfies
	\begin{equation}\label{eq:BC-partition-bounds}
		\psi_j\ge0,\qquad
		\sum_{j\in\mathbb N}\psi_j=1,\qquad
		\sup_w\sum_{j\in\mathbb N}|\partial\psi_j(w)|\le C.
	\end{equation}
	By \eqref{eq:decomposition} and \cite[(5-4)]{HV23},
	$f_1\in C^2\cap L^\infty$ and
	$|\bar\partial f_1|+M_{q,\varepsilon}(f_2)\lesssim A$.
	Together with \eqref{eq:BC-radius-comparison}, this gives
	\begin{equation}\label{eq:BC-Lorentz-decomposition}
		\|\bar\partial f_1\|_Y
		+\|M_{q,\varepsilon}(f_2)\|_Y+\|A\|_Y
		\lesssim\|G_{q,1}(f)\|_Y.
	\end{equation}
	Lemma~\ref{lem:Lorentz-derivatives} and
	\eqref{eq:BC-Lorentz-decomposition} therefore imply
	\begin{equation}\label{eq:BC-Lorentz-gradient}
		\|\partial f_1\|_Y
		\lesssim\|G_{q,1}(f)\|_Y.
	\end{equation}
	Fix $z\in\C^n$ and choose $k$ with $\psi_k(z)>0$.
	Put $J_z:=\{j\in\mathbb N:
	\operatorname{supp}\psi_j\cap B(z,2\varepsilon)
	\ne\varnothing\}$. For $j\in J_z$, we have
	$|a_j-z|<1/2+2\varepsilon=5/8$ and $|a_k-z|<1/2$.
	Hence
	\[
	B(z,4\varepsilon)\subset B(a_j,1)\cap B(a_k,1),
	\qquad
	B(a_j,1)\cup B(a_k,1)\subset B(z,2).
	\]
	By the choice of the local approximants in
	\cite[proof of Lemma~3.6]{HV23} and holomorphic interior estimates,
\begin{align}\label{eq:BC-holomorphic-comparison11}
			\sup_{B(z,2\varepsilon)}
			\bigl(|h_j-h_k|+|\partial(h_j-h_k)|\bigr)
			&\lesssim M_{q,4\varepsilon}(h_j-h_k)(z)\notag\\
			&\lesssim M_{q,1}(f-h_j)(a_j)+M_{q,1}(f-h_k)(a_k)\notag\\
			&\lesssim G_{q,1}(f)(a_j)+G_{q,1}(f)(a_k)
			\lesssim A(z).
		\end{align}
	On $B(z,2\varepsilon)$, \eqref{eq:BC-partition-bounds} gives
	\[
	\partial f_1-\partial h_k
	=\sum_{j\in J_z}
	\bigl(\psi_j\partial(h_j-h_k)
	+\partial\psi_j(h_j-h_k)\bigr).
	\]
	Thus \eqref{eq:BC-partition-bounds} and
	\eqref{eq:BC-holomorphic-comparison11} yield
	\begin{equation}\label{eq:BC-local-error}
		\sup_{B(z,2\varepsilon)}
		|\partial f_1-\partial h_k|\lesssim A(z).
	\end{equation}
	The submean inequality for $\partial h_k$ and
	\eqref{eq:BC-local-error} imply
	\[
	\sup_{B(z,\varepsilon)}|\partial h_k|
	\lesssim M_{1,2\varepsilon}(\partial h_k)(z)
	\lesssim\mathcal M(|\partial f_1|)(z)+A(z).
	\]
	Applying \eqref{eq:BC-local-error} once more gives
\begin{equation}\label{eq:local-transfer}
	M_{q,\varepsilon}(\partial f_1)(z)
	\le
	\sup_{B(z,\varepsilon)}|\partial h_k|
	+\sup_{B(z,\varepsilon)}|\partial f_1-\partial h_k|
	\lesssim\mathcal M(|\partial f_1|)(z)+A(z).
\end{equation}
	By Lemma~\ref{lem:known-Lorentz-calculus},
	\eqref{eq:BC-Lorentz-decomposition}, and
	\eqref{eq:BC-Lorentz-gradient}, we conclude that
	\begin{equation}\label{eq:BC-Lorentz-local-mean}
		\|M_{q,\varepsilon}(\partial f_1)\|_Y
		\lesssim\|G_{q,1}(f)\|_Y.
	\end{equation}
	
	Finally, $\bar\partial\bar f_1=\overline{\partial f_1}$ and
	\eqref{eq:radius-local-distance} give
	\begin{equation}\label{eq:BC-Lorentz-Henkin}
		G_{q,\varepsilon/2}(\bar f)
		\lesssim
		M_{q,\varepsilon}(\partial f_1)
		+M_{q,\varepsilon}(f_2).
	\end{equation}
	Combining \eqref{eq:BC-Lorentz-Henkin},
	\eqref{eq:BC-Lorentz-decomposition},
	\eqref{eq:BC-Lorentz-local-mean}, and
	\eqref{eq:BC-radius-comparison} proves
	\eqref{eq:BC-Lorentz-transfer}.
	Interchanging $f$ and $\bar f$ and applying the corresponding
	membership and norm criteria proves both assertions.
	This completes the proof.
\end{proof}

\subsection{Sharpness of the exponent ranges}

The following examples show that the lower bounds
on $s$ in \eqref{eq:BC-range} are sharp.

\begin{theorem}\label{thm:BC-counter}
	Let $n\ge1$, $0<p<\infty$, $1\le q<\infty$, and $0<s<\infty$.
	Suppose that $\operatorname{Hess}_{\mathbb R}\varphi\simeq E$
	and $(p,q,s)\notin\mathcal B$.
	Then there exist $f,g\in C^\infty(\C^n)\cap L^\infty(\C^n)$ such that
	\[
	H_f\in\mathfrak A_s(\F p,\Lw q),
	\qquad
	H_{\bar f}\notin\mathfrak A_s(\F p,\Lw q),
	\]
	and
	\[
	H_g\in\mathfrak A_{s,\infty}(\F p,\Lw q),
	\qquad
	H_{\bar g}\notin\mathfrak A_{s,\infty}(\F p,\Lw q).
	\]
\end{theorem}

\begin{proof}
	By \eqref{eq:BC-range}, our assumption gives
	$s\le(1+\eta)^{-1}$ when $p\le q$ and
	$s\le r/(r-1)$ when $p>q$.
	Hence \eqref{eq:lorentz-index} yields $0<\rho\le1$.
	For each $Y\in\{\mathcal L^{\rho,s},\mathcal L^{\rho,\infty}\}$,
	we construct $u_j\in C_c^\infty(\C^n)$ such that
	\begin{equation}\label{eq:BC-block-normalization}
		\|u_j\|_\infty+\|G_{q,1}(u_j)\|_Y\le2^{-j^2},
		\qquad
		\|G_{q,1}(\bar u_j)\|_Y\ge j.
	\end{equation}
	
	Write $z=(z_1,z')\in\C\times\C^{n-1}$,
	omitting transverse factors when $n=1$.
	Fix $\theta\in C^\infty(\mathbb R)$ such that
	$0\le\theta\le1$, $\theta(t)=1$ for $t\le0$,
	and $\theta(t)=0$ for $t\ge1$.
	The Henkin estimate used in \eqref{eq:radius-local-distance}, applied
directly to a smooth function $u$, gives
\begin{equation}\label{eq:BC-smooth-local-distance}
G_{q,1}(u)(z)\lesssim M_{q,2}(\bar\partial u)(z)
\lesssim\sup_{B(z,2)}|\bar\partial u|.
\end{equation}
	For $R>4$, define
	\[
	u_R(z):=
	\begin{cases}
		\displaystyle
		\frac{\theta(2-|z_1|)}{z_1}
		\theta(|z_1|-R)\theta(|z'|-R),
		& |z_1|>1,\\[2pt]
		0,& |z_1|\le1.
	\end{cases}
	\]
	Then $u_R\in C_c^\infty(\C^n)$ and
	$\|u_R\|_\infty\le1$.
	All implicit constants below are independent of $R$ and $j$.
	The three cutoff derivatives have size $O(1)$, $O(R^{-1})$,
and $O((1+|z_1|)^{-1})$, respectively, on their transition regions.
By \eqref{eq:BC-smooth-local-distance},
	\[
	G_{q,1}(u_R)(z)\lesssim
	\bigl(\chi_{\{|z_1|<4\}}
	+R^{-1}\chi_{\{||z_1|-R|<3\}}\bigr)
	\chi_{\{|z'|<R+3\}}
	+\frac{\chi_{\{|z_1|<R+3,\;||z'|-R|<3\}}}{1+|z_1|},
	\]
	with the last term omitted when $n=1$.
For $n\ge2$, denote that term by $A_R$. Its transverse shell has
thickness $O(1)$, so
\[
v\{A_R>\lambda\}\lesssim R^{2n-3}\min\{R^2,\lambda^{-2}\},
\qquad 0<\lambda<1,
\]
and $A_R\le1$. Thus $\|A_R\|_Y\lesssim R^{(2n-1)/\rho-1}$.
The first two terms are bounded functions on sets of volume
$O(R^{2n-2})$ and $O(R^{2n-1})$, respectively. Consequently,
	\begin{equation}\label{eq:BC-cutoff-upper}
		\|u_R\|_\infty+\|G_{q,1}(u_R)\|_Y
		\lesssim R^{(2n-2)/\rho}+R^{(2n-1)/\rho-1}
		\lesssim R^{(2n-1)/\rho-1}.
	\end{equation}
	
	For $3<|z_1|<R/2$ and $|z'|<R/2$, uniformly for
	$w\in B(0,1)$,
	\[
	\overline{u_R(z+w)}
	=\frac1{\bar z_1+\bar w_1}
	=\bar z_1^{-1}-\bar z_1^{-2}\bar w_1
	+O(|z_1|^{-3}).
	\]
	Since the holomorphic subspace is closed in $L^q(B(0,1))$,
$d_q:=G_{q,1}(\bar z_1)(0)>0$. Removing the constant term gives
\[
G_{q,1}(\bar u_R)(z)\ge d_q|z_1|^{-2}-C|z_1|^{-3}
\gtrsim|z_1|^{-2}
\]
for $r_q<|z_1|<R/2$, where $r_q>3$ is fixed.
	Consequently,
	\begin{equation}\label{eq:BC-cutoff-rearrangement}
		(G_{q,1}(\bar u_R))^*(t)
		\gtrsim\frac{R^{2n-2}}{t},
		\qquad
		CR^{2n-2}<t<cR^{2n}.
	\end{equation}
	For $0<\rho<1$, the lower bound on
	$\{R/4<|z_1|<R/2,\ |z'|<R/2\}$ gives
	\[
	\|G_{q,1}(\bar u_R)\|_{\mathcal L^{\rho,\infty}}
	\gtrsim R^{-2}
	\left(R^{2n-2}\int_{R/4}^{R/2}t\,dt\right)^{1/\rho}
	\simeq R^{2n/\rho-2}.
	\]
	For $\rho=1$, \eqref{eq:BC-cutoff-rearrangement} yields
	\[
	\|G_{q,1}(\bar u_R)\|_{\mathcal L^{1,s}}^s
	\gtrsim
	R^{(2n-2)s}
	\int_{CR^{2n-2}}^{cR^{2n}}\frac{dt}{t}
	\simeq R^{(2n-2)s}\log R.
	\]
	Together with \eqref{eq:BC-cutoff-upper} and
	$\mathcal L^{\rho,s}\subset\mathcal L^{\rho,\infty}$,
	these estimates imply, as $R\to\infty$,
	\begin{equation}\label{eq:BC-cutoff-ratio}
		\frac{\|G_{q,1}(\bar u_R)\|_Y}
		{\|u_R\|_\infty+\|G_{q,1}(u_R)\|_Y}
		\gtrsim
		\begin{cases}
			R^{1/\rho-1}\longrightarrow\infty,
			&0<\rho<1,\\[2pt]
			(\log R)^{1/s}\longrightarrow\infty,
			&\rho=1,\quad Y=\mathcal L^{1,s}.
		\end{cases}
	\end{equation}
	
	It remains to consider $Y=\mathcal L^{1,\infty}$.
Choose an even $\psi\in C_c^\infty(-1/4,1/4)$ such that
$\psi\ge0$ and $\int_{\mathbb R}\psi(u)\,du=1$.
Define
\[
b_R(t):=
\int_{-1/4}^{1/4}
\max\{\log R-|t-u-\log R|,0\}\,\psi(u)\,du,
\qquad t\in\mathbb R.
\]
	Then $\operatorname{supp}b_R
	\subset[-1/4,2\log R+1/4]$, $|b_R'|\le1$,
	and $b_R(t)=t$ for $1/4<t<\log R-1/4$.
	Moreover,
	\[
	\sup_t e^{-t}b_R(t)\lesssim1,
	\qquad
	\sup_t e^tb_R(t)\lesssim R^2.
	\]
	Define $v_R(0,z'):=0$ and, for $z_1\ne0$,
	\[
	v_R(z):=
	R^{4-4n}\frac{b_R(\log|z_1|)}{z_1}
	\theta\!\left(\frac{|z'|}{R^2}-1\right).
	\]
	Then $v_R\in C_c^\infty(\C^n)$ and
	$\|v_R\|_\infty\lesssim R^{4-4n}\le1$.
	Writing $r=|z_1|>0$, differentiation gives
\[
\bar\partial_{z_1}v_R
=R^{4-4n}\frac{b_R'(\log r)}{2r^2}
\theta\!\left(\frac{|z'|}{R^2}-1\right).
\]
The transverse derivative is bounded by
$CR^{2-4n}b_R(\log r)/r\le CR^{4-4n}r^{-2}$,
because $rb_R(\log r)\lesssim R^2$ and this term is omitted for $n=1$.
All derivatives vanish unless $e^{-1/4}\le r\le e^{1/4}R^2$
and $|z'|\le2R^2$. Hence \eqref{eq:BC-smooth-local-distance} gives
	\[
	G_{q,1}(v_R)(z)
	\lesssim
	\frac{R^{4-4n}}{(1+|z_1|)^2}
	\chi_{\{|z_1|<3R^2,\;|z'|<3R^2\}}(z).
	\]
	Thus, for $\lambda>0$,
	\[
	v\bigl(\{G_{q,1}(v_R)>\lambda\}\bigr)
	\lesssim R^{4n-4}\frac{R^{4-4n}}{\lambda}
	=\lambda^{-1},
	\]
	and hence
	\begin{equation}\label{eq:BC-weak-endpoint-upper}
		\|v_R\|_\infty+
		\|G_{q,1}(v_R)\|_{\mathcal L^{1,\infty}}
		\lesssim1.
	\end{equation}
	
	For sufficiently large $R$, if
	$R/e^2<|z_1|<R/e$ and $|z'|<R^2/2$, then
	\[
	v_R(z+w)
	=R^{4-4n}\frac{\log|z_1+w_1|}{z_1+w_1},
	\qquad w\in B(0,1).
	\]
	Consequently,
	\[
	\partial_{z_1}v_R(z)
	=R^{4-4n}\frac{1/2-\log|z_1|}{z_1^2},
	\qquad
	\sup_{\zeta\in B(z,1)}
	|\operatorname{Hess}_{\mathbb R}v_R(\zeta)|
	\lesssim R^{1-4n}\log R.
	\]
	Taylor's formula gives, uniformly for $w\in B(0,1)$,
	\[
	\overline{v_R(z+w)}
	=
	\overline{v_R(z)}
	+\partial_{z_1}\bar v_R(z)\,w_1
	+\overline{\partial_{z_1}v_R(z)}\,\bar w_1
	+O(R^{1-4n}\log R).
	\]
	Removing the first two terms and using the same $d_q>0$ gives
\[
G_{q,1}(\bar v_R)(z)
\ge d_q|\partial_{z_1}v_R(z)|-CR^{1-4n}\log R
\gtrsim R^{2-4n}\log R
\]
for sufficiently large $R$.
	Combining this with \eqref{eq:BC-weak-endpoint-upper}, we obtain
	\begin{equation}\label{eq:BC-weak-endpoint-ratio}
		\frac{\|G_{q,1}(\bar v_R)\|_{\mathcal L^{1,\infty}}}
		{\|v_R\|_\infty+\|G_{q,1}(v_R)\|_{\mathcal L^{1,\infty}}}
		\gtrsim
		R^{2-4n}\log R\,
		R^{4n-4}\int_{R/e^2}^{R/e}t\,dt
		\simeq\log R\longrightarrow\infty.
	\end{equation}
	
	By \eqref{eq:BC-cutoff-ratio} and
	\eqref{eq:BC-weak-endpoint-ratio},
	choose $w_j\in\{u_R:R>4\}$ when
	$Y\ne\mathcal L^{1,\infty}$ and
	$w_j\in\{v_R:R>4\}$ otherwise, such that
	\[
	\frac{\|G_{q,1}(\bar w_j)\|_Y}
	{\|w_j\|_\infty+\|G_{q,1}(w_j)\|_Y}
	\ge j2^{j^2}.
	\]
	Define
	\[
	u_j:=
	\frac{2^{-j^2}w_j}
	{\|w_j\|_\infty+\|G_{q,1}(w_j)\|_Y}.
	\]
	Then \eqref{eq:BC-block-normalization} holds.
	As in \cite[Sections~6--7]{Xia26}, choose $x_j\in\C^n$
	such that
	\[
	\operatorname{dist}
	(x_j+\operatorname{supp}u_j,\,
	x_k+\operatorname{supp}u_k)>4
	\quad(j\ne k),
	\qquad
	F(z):=\sum_{j\ge1}u_j(z-x_j).
	\]
	The sum is locally finite, so $F\in C^\infty(\C^n)$ and
	$\|F\|_\infty\le\sup_j\|u_j\|_\infty\le1$.
	Each unit ball meets at most one translated support and hence
	\[
	G_{q,1}(F)(z)
	=\sum_{j\ge1}G_{q,1}(u_j)(z-x_j),
	\qquad
	G_{q,1}(\bar F)(z)
	=\sum_{j\ge1}G_{q,1}(\bar u_j)(z-x_j).
	\]
	By translation invariance,
	\eqref{eq:rearrangement-sum-norm},
	\eqref{eq:rearrangement-sum},
	\eqref{eq:BC-block-normalization}, and Fatou's lemma,
	\[
	\|G_{q,1}(F)\|_Y
	\lesssim
	\sum_{j\ge1}C^j\|G_{q,1}(u_j)\|_Y
	\le\sum_{j\ge1}C^j2^{-j^2}<\infty,
	\]
	whereas
	\[
	\|G_{q,1}(\bar F)\|_Y
	\ge\sup_{j\ge1}\|G_{q,1}(\bar u_j)\|_Y
	\ge\sup_{j\ge1}j=\infty.
	\]
	
	For $Y=\mathcal L^{\rho,s}$ and
	$Y=\mathcal L^{\rho,\infty}$, denote the resulting symbols
	by $f$ and $g$, respectively.
	Since $f,g\in C^\infty(\C^n)\cap L^\infty(\C^n)
	\subset\mathcal S$,
	Theorems~\ref{thm:increasing}\textup{(i)}
	and~\ref{cor:weak-approximation-ideals}\textup{(i)}
	prove the assertions.
\end{proof}

\section{Further remarks}
The proofs in this paper combine local holomorphic approximation,
weighted $L^p$ estimates for the $\bar\partial$-equation,
and a comparison with diagonal operators.
We believe that, with suitable modifications, these methods
can be used to establish analogous $(p,q)$ approximation-ideal
criteria and prescribed-decay results in other settings.
These include Bergman spaces on the unit disk
\cite{Luecking92} and on bounded strongly pseudoconvex domains
in $\C^n$ \cite{Li94}, Bergman spaces with exponential weights
on the unit disk \cite{HLS19}, and the weighted Fock spaces
introduced in \cite{LWZ26}.

\end{document}